\documentclass[11pt]{article}
\input amssymb.sty
\input amssym.def
\input amssym
\input epsf

\newtheorem{theorem}{Theorem}[section]
\newtheorem{lemma}[theorem]{Lemma}
\newtheorem{proposition}[theorem]{Proposition}
\newtheorem{corollary}[theorem]{Corollary}
\newtheorem{conjecture}[theorem]{Conjecture}
\newtheorem{definition}[theorem]{Definition}
\newtheorem{theoremconstruction}[theorem]{Theorem-Construction}

\begin{document}
\newcommand{\be}{\begin{equation}}
\newcommand{\ee}{\end{equation}}
\newcommand{\bt}{\begin{theorem}}
\newcommand{\et}{\end{theorem}}
\newcommand{\bd}{\begin{definition}}
\newcommand{\ed}{\end{definition}}
\newcommand{\bp}{\begin{proposition}}
\newcommand{\ep}{\end{proposition}}
\newcommand{\bl}{\begin{lemma}}
\newcommand{\el}{\end{lemma}}
\newcommand{\bc}{\begin{corollary}}
\newcommand{\ec}{\end{corollary}}
\newcommand{\bcon}{\begin{conjecture}}
\newcommand{\econ}{\end{conjecture}}
\newcommand{\la}{\label}
\newcommand{\Z}{{\Bbb Z}}
\newcommand{\R}{{\Bbb R}}
\newcommand{\Q}{{\Bbb Q}}
\newcommand{\C}{{\Bbb C}}
\newcommand{\hra}{\hookrightarrow}
\newcommand{\lra}{\longrightarrow}
\newcommand{\lms}{\longmapsto}

\begin{titlepage}
\title{Hodge correlators II}
\author{A.B. Goncharov }
\date{\it To Pierre Deligne for his 65th birthday}

\end{titlepage}
\maketitle

\tableofcontents

\section{Introduction} 

\subsection{Summary} Let $X$ be a compact Kahler manifold of dimension $n$. 
Denote by ${\cal C}{\cal C}_\bullet(X)$ 
the cyclic homology complex of the reduced cohomology algebra of $X$. 
We define a linear map, the {\it Hodge correlator map}: 
\be \la{1}
{\rm Cor}_{\cal H}: H_0\Bigl({\cal C}{\cal C}_\bullet(X)\otimes H_{2n}(X)[-2]\Bigr) \lra \C.
\ee
Evaluating this map on cycles  
we get complex numbers, called 
{\it Hodge correlators}.

We introduce 
a Feynman integral related to $X$. We show that 
its correlators, 
defined via the standard perturbative series 
expansion procedure, are well defined, and identify them with the Hodge correlators. 
We show that the Hodge correlator map defines a functorial real mixed Hodge structure of the 
rational homotopy type of $X$.   The Hodge correlators are the homotopy periods of $X$.

\vskip 2mm
The category of real mixed Hodge structures (MHS) is canonically identified 
with the category of representations of the Hodge Galois group. 
We show that the Hodge correlator map provides  
an explicit  action 
of the Hodge 
Galois group on the formal neighborhood of the trivial local system on $X$. 
It
encodes a functorial homotopy action 
of the Hodge 
Galois group on the rational homotopy type of $X$, 
describing the real MHS on the latter.  

A real MHS on the rational homotopy type of a  
complex algebraic variety was defined, by different methods, 
 by J. Morgan \cite{M} and R. Hain \cite{H}. Our construction should give
 the same real MHS. We prove this for the 
real MHS on the 
pronilpotent completion of the fundamental group of $X$. 

\vskip 3mm
Now let $X$ be a regular projective algebraic 
variety over a field $k$. Assuming the motivic formalism,  we define 
{\it motivic correlators of $X$}. They lie in the motivic Lie coalgebra of $k$. 
Given an embedding $k \hra \C$, the periods 
of the Hodge realization of the motivic correlators are the Hodge correlators. 
The motivic correlators come together with an explicit formula for their coproduct 
in the motivic Lie coalgebra of $k$. 
This is one of the main advantages of presenting homotopy periods 
of $X(\C)$ as periods of the motivic correlators. Using this one can perform 
{arithmetic analysis} 
of the homotopy periods of varieties defined over number fields. 
Summarising: 

\begin{figure}[ht]
\centerline{\epsfbox{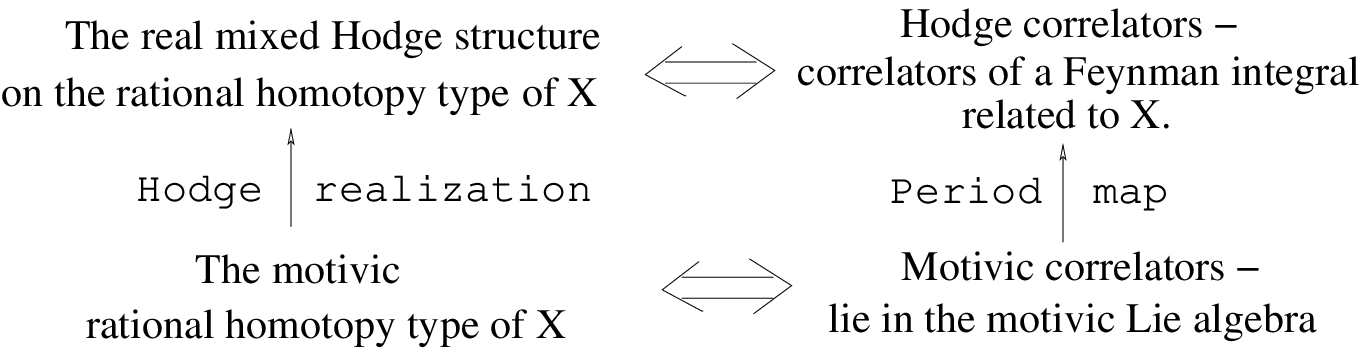}}
\label{hc1}
\end{figure}
A similar  picture for a smooth 
but not necessarily compact 
complex curve $X$ was obtained in \cite{G1}. 
In that case 
there are natural harmonic representatives for ${\rm gr}^WH^{\ast}(X)$. 
This allows to avoid discussion of DG objects, and makes the story simpler. 

The results of this paper admit a generalization where 
$X$ is replaced by an arbitrary regular complex algebraic variety, 
which we hope to  discuss elsewhere.  

\vskip 2mm
In the rest of the Introduction we give a rather detailed account 
of our constructions and results. Complete details are available in
 the main body of the paper. 

\subsection{Hodge correlators for a compact Kahler manifold}  
Given a graded vector space $V$, denote by ${\cal C}_{V}$ the cyclic tensor envelope of 
$V$:
$$
{\cal C}_{V}:= \oplus_{m=0}^{\infty}\Bigl(\otimes^mV\Bigr)_{\Z/m\Z}
$$
where the subscript $\Z/m\Z$ denotes the coinvariants of the cyclic shift 
$$
v_1 \otimes v_2 \otimes ... \otimes v_m \lms 
(-1)^s v_m \otimes v_1 \otimes ... \otimes v_{m-1}, \qquad s = {\rm deg}(v_m)
\sum_{i=1}^{m-1}{\rm deg}(v_i). 
$$
We denote by $(v_1 \otimes v_2 \otimes ... \otimes v_m)_{\cal C}$ 
the projection of the element $v_1 \otimes v_2 \otimes ... \otimes v_m$ to ${\cal C}_{V}$. 
Let
$$
\overline {\rm H}_*(X):= \oplus_{i=1}^{2n-1}{\rm H}_i(X),
$$
be the reduced cohomology algebra of $X$. Let us shift it to the right by one:
$$
{\Bbb H}_*:= \overline {\rm H}_*(X)[-1]. 
$$
It sits 
in the degrees $[-(2n-1), -1]$.  
We also need the dual object:
  $$
{\Bbb H}^*:= \overline {\rm H}^*(X)[1], \qquad 
\overline {\rm H}^*(X):= \oplus_{i=1}^{2n-1}{\rm H}^i(X).
$$
The spaces ${\cal C}_{{\Bbb H}^*}$ and ${\cal C}_{{\Bbb H}_*}$ are graded dual 
to each other. 

\noindent
Finally, let us introduce the shifted by two fundamental cohomology class of $X$:
\be \la{sfc}
 {\cal H}:= H^{2n}(X)[2]. 
\ee
 Let  $ {\cal H}^{\vee}$ is its dual.

\vskip 3mm 
Choose a base point $a \in X$. 
Choose a splitting of the de Rham complex of $X$, bigraded in the usual way, into 
an arbitrary subspace ${\cal H}ar_X$ isomorphically projecting onto 
the cohomology of $X$ ("harmonic forms'') 
and its orthogonal complement. 
We take the $\delta$-function $\delta_a$ 
at the point $a \in X$ 
as a representative of the fundamental class.

Our first goal is to define a degree zero linear map, 
a precursor of the Hodge correlator map: 
\be \la{1.17.08.10q}
{\rm Cor}^*_{{\cal H}, a}: {\cal C}_{{\Bbb H}^*} \otimes 
{\cal H}^{\vee} \lra \C.
\ee

\vskip 2mm
We need a Green current $G_a(x,y)$ on $X \times X$. 
It satisfies the differential equation 
\be \la{GEQ}
(2\pi i)^{-1}\overline \partial 
\partial G_a(x,y) = \delta_{\Delta} - P_{\rm Har}
\ee
where $\delta_{\Delta}$ is the $\delta$-function of the diagonal, and 
$P_{\rm Har}$ is the Schwarz kernel of the projector onto the space of harmonic forms, 
realized by an $(n,n)$-form on $X \times X$. Given a basis $\{\alpha_i\}$ 
in the space of ``harmonic forms'' of dimensions $\not = 0, 2n$, and the dual basis 
$\{\alpha^{\vee}_i\}$, we have 
$$
P_{\rm Har} = \delta_a \otimes 1 + 1 \otimes \delta_a + \sum \alpha^{\vee}_i \otimes 
\alpha_i, \qquad \int_X  \alpha_i \wedge 
\alpha_j^{\vee} = \delta_{ij}.
$$
 The two currents on the right hand side of (\ref{GEQ}) represent 
the same cohomology class, so the equation has a solution by the 
$\overline \partial \partial$-lemma. 
Let us choose such a solution. 

\vskip 2mm
\begin{figure}[ht]
\centerline{\epsfbox{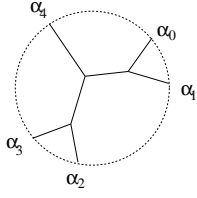}}
\caption{A plane  trivalent tree decorated by harmonic forms $\alpha_i$.}
\label{hf5}
\end{figure}
Let us pick a cyclic tensor product 
of harmonic forms
$$
W= (\alpha_0 \otimes \ldots \otimes \alpha_m)_{\cal C}.  
$$
 Take a plane trivalent tree 
$T$ decorated by $W$, see Fig \ref{hf5}. Let us  cook up 
a  top degree current 
on 
\be \la{4.30.08.1}
X ^{\{\mbox{\rm internal vertices of $T$}\}}.
\ee
 Let ${\cal A}^*_M[-1]$ be the de Rham complex of $M$, 
shifted by $1$ to the right. 
There is a degree zero linear map 
\be \la{xixi}
\xi: S^{*}({\cal A}^*_M[-1]) \to {\cal A}^*_M[-1], \qquad \varphi_0 \cdot \ldots \cdot
\varphi_m 
\lms {\rm Sym}_{m+1}\Bigl(\varphi_0 \wedge d^\C \varphi_1 \wedge \ldots  \wedge d^\C 
\varphi_m\Bigr). 
\ee 
Take a decorating harmonic form $\alpha_i$ assigned to 
an external edge $E_i$ of the tree $T$. Put the form $\alpha_i$ 
 to the copy of $X$ assigned  
to the internal vertex of the edge $E_i$. 
Apply the operator $\xi$  to the wedge product of 
the Green currents assigned to  
the edges of $T$. Finally,  multiply 
the obtained currents on (\ref{4.30.08.1}) 
with an appropriate sign, discussed in Section 2.3. One shows that we get a current. 
Integrating it over (\ref{4.30.08.1})  we 
get a number assigned to $T$. Taking the sum over all trees $T$, we get 
a complex number ${\rm Cor}^*_{{\cal H}, a}(W)$. Altogether, we 
get the map (\ref{1.17.08.10q}). One checks that its degree is zero.

\vskip 3mm
The $\cup$-product in ${\rm H}^{*}(X)$, followed by the projection to $\overline {\rm H}^*(X)$, 
induces an algebra structure on $\overline {\rm H}^*(X)$. 
There is a differential $\delta$ on ${\cal C}_{{\Bbb  H}^*}$,  
provided by products of neighbors in a cyclic word:
$$
\delta (\overline \alpha_0 \otimes \ldots \otimes \overline \alpha_m)_{\cal C} = 
{\rm Cycle}_{m+1}
(-1)^{|\alpha_m|}(\overline \alpha_0 \otimes \ldots \otimes \overline \alpha_{m-2} 
\otimes \overline {\alpha_{m-1}\cup \alpha_m})_{\cal C}.   
$$
Here ${\rm Cycle}_{m+1}$ means the sum of  cyclic shifts, 
$\overline \alpha \in {\Bbb H}^*$ is the shifted by one 
element $\alpha \in \overline {\rm H}^*(X)$, and $|\overline \alpha|$ is its degree.
The complex $({\cal C}_{{\Bbb  H}^*}, \delta)$ is nothing else but 
the {\it cyclic homology complex} for the algebra 
$\overline {\rm H}^*(X)$. 

\vskip 2mm
There is a subspace of {\it shuffle relations} in 
${\cal C}_{{\Bbb H}_*}$  generated by the 
elements
$$
\sum_{\sigma \in \Sigma_{p,q}}\pm (v_0 \otimes v_{\sigma(1)} 
\otimes \ldots \otimes v_{\sigma(p+q)})_{\cal C}, 
\qquad p, q \geq 1, 
$$ 
where the sum is over all $(p,q)$-shuffles, and 
the signs are given by the standard sign rule. 
Set
$$ 
{\cal C}{{\cal L}ie}^{\vee}_{{\Bbb H}^*}:= 
\frac{ {\cal C}_{{\Bbb H}^*}}{\mbox{Shuffle relations}};
\qquad 
{\cal C}{{\cal L}ie}_{{\Bbb H}_*}:= \mbox{the dual of 
${\cal C}{{\cal L}ie}^{\vee}_{{\Bbb H}^*}$}.
$$

  Let ${\cal L}ie_{{\Bbb H}_*}$ be the free graded Lie algebra 
 generated by ${\Bbb H}_*$. The space ${\cal C}{{\cal L}ie}_{{\Bbb H}_*}$ 
is the projection to the cyclic envelope of the space ${\cal L}ie_{{\Bbb H}_*} \otimes {\Bbb H}_*$.

\bt \la{MAATH} The map (\ref{1.17.08.10q}) induces a well defined linear map, 
called the {\rm Hodge correlator map}: 
\be \la{1.17.08.10as}
{\rm Cor}^*_{{\cal H}, a}: H^0_{\delta}\Bigl(
{\cal C}{\cal L}ie^{\vee}_{{\Bbb H}^*}\otimes {\cal H}^{\vee}\Bigr) \lra \C. 
\ee
-- It 
does not depend on the choices involved in the definition of the map 
(\ref{1.17.08.10q}). 
\et
We prove this in Section \ref{hf4.1ref}, 
as a consequence of an equivalent Theorem 
\ref{3/11/07/101q} proved there. 
\vskip 3mm
Now let $X$ be the set of complex 
points of an algebraic variety $X_{/\Q}$ over $\Q$. Then 
the rational de Rham cohomology $H^*_{\rm DR}(X_{/\Q}, \Q)$
 induce a $\Q$-rational structure on the 
space on the left in (\ref{1.17.08.10as}).  So we get a collection 
of numbers, the images  
of the elements of this $\Q$-vector space under the Hodge correlator map. 
We  call them  {\it the real periods of the homotopy type of $X$}.

\vskip 2mm
Dualising map (\ref{1.17.08.10as}) we get a homology class, 
{\it the Hodge correlator class}: 
\be \la{HCCL}
{\bf H}^*_{X, a}\in H_0^{\delta}\Bigl(
{\cal C}{\cal L}ie_{{\Bbb H}_*}\otimes {\cal H}\Bigr).
\ee

\paragraph{Describing special derivations.} 
Let 
$
\delta: {\Bbb H}_* \lra \Lambda^2{\Bbb H}_*
$ 
be the degree $1$ map  dualising the product on $\overline {\rm H}^*(X)$. 
See Section 3.2 for a discussion of signs. 
It gives rise to a differential $\delta$ 
of the free Lie algebra 
 ${\cal L}ie_{{\Bbb H}_*}$, providing it with a DG Lie algebra structure.

The image of the shifted by one fundamental class ${\rm H}_{2n}(X)[-1]$ 
under the map dualising the commutative product map 
${\rm H}^*(X) \otimes {\rm H}^*(X) \to {\rm H}^{2n}(X)$ 
provides a special element ${\rm S} \in {\cal L}ie_{{\Bbb H}_{*}}$, 
better understood as a map
$$
{\rm S}: {\rm H}_{2n}(X)[-1] \lra [{\Bbb H}_*(X), {\Bbb H}_*(X)] \subset {\cal L}ie_{{\Bbb H}_{*}}.
$$

\bd
A derivation of the Lie subalgebra ${\cal L}ie_{{\Bbb H}_{*}}$ 
is special if it kills the element ${\rm S}$. 
We denote by ${\rm Der}^S({\cal L}ie_{{\Bbb H}_*})$ the Lie algebra of all special derivations of 
the Lie algebra ${\cal L}ie_{{\Bbb H}_*}$. 
\ed

\vskip 2mm

The commutator with the differential $\delta$ on ${\cal L}ie_{{\Bbb H}_*}$  provides the Lie algebra 
${\rm Der}^S{\cal L}ie_{{\Bbb H}_*}$ with a DG Lie algebra structure.  
In Section 3.2 we show the following

\bp \la{1.17.08.10as1sd} (i) There is a natural isomorphism 
\be \la{1.17.08.10as1}
\theta: {\cal C}{\cal L}ie_{{\Bbb H}_*}\otimes {\cal H} \stackrel{\sim}{\lra}
{\rm Der}^S{\cal L}ie_{{\Bbb H}_*}. 
\ee
(ii) There is a Lie algebra structure on  
${\cal C}{\cal L}ie_{{\Bbb H}_*}\otimes {\cal H}$ which, 
together with the differential $\delta$, makes it into a DG Lie algebra. The 
isomorphism (\ref{1.17.08.10as1})
is an isomorphism of DG Lie algebras.
\ep
We prove it in Section 3.2. 
In particular, there is  an isomorphism of Lie algebras
\be \la{1.17.08.10as12}
\theta_0: H_0^\delta\Bigl({\cal C}{\cal L}ie_{{\Bbb H}_*}\otimes {\cal H}\Bigr) \stackrel{\sim}{\lra} 
H_0^\delta\Bigl({\rm Der}^S{\cal L}ie_{{\Bbb H}_*}\Bigr). 
\ee
It  is crucial for our definition of motivic correlators.

Here is a sketch of constructions leading to a proof of Proposition \ref{1.17.08.10as1sd}. 
First, given an element $F\otimes {\cal H} \in {\cal C}_{{\Bbb H}_*}\otimes {\cal H}$, 
there is a derivation $\theta_{F\otimes {\cal H}}$ of the tensor algebra 
${\rm T}_{{\Bbb H}_*}$ acting on the generators as follows:
$$
\theta_{F\otimes {\cal H}}: 
p \lms \sum_q \frac{\partial F}{\partial q}\langle q \cap {\cal H} \cap p \rangle. 
$$
Here $\frac{\partial F}{\partial q}$ is the non-commutative 
partial derivative, e.g. $$
\frac{\partial}{\partial q_1}(q_1 q_2 q_1 q_3)_{\cal C} = q_2 q_1 q_3 + q_3q_1 q_2 ,
$$ the sum is over a basis $\{q\}$ in ${\Bbb H}_*$, 
and $\langle q \cap {\cal H} \cap p \rangle$ 
is the iterated cap product with the fundamental class ${\cal H}$. 
Second, the condition  
$F \in {\cal C}{\cal L}ie_{{\Bbb H}_*}\otimes {\cal H}$ just means that 
$\theta_{F\otimes {\cal H}}$ preserves the Lie subalgebra ${\cal L}ie_{{\Bbb H}_*}$. 
Finally,  $\theta_{F\otimes {\cal H}}({\rm S}) =0$, i.e. it is a special derivation.

\begin{figure}[ht]
\centerline{\epsfbox{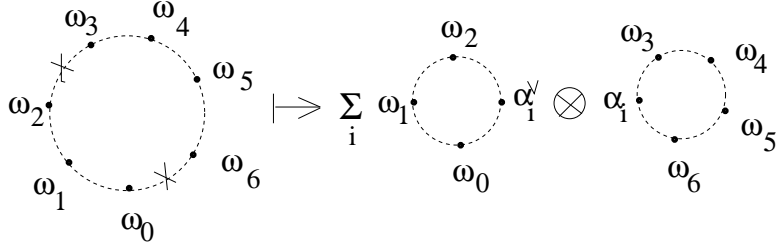}}
\caption{The Lie cobracket on ${\cal C}_{{\Bbb H}^*}\otimes {\cal H}^{\vee}$. 
The even factors ${\cal H}^{\vee}$ are not shown.}
\label{hf7}
\end{figure}

\vskip 2mm
There is a Lie cobracket on the space 
$
{\cal C}_{{\Bbb H}^*}\otimes {\cal H}^{\vee},  
$ 
making it into  a graded Lie coalgebra, and 
$({\cal C}_{{\Bbb H}^*}\otimes {\cal H}^{\vee}, \delta)$ into
 a DG Lie coalgebra. Namely, take a cyclic word, 
cut it into two pieces, insert the Casimir element 
$\sum_i\alpha_i^{\vee} \otimes \alpha_i$ 
and take the sum over all cuts with appropriate signs, as on Fig \ref{hf7}. 
Thus the dual 
${\cal C}_{{\Bbb H}_*}\otimes {\cal H}$ is a DG Lie algebra. 
Furthermore, 
${\cal C}{{\cal L}ie}^{\vee}_{{\Bbb H}^*} \otimes {\cal H}^{\vee}$ 
inherits a DG Lie coalgebra structure. 
So  
${\cal C}{{\cal L}ie}_{{\Bbb H}_*}\otimes {\cal H}$ is a DG Lie algebra.
Therefore   
$H^0_{\delta}\Bigl({\cal C}{{\cal L}ie}^{\vee}_{{\Bbb H}^*}\otimes {\cal H}^{\vee}\Bigr)$  
is a Lie coalgebra.

\paragraph{Functoriality of the Hodge correlator classes.} Given a map $f:X \to Y$ of Kahler manifolds, 
the map $f_*: H_*(X) \to H_*(Y)$ induces a Lie algebra map 
$$
f_*: {\cal C}{\cal L}ie_{{\Bbb H}_*(X)} \to 
{\cal C}{\cal L}ie_{{\Bbb H}_*(Y)}. 
$$
It commutes with the differentials (Lemma \ref{6.13.08.1}). Here is our key result:

\bt \la{PPPq} 
For any map $f:X \to Y$ of Kahler manifolds 
the derivation ${\bf H}^*_{Y, f(a)}$ preserves the 
Lie subalgebra $f_*({\cal L}ie_{{\Bbb H}_*(X)})$, and its restriction to the latter 
is  homotopic to  $f_*({\bf H}^*_{X, a})$. 
\et

Theorem 
\ref{MAATH} can be viewed as a special case of 
Theorem \ref{PPPq} when $f:X \to X$ is the identity map, but we use
two different splittings/Green currents to 
define the Hodge correlators on $X$. We prove Theorem \ref{PPPq} in Section \ref{hf4.2ref}. 
We prove a more precise Theorem \ref{7.12.08.1}, calculating the coboundary explicitly.\footnote{There 
is a different, less computational proof, which we discuss elsewhere. 
We also hope to give a different proof using the Feynman integral from Section \ref{hf1.6}  and 
the BV formalism elsewhere.}

\vskip 2mm

{\bf Remark}. Let $({\cal L}_\bullet, \delta)$ be a DG Lie algebra with a differential 
$\delta$. Denote by ${\rm Der}{\cal L}_\bullet$ the Lie algebra of derivations 
of the Lie algebra ${\cal L}_\bullet$. The commutator with the differential $\delta$ is a 
differential on the Lie algebra ${\rm Der}{\cal L}_\bullet$, making it into a DG Lie algebra.
 Then one has 
the following elementary Lemma.
\bl 
The Lie algebra $H_0^{\delta}\Bigl({\rm Der}{\cal L}_\bullet\Bigr)$ consists of 
the derivations of the Lie algebra ${\cal L}_\bullet$ 
preserving the grading, commuting  with 
the differential $\delta$, and defined up to a homotopy. 
\el
 The degree zero cycles 
are the derivations preserving the grading and commuting 
with the differential. 
The degree zero boundaries are the differentials homotopic to zero. 

Applying the map (\ref{1.17.08.10as12}) to the Hodge correlator class (\ref{HCCL}) 
class we arrive at  an element
\be \la{6.7.08.1}
{\bf H}^*_{X, a}\in H_0^{\delta}\Bigl({\rm Der}^S{\cal L}ie_{{\Bbb H}_*}\Bigr).
\ee

In the next two subsections we explain how this element 
defines a functorial real mixed Hodge 
structure on the rational homotopy type of compact Kahler manifolds.

\subsection{Rational homotopy type DG Lie algebras} 
Let $M$ be a smooth simply connected 
manifold, $a\in M$. According to Quillen \cite{Q} and Sullivan   
\cite{S}, 
the rational homotopy type of $M$ is described by a DG Lie algebra 
over $\Q$ concentrated in degrees 
$0, -1, -2, ...$. It is well defined up to a quasiisomorphism. 
There is a similar DG Lie algebra 
${\rm L}(M, a)$ 
over $\Q$ in the case when $M$ is not necessarily simply 
connected. Its zero cohomology $H^0{\rm L}(M, a)$  is  a Lie algebra, 
isomorphic to  the pronilpotent completion 
of $\pi_1(M, a)$.

\vskip 3mm
The DG Lie algebra ${\rm L}(M, a)\otimes \C$ appears as follows. 
Let ${\rm DGCom}$ (respectively ${\rm DGCoLie}$) be the category of 
differential graded commutative algebras (respectively Lie coalgebras) 
over $\Q$. 
Then there is a pair of adjoint functors 
(the Bar construction followed to the projection to the 
indecomposables, and the Chevalley complex)
$$
{\cal B}: {\rm DGCom} \lra {\rm DGCoLie}, \qquad  
{\cal S}: {\rm DGCoLie} \lra {\rm DGCom}, 
$$
Take the de Rham DGCom ${\cal A}^*(M)$ over $\C$. Given a point $a \in M$, 
the evaluation map  ${\rm ev}_a: 
{\cal A}^*(M)\to \C$ provides  a DGCom denoted 
${\cal A}^*(M, a)$. Applying to it the 
functor ${\cal B}\otimes \C$ we get a DGCoLie quasiisomorphic to the 
Lie coalgebra dual to  
${\rm L}(M, a)\otimes \C$. 

Suppose now that $M$ is a compact Kahler manifold. 
Then, thanks to the formality theorem \cite{DGMS}, 
the DG commutative algebra ${\cal A}^*(M, a)$ is quasiisomorphic to 
its cohomology, the reduced cohomology  $\widetilde H^*(M, \C)$ 
of $X$. So applying the 
functor ${\cal B}$ to the rational reduced cohomology $\widetilde H^*(M, \Q)$ 
we get the minimal model of the rational homotopy type DG Lie algebra. 
We reserve from now on the notation  ${\rm L}(M, a)$ to this DG Lie algebra.

\vskip 3mm
Let $X$ be a smooth irreducible projective 
variety over a field $k$. Let us choose an 
embedding $\sigma: k \subset \C$. Then the
 rational homotopy type DG Lie algebra 
${\rm L}(X(\C), a)$ is supposed to be 
an object of motivic origin: 
One should exist a DG Lie algebra ${\rm L}_{\cal M}(X, a)$ 
in the yet hypothetical abelian tensor 
category of mixed motives over $k$, 
whose Betti realization is the DG Lie algebra 
${\rm L}(X(\C), a)$. The Hodge realization ${\rm L}_{\rm Hod}(X, a)$ is known thanks to 
Morgan \cite{M}.

\vskip 2mm
Let ${\cal L}ie_{\widetilde {\Bbb H}_{*}}$ be 
the free graded Lie algebra 
 generated by 
$$
\widetilde {\Bbb H}_{*}:= {\rm H}_{>0}(X)[-1], 
$$ 
sitting in degrees $[-(2n-1), 0]$. The Lie algebra ${\cal L}ie_{\widetilde {\Bbb H}_{*}}$
has DG Lie algebra with the differential dualising the product on ${\rm H}^{>0}(X)$. 
The formality theorem \cite{DGMS} implies:

\begin{theorem} \label{5.12.07.1}
The associate graded 
for the weight filtration 
${\rm gr}^W{\rm L}(X, a)$ 
of the rational homotopy type  DG Lie algebra  is 
quasiisomorphic to the DG Lie algebra ${\cal L}ie_{\widetilde {\Bbb H}_*}$.
\end{theorem}

The zero cohomology $H^0{\rm L}_{\cal M}(X, a)$ 
is supposed to be a Lie algebra 
isomorphic to the (conjectural) unipotent motivic fundamental group of $X$. 
The latter was defined unconditionally for 
unirational $X$'s over number fields \cite{DG}.

\subsection{Hodge correlators and a real MHS on the rational homotopy type}

\paragraph{The real Hodge Galois group.} The category of real mixed Hodge structures is an abelian  
tensor category \cite{D}. 
The associate graded for the weight filtration provides a fiber functor 
$H \lms {\rm gr}^WH$
 on this category. Therefore it is canonically 
equivalent to the category of representations 
of a pro-algebraic group over $\R$, called the Hodge Galois group $G_{\rm Hod}$. 
It  
is isomorphic to 
a semidirect product of $\C^*_{\C/\R}$ and a prounipotent algebraic group 
$U_{\rm Hod}$: 
$$
0 \lra U_{\rm Hod} \lra G_{\rm Hod} \lra \C^*_{\C/\R} \lra 0, \qquad 
\C^*_{\C/\R} \hra G_{\rm Hod}.  
$$
The action of the group $\C^*_{\C/\R}$ provides 
the Lie algebra of the unipotent group 
${\rm U}_{\rm Hod}$ with a structure of a Lie algebra in the category of 
direct sums of pure 
real Hodge structures. We denote this Lie algebra by ${\rm L}_{\rm Hod}$. 
So $G_{\rm Hod}$-modules are nothing else but 
${\rm L}_{\rm Hod}$-modules in the category of real Hodge structures. 
 
According to Deligne \cite{D2}, 
the Lie algebra ${\rm L}_{\rm Hod}$  
is generated by certain generators 
$g_{p,q}$, $p,q \geq 1$, satisfying the only relation $\overline g_{p,q} = -g_{q, p}$. 
In Section 4 of \cite{G1} we defined a 
different set of generators, denoted below 
by $G^*_{p,q}$ and called the canonical generators --  
they behave nicer in families.\footnote{There are, in fact, two natural normalizations of the generators, denoted by $G^*_{p,q}$ and $G_{p,q}$ -- 
they differ by certain binomial coefficients depending 
on $p,q$. Both play a role in our story. 
We use $\ast$ in the objects related to the $G_{p,q}$-generators, 
e.g. ${\rm Cor}_{\cal H}^*$. 
 The Hodge correlators 
related to the generators $G^*_{p,q}$ are defined by using the operator 
$\omega$, see Section 2.1, instead of the operator $\xi$ in (\ref{xixi}).} 
In the mixed Tate case (i.e. $p=q$) they essentially 
coincide with the ones defined 
in \cite{L}. 
So  a real mixed Hodge structure is nothing else but the following data: 

\begin{itemize} \la{!!}
\item
a real Hodge structure $V$, i.e. a real vector space $V$, 
whose complexification $V_\C$ is bigraded, 

\item a 
collection of imaginary operators  $\{G^*_{p,q}\}$ on $V_\C$, $p,q >0$, of 
the bidegree by $(-p, -q)$, such that $\overline G^*_{p,q} =  - G^*_{q, p}$. 
\end{itemize}
Indeed, a representation  of $\C^*_{\C/\R}$ in $V$ 
is the same thing as a real Hodge structure in $V$, and the operators  
$G^*_{p,q}$ are the images of the canonical generators 
of   ${\rm L}_{\rm Hod}$. We 
encode a collection of the operators $\{G^*_{p,q}\}$ by a single operator
$$
G^* = \sum_{p,q>0}G^*_{p,q}.
$$

\paragraph{The Hodge Galois group action on the rational homotopy type DG Lie algebra.}
\begin{theoremconstruction} \la{4.25.08.22}
There is a canonical map 
\be \la{4.25.08.21qq}
{\rm L}_{\rm Hod} \lra H_0^{\delta} 
\Bigl({\rm Der}^S{\cal L}ie_{{\Bbb H}_*}\Bigr), \qquad G^* \lms {\Bbb H}^*_{X,a}, 
\ee 
The action of the pure Hodge group ${\Bbb C}_{\C/\R}^*$ 
is provided by  the real  Hodge structure  
on $H^{*}(X)$.
 \end{theoremconstruction}
This 
leads to a real MHS on $H^{\delta}_*{\rm L}(X, a)$ as follows. 

\bd \la{spes}
A derivation of the graded Lie algebra ${{\cal L}ie}_{\widetilde {\Bbb H}_*}$ 
is {\it special} if it kills both the 
generator $H_{2n}(X)[-1]$ and the element ${\rm S} = \delta (H_{2n}(X)[-1])$. 
\ed
Denote by ${\rm Der}^S{{\cal L}ie}_{\widetilde {\Bbb H}_*}$ the Lie algebra of all 
special derivations. The inclusion 
${\Bbb H}_*\hra \widetilde {\Bbb H}_*$ provides  an isomorphism of graded Lie algebras
\be \la{mapqc}
{\rm Der}^S{{\cal L}ie}_{{\Bbb H}_*}\stackrel{\sim}{\lra} {\rm Der}^S{{\cal L}ie}_{\widetilde {\Bbb H}_*}.
\ee
One easily checks that the differentials preserve the subalgebras of special derivations. 
So the map 
(\ref{mapqc}) is an isomorphism of DG Lie algebras.
Therefore it induces an isomorphism
$$
H_0^\delta\Bigl({\rm Der}^S{\cal L}ie_{{\Bbb H}_*}\Bigr) =  
H_0^\delta\Bigl({\rm Der}^S{{\cal L}ie}_{\widetilde {\Bbb H}_*}\Bigr). 
$$
The element ${\Bbb H}^*_{X,a}$ lives in the Lie algebra on the left. 
We transform it to an element $\widetilde {\Bbb H}^*_{X,a}$ of the Lie algebra on the right. 
So we get 
a map 
\be \la{4.221qqw}
{\rm L}_{\rm Hod} \lra H_0^{\delta} 
\Bigl({\rm Der}^S{\cal L}ie_{\widetilde {\Bbb H}_*}\Bigr), \qquad G^* \lms \widetilde {\Bbb H}^*_{X,a}, 
\ee 
The Lie algebra  $H_0^{\delta} 
\Bigl({\rm Der}^S{\cal L}ie_{\widetilde {\Bbb H}_*}\Bigr)$ acts 
by derivations of ${\cal L}ie_{\widetilde {\Bbb H}_*}$, defined up to a homotopy. 
Finally,  ${\cal L}ie_{\widetilde {\Bbb H}_*}$ is 
isomorphic to 
 ${\rm gr}^W{\rm L}(X, a)$ by Theorem \ref{5.12.07.1}. 

Summarising, we defined a real MHS on ${\rm L}(X, a)$. 
When the point $a$ varies, we get a variation of real MHS.  
The proof of this is completely similar to the one given in \cite{G1} 
in the case of curves, and thus is omitted. 
 
\vskip 3mm
In particular, we get a real MHS on the pronilpotent completion 
$\pi_1^{\rm nil}(X,a)$. 
\bt \la{FMHSR} Let $X$ be a complex regular projective  variety. Then 
the  real MHS on $\pi_1^{\rm nil}(X,a)$ provided by the 
Hodge correlators coincides with the classical one. 
\et
{\bf Proof}. We proved in \cite{G1} that this is true for any smooth curve $X$. 
The Hodge correlator real MHS on $\pi_1^{\rm nil}(X,a)$ is 
functorial by Theorem \ref{PPPq}. 
By the Lefschetz theorem there exists a plane section $Y$ of $X$ 
such that $\pi_1(Y)$ surjects onto $\pi_1(X)$. 
This implies Theorem \ref{FMHSR}.  

\vskip 3mm
The real MHS on $\pi_1^{\rm nil}(X,a)$ can be  described 
via Chen's iterated integrals. Theorem \ref{FMHSR} 
tells that the Hodge correlators provide the same 
real periods as iterated integrals. It would be interesting to 
relate 
the Hodge correlators to the iterated integrals directly. 
I know how to do this only in the simplest cases.

\subsection{Motivic correlators for a regular projective variety}
\paragraph{Motivic Galois groups.} One expects that there exists an abelian tensor category 
${\cal M}{\cal M}_k$ of mixed motives over a field $k$. 
Each object $M$ of this category should 
carry a canonical weight filtration $W_{\bullet}$. The functor 
$$
\omega: {\cal M}{\cal M}_k \stackrel{?}{\lra} {\cal P}{\cal M}_k, \qquad
M \lms {\rm gr}^{W}M
$$
should be a fiber functor with values in the semisimple tensor category 
${\cal P}{\cal M}_k$ of pure motives. The algebraic group scheme of all 
automorphisms of this fiber functor is called the motivic Galois group 
$G_{\rm Mot}$. Thanks to the Tannakian formalism, the functor $\omega$ 
then provides a canonical equivalence between 
the category of all mixed motives and the category 
of representations of $G_{\rm Mot}$ in the category of pure motives. 
So the  motivic Galois group 
$G_{\rm Mot}$ contains all the information telling us how the 
mixed motives are obtained from the pure ones. 

The motivic Galois group $G_{\rm Mot}$ is a prounipotent algebraic group 
in the tensor category of pure motives. 
We denote by ${\rm L}_{\rm Mot}$ its Lie algebra, and by ${\cal L}_{\rm Mot}$ 
the dual Lie coalgebra -- see, say,  \cite{G4} for more details.

\paragraph{Motivic correlators.} The motivic Galois group should act  
by symmetries of the rational motivic homotopy type of  a variety
$X$. The latter, however, is well defined only up to 
a quasiisomorphism. Here is a precise statement:

{\it One should exist a 
homomorphism  of Lie algebras in the category of pure motives} 
\be \la{4.25.08.21g}
{\rm L}_{\rm Mot} \stackrel{?}{\lra} H_0^{\delta} 
\Bigl({\rm Der}({\rm gr}^W{\rm L}_{\cal M}(X,a))\Bigr). 
\ee

The $\R$-Hodge realization of the 
motivic Lie algebra ${\rm L}_{\rm Mot}$ 
is the Hodge Lie algebra ${\rm L}_{\rm Hod}$. So the 
 $\R$-Hodge realization of the map (\ref{4.25.08.21g}) should give the map (\ref{4.25.08.21qq}).

\vskip 2mm
Denote by ${\cal L}ie_{{\rm Mot}({\Bbb H}_*)}$ the free graded Lie algebra 
generated by the reduced 
homology motive of $X$ shifted to the right by $1$. It is a DG Lie algebra. 
Similarly, there is a slightly bigger DG Lie algebra 
${\cal L}ie_{ {\rm Mot}(\widetilde {\Bbb H}_*)}$. Theorem \ref{5.12.07.1} 
tells that one should exist a canonical isomorphism
\be \la{CIS}
{\cal L}ie_{{\rm Mot}(\widetilde {\Bbb H}_*)}
 \stackrel{?}{=}{\rm gr}^W{\rm L}_{\cal M}(X,a).
\ee 
Just like in Definition \ref{spes}, one defines  special derivations of the Lie algebra 
on the left, and hence on the right of (\ref{CIS}).  
The motivic Galois Lie algebra should act by special derivations 
of the  Lie algebra 
${\rm gr}^W{\rm L}^{\rm Mot}(X,a)$. 
There is a canonical isomorphism
\be \la{1.17.08.10as12s}
{\rm Der}^S{\cal L}ie_{{\rm Mot}({\Bbb H}_*)} \stackrel{\sim}{\lra} 
{\rm Der}^S{\cal L}ie_{{\rm Mot}(\widetilde {\Bbb H}_*)}. 
\ee
Therefore  the map 
(\ref{4.25.08.21g}) plus isomorphisms (\ref{1.17.08.10as12s}) and  
(\ref{CIS}) provide a canonical map of Lie algebras
\be \la{4.25.08.21q}
{\rm L}_{\rm Mot} \stackrel{?}{\lra} H_0^{\delta} 
\Bigl({\rm Der}^S{\cal L}ie_{{\rm Mot}({\Bbb H}_*)}\Bigr). 
\ee
Dualizing the map (\ref{4.25.08.21q}) and using the motivic version of 
isomorphism (\ref{1.17.08.10as12}) we get a morphism of Lie coalgebras, 
called the 
{\it motivic correlator map} for  $X$:
\be \la{4.25.08.21rty}
{\rm Cor}_{\rm Mot}: 
H^0_{\delta}\Bigl({\cal C}{\cal L}ie^{\vee}_{{\rm Mot}({\Bbb H}_*)}\otimes {\cal H}^{\vee} \Bigr)
\stackrel{?}{\lra}
 {\cal L}_{\rm Mot}.
\ee
Its source is easily 
defined in terms of the cohomology motive of $X$.  The target is remarkable: 
it carries all the information 
how the mixed motives are glued from the pure ones. The cobracket 
of an element in ${\cal L}_{\rm Mot}$  contains all the information about 
the complexity of the element. Therefore, since 
(\ref{4.25.08.21rty}) is a Lie coalgebra map,  and 
the cobracket in the left Lie coalgebra  
is very transparent, see Fig \ref{hf7}, we get a full control 
of the arithmetic complexity of motivic correlators.

\paragraph{Relating motivic and Hodge correlators.} 
The hypothetical motivic Lie coalgebra ${\cal L}_{\rm Mot}$ 
is an inductive limit of objects in the  
semisimple category ${\cal P}{\cal M}_k$ 
of pure motives over $k$. 
One can decompose it into the isotipical components:
$$
{\cal L}_{\rm Mot} \stackrel{?}{=} \bigoplus_{M \in {\rm Iso}({\cal P}{\cal M}_k)}
{\cal L}_{M}\otimes M_{}^{\vee}, \qquad 
{\cal L}_{M} = {\rm Hom}_{{\cal P}{\cal M}_k}(M^{\vee}, {\cal L}_{\rm Mot}).
$$
Here $M$ runs through the isomorphism classes of simple objects in 
the category ${\cal P}{\cal M}_k$, $M^{\vee}$ is the dual motive, 
and ${\cal L}_{M}$ 
is a $\Q$-vector space. 

One can do this in the Hodge realization. 
Namely, denote by ${\cal L}_{\rm Hod}$ the $\R$-Hodge Lie coalgebra. 
Then there is a decomposition into the isotipical components:
$$
{\cal L}_{\rm Hod} = \bigoplus_{H \in {\rm Iso}({\cal P}{\rm Hod}_\R)}
{\cal L}_{H}\otimes H_{}^{\vee}.
$$

\bp \la{cpm} 
There is a canonical linear map, called {\it the real period map}: 
$$
{\rm P}_\R: {\cal L}_{\rm Hod} \lra i\R.
$$
\ep

The Hodge and motivic correlator are compatible as follows. 
\bp \la{10.21.08.1}
Let $X$ be a regular projective variety over $\R$. Then 
the composition 
$$
H^0_{\delta}\Bigl({\cal C}{\cal L}ie^{\vee}_{{\rm Hod}({\Bbb H}_*)}\otimes 
{\cal H}^{\vee} \Bigr)~
\stackrel{{\rm Cor}_{\rm Hod}}{\lra} ~{\cal L}_{\rm Hod} ~
\stackrel{{\rm P}_\R}{\lra} ~i\R.
$$
coincides with  the Hodge correlator map 
${\rm Cor}_{\cal H}$. 
\ep
Therefore, given a regular projective variety $X$ over $\Q$, and 
assuming the existence the motivic correlator map  ${\rm Cor}_{\rm Mot}$, 
the Hodge correlator map coincides with the composition
$$
H^0_{\delta}\Bigl({\cal C}{\cal L}ie^{\vee}_{{\rm Mot}({\Bbb H}_*)}\otimes 
{\cal H}^{\vee} \Bigr)~
\stackrel{{\rm Cor}_{\rm Mot}}{\lra}~ {\cal L}_{\rm Mot} ~
\stackrel{r_{\rm Hod}}{\lra} ~
{\cal L}_{\rm Hod} ~\stackrel{{\rm P}_\R}{\lra}~ i\R.
$$
Here $r_{\rm Hod}$ is the $\R$-Hodge realization functor. 

\vskip 3mm
\subsection{A Feynman integral for Hodge correlators} \la{hf1.6} 
Let ${\cal L}_1$ and ${\cal L}_2$ be two DG Lie algebras. 
Let us forget for a moment about their differentials, treating them as graded Lie superalgebras. 
Then there is an affine ${\Bbb G}_m$-superscheme\footnote{That is, an affine
 superscheme equipped with an action of the algebraic group ${\Bbb G}_m$} of graded 
Lie superalgebra maps ${{\cal H}om}_{\rm Lie}({\cal L}_1, {\cal L}_2)$. 
Its points with values in a graded supercommutative algebra $A$ are defined by 
$$
{{\cal H}om}_{\rm Lie}({\cal L}_1, {\cal L}_2)(A):= {\rm Hom}_{\rm Lie}({\cal L}_1, {\cal L}_2\otimes A).
$$
Denote by ${\cal L}ie_{V}$ the free graded Lie superalgebra 
generated by a finite dimensional graded vector  space $V$. 
Let ${\cal G}$ be a graded Lie superalgebra. 
Then there is an isomorphism of affine ${\Bbb G}_m$-superschemes
\be \la{4.20.08.6}
{{\cal H}om}_{\rm Lie}({\cal L}ie_{V^{*}}, {\cal G}) = V \otimes {\cal G}.
\ee
Suppose now that ${\cal L}ie_{V^*}$ has a differential 
$\Delta$ providing it with a  DG Lie algebra structure. 
Since $\Delta$ is an odd (degree $1$) infinitesimal automorphism 
of the Lie superalgebra ${\cal L}ie_{V^*}$, 
 it is transformed by isomorphism (\ref{4.20.08.6}) 
to an infinitesimal 
automorphism of the affine  ${\Bbb G}_m$-superscheme $V \otimes {\cal G}$, that is 
to a homological vector field $v_{{\Delta}}$ on 
$V \otimes {\cal G}$, making it into an affine DG scheme -- 
see Section \ref{hf5.2ref} for background. 

Furthermore, let ${\cal D}$ be a derivation 
of the graded Lie algebra ${\cal L}ie_{V}$ commuting with  
$\Delta$. Then isomorphism (\ref{4.20.08.6}) transforms it  
to a homological vector field $v_{{\cal D}}$ commuting with 
$v_{{\Delta}}$. 
 
\vskip 3mm
{\bf Example.} 
Let $({\cal L}, d)$ be a DG Lie algebra. Then ${\cal L}[1]$ is a DG space 
with the homological vector field given by the Chern-Simons vector field 
$Q_{\rm CS}$ 
(see Section \ref{hf6.1ref}):
$$
\stackrel{{\cdot}}{\alpha} = Q_{\rm CS}(\alpha) = 
d \alpha + \frac{1}{2}[\alpha, \alpha]. 
$$
In particular,  given a DG commutative algebra $A$  and a Lie algebra ${\cal G}$, 
there is  a DG Lie algebra 
$A \otimes {\cal G}$ with the Lie bracket 
$[a_1 \otimes g_1, a_2 \otimes g_2]:= a_1a_2 \otimes [g_1, g_2]$, and hence 
a DG scheme  $A[1] \otimes {\cal G}$.

\vskip 3mm
The functor ${\cal B}: {\rm DGCom} \to {\rm DG Colie}$ takes an 
$A \in {\rm DGCom}$ to 
 the cofree Lie coalgebra ${\rm CoLie}(A[1])$ 
cogenerated by the graded space $A[1]$, with a certain differential $\Delta$.  

The functor ${\cal B}$  is characterised by the following universality property: 
\begin{itemize}
\item One has a functorial isomorphism of the affine DG schemes:
\be \la{DGISO}
{{\cal H}om}_{\rm CoLie}({\cal G}^{\ast}, {\cal B}(A)) = A[1] \otimes {\cal G}.
\ee
\end{itemize}
In other words, the functorial isomorphism (\ref{DGISO})  
transforms  
the differential $\Delta$ on ${\cal B}(A)$ 
to the 
Chern-Simons homological vector field $Q_{\rm CS}$ on $A[1] \otimes {\cal G}$.  
\vskip 3mm

Applying (\ref{DGISO}) to the algebra 
$A = \overline {\rm H}^*(X)$, we get  
an isomorphism of DG schemes
\be \la{4.20.08.5}
{{\cal H}om}_{\rm Lie}({\cal L}ie_{{\Bbb H}_*}, {\cal G}) = 
{\Bbb H}^*(X) \otimes {\cal G} =: {\cal G}_{{H}}[1].
\ee 
The differential $\delta$ on the DG Lie algebra ${\cal L}ie_{{\Bbb H}_*}$ 
is transformed to the Chern-Simons homological 
vector field $Q_{CS}$ on ${\cal G}_{{H}}[1]$. 

\vskip 3mm
Let ${\cal G}$ be a Lie 
algebra with a  non-degenerate 
scalar product $Q(\ast, \ast)$. 
For example one can take the Lie algebra ${\rm Mat}_N$ of 
$N\times N$ matrices 
with $Q(A,B):= {\rm Tr}(AB)$. 

The Hodge correlator ${\bf H}^*_{X,a}$, see (\ref{6.7.08.1}),  is a derivation of 
${\cal L}ie_{{\Bbb H}_*}$ commuting with $\delta$ and defined up to a homotopy. 
So it is transformed by isomorphism (\ref{4.20.08.5}) 
to a homological vector field ${\rm Q}_{\rm Hod}$ on ${\cal G}_{{H}}[1]$, 
commuting with the Chern-Simons vector field ${\rm Q}_{\rm CS}$ and 
well defined up to commutators 
$[{\rm Q}_{\rm CS}, R]$. 
Our next goal 
is to show 
that we get a Hamiltonian vector field, whose Hamiltonian 
is given by a Feynman integral introduced below.

\vskip 3mm

Given  ${\cal G}$-valued differential forms $\varphi_i$ on $X$, 
we define a differential  form
$$
\langle \varphi_1, \varphi_2,\varphi_3\rangle:= Q(
\varphi_1, [\varphi_2,\varphi_3]\rangle.
$$ 

Choose a splitting of the bigraded de Rham complex  ``harmonic'' forms ${\cal H}ar_X$ 
and its orthogonal complement.  
Let ${\rm Ker}~d^\C\subset {\cal A}^{*}(X)$ be the kernel of $d^\C$. 
Let us define, using the splitting,  a function on the functional space 
$$
{\rm Ker}~d^\C\otimes {\cal G}[1].   
$$ 
Let $\psi \in {\rm Ker}~d^\C\otimes {\cal G}[1]$. 
Write ${\rm Ker}~d^\C
 = {\cal H}ar_X \oplus {\rm Im}~d^\C$, and 
$$
\psi = \psi_0 + \alpha,  
\quad \psi_0 \in {\rm Im}~d^\C\otimes {\cal G}[1] , \quad \alpha \in 
{\cal H}ar_X\otimes {\cal G}[1].
$$
Choose $\varphi$ such that $d^\C\varphi = \psi_0$. Then there is an  action: 
\be \la{ALP}
S(\psi) = \int_X\frac{1}{2}( \varphi, d\psi ) + 
\frac{1}{6}\langle \varphi, \psi,
 \psi \rangle, \qquad d^\C\varphi = \psi_0.
\ee
Observe that $S(\psi)$ is independent of 
the choice of $\varphi$. Indeed, 
 $\int_X \langle d^\C \eta, d \psi\rangle  = 
\int_X \langle d^\C \eta, \psi, \psi\rangle =0$ since $d^\C\psi=0$. 
Notice that $d^\C$ (Lagrangian in (\ref{ALP})) is very close to 
the Chern-Simons Lagrangian: 
$$
d^\C\Bigl(\frac{1}{2}(\varphi, d\psi ) + 
\frac{1}{6}\langle \varphi, \psi,
 \psi \rangle\Bigr) = \frac{1}{2}\langle \psi_0, d\psi \rangle + 
\frac{1}{6}\langle \psi_0, \psi,
 \psi \rangle
$$
We  would like to ``integrate'' 
$e^{ iS(\psi)}$ along the fibers of the natural projection 
provided by the splitting of the de Rham complex: 
$$
{\rm Ker}~d^\C\otimes {\cal G}[1] \lra H^*(X)\otimes {\cal G}[1],
 \qquad \psi = \psi_0 + \alpha \lms \alpha.  
$$
However the quadratic part 
$(\varphi, d\psi)$ of the Lagrangian 
is highly degenerate. Precisely, by  the classical 
Hodge theory 
there exists a subspace $F \subset {\cal A}^{*,*}_X$ 
such that there is an isomorphism of bicomplexes
\be \la{3.3.08.10}
{\cal A}^{*,*}_X = {\cal H}ar_X ~~\bigoplus ~~
\begin{array}{ccc}
d^\C F&\stackrel{d}{\lra}
 &d^\C d F \\
d^\C \uparrow &&d^\C \uparrow \\
F &\stackrel{d}{\lra} & d F 
\end{array}  
\ee
Here all arrows in the square are isomorphisms. 
Then 
$$
{\rm Ker}~d^\C = {\cal H}ar_X \oplus d^\C F \oplus d^\C d F, 
$$
and kernel of the quadratic form 
$(\varphi, d\psi)$ is the subspace 
${\cal H}ar_X \oplus 
d^\C d F$. The form is non-degenerate 
on $d^\C F$.  
We integrate over  
$\psi \in d^\C F$, i.e. over a ``a quarter $F$ of the de Rham complex'',
getting a function, the Hodge correlator, which  
depends on the decomposition (\ref{3.3.08.10}):
\be \la{4.18.08.1}
\widetilde {\rm Cor}_{\cal H}^{*}(\alpha):= 
\int e^{iS(\alpha + d^\C \varphi)}{\cal D}\varphi, \qquad \varphi\in F.
\ee
Elaborating 
$$
S(\alpha + d^\C \varphi):= \int_X 
\frac{1}{2}(\varphi, dd^\C \varphi) 
+ \frac{1}{6}
\langle \varphi, d^\C  \varphi, d^\C  \varphi\rangle 
+ \frac{1}{3}\langle \alpha, \varphi, d^\C \varphi \rangle 
+ \frac{1}{6}\langle \alpha, \alpha, \varphi \rangle,  
$$
we use the perturbative series expansion 
on the tree level to make sense of the integral, treating 
$\alpha$ as a small parameter -- see Section \ref{hf5.3sec}. 
It involves only plane trivalent trees decorated by 
the harmonic form $\alpha$, i.e. it amounts to minimizing 
the action:
$$
\widetilde {\rm Cor}_{\cal H}^{*}(\alpha)= {\rm Min}_{\varphi \in F}
S(\alpha + d^\C \varphi).
$$

The space  
${\cal G}_{{\rm H}}[1]$ has an even  symplectic structure 
$$
\omega(\overline h_1 \otimes g_1, \overline h_2 \otimes g_2):=  (-1)^{|h_1|}Q(g_1, g_2) \int_X h_1 \wedge h_2.
$$
 So function (\ref{4.18.08.1}) gives rise to a 
Hamiltonian vector field  on  
${\cal G}_{{\rm H}}[1]$:
\be \la{18.4.08.8}
\widetilde {\rm Q}^{*}_{\rm Hod}
 \in S^{>0}({{\cal G}_{{\rm H}}[1]}^{\vee})\otimes
 {\cal G}_{{\rm H}}[1].
\ee
Here on the right stands 
the reduced Chevalley cochain complex 
of the graded Lie algebra ${\cal G}_{{\rm H}}$ with the adjoint 
coefficients. It has a differential $\delta$ \, 
given by the commutator with $Q_{\rm CS}$.

\begin{theorem} \la{18.4.08.10}
i) The element (\ref{18.4.08.8}) provides a 
well defined, functorial  cohomology class
\be \la{CLass}
{\rm Q}^*_{\rm Hod} \in H_0^{\delta}\Bigl({\rm Der}{\cal G}_{{\rm H}}\Bigr). 
\ee

ii) The class  (\ref{CLass}) is the image of the class  ${\bf H}^*_{X,a}$ 
under the isomorphism (\ref{4.20.08.5}). 
\et

The part i) means that the cohomology class (\ref{18.4.08.8})
 does not depend 
on the choice of the subspaces $F$ and ${\cal H}ar_X$, as well as on 
the choice of a Green current involved in 
its perturbative series expansion 
definition. The part i) is deduced  from ii) and  Theorem \ref{MAATH}. 
The part ii) is obtained by comparing the  
correlators from Sections 2.3/\ref{hf6.a} with the ones from Section \ref{hf5.4ref}.

\vskip 3mm
We view the class $Q^*_{\rm Hod}$ as a  vector field on the DG space ${\cal G}_H[1]$
 -- we call it the {\it Hodge vector field}. 
The part i) of Theorem \ref{18.4.08.10} just means that 
it commutes with the Chern-Simons field ${\rm Q}_{CS}$, 
and is well defined up to commutators with ${\rm Q}_{CS}$. Let $\varepsilon$ be an odd element 
of degree $1$, so $\varepsilon^2=0$. Then the part i) of Theorem \ref{18.4.08.10} 
is equivalent to the

\bt \la{DEFO}
The vector field ${\rm Q}_{CS} + \varepsilon Q^*_{\rm Hod}$ is a homological vector field on the 
 DG space ${\cal G}_H[1]$, well defined up to commutators $[{\rm Q}_{CS}, \ast]$.

So  a real mixed Hodge structure on the rational homotopy type of $X$ is the same thing as 
a deformation of the DG space ${\cal G}_H[1]$ over the odd line ${\rm Spec}\R[\varepsilon]$,
 functorial in ${\cal G}$, considered up to an isomorphism. 
\et

\paragraph{The Hodge Galois group and complex local systems}

An $\R$-mixed Hodge structure on the pronilpotent fundamental group 
$\pi_1^{\rm nil}(X,x)$ tell us  that there is an 
extension of Lie algebras, where ${\cal G}_{\rm Hod}$ is the Lie algebra of the 
Hodge Galois group $G_{\rm Hod}$: 
\be \la{6.6.08.1}
0 \lra \pi_1^{\rm nil}(X,x) \lra \pi_1^{\rm Hod}(X,x)
\lra {\cal G}_{\rm Hod} \lra 0.  
\ee
A  base point $x\in X$ provides a  splitting 
$s_x: {\cal G}_{\rm Hod} \to \pi_1^{\rm Hod}(X,x)$, 
and hence  an action of ${\cal G}_{\rm Hod}$ on 
the Lie algebra $\pi_1^{\rm nil}(X,x)$, and thus a  
real MHS on $\pi_1^{\rm nil}(X,x)$.

It follows that the Lie algebra ${\cal G}_{\rm Hod}$ acts  
on representations of the Lie algebra  $\pi_1^{\rm nil}(X,x)$. The latter 
are automatically nilpotent, and  
form the formal neighborhood of the trivial local system. 
So ${\cal G}_{\rm Hod}$ acts on 
the formal neighborhood of the trivial local system. 
In Section \ref{hf6.5ref} we describe this action explicitly by using
 the class (\ref{CLass}). 





 \paragraph{Acknowledgments.} I am very grateful to 
Maxim  Kontsevich for useful discussions.

A part of this work 
 was done    
during my visits to the IHES at the Summer of 
2007 and MPI(Bonn) at the Summer of 2008. 
I am grateful to these institutions for the hospitality and support. 
I was supported by the NSF grant DMS-0653721. 
I am grateful to the referee for useful comments.


 \label{HF1sec}
 \section{Hodge correlators for compact Kahler manifolds} \la{hf2sec}

\subsection{The polydifferential operator $\omega$ } 
In this subsection we recall some key definitions from 
\cite{G1}. The proofs can be found in {\it loc. cit.}. 

Let $\varphi_0$, ..., 
$\varphi_{m}$ be smooth forms on a complex manifold $M$. 
We define a form $\omega(\varphi_0, ..., \varphi_{m})$ as follows. 
Given forms  $f_i$ of degrees ${\rm deg}(f_i)$ and a function $F(f_1, \ldots , f_m)$ 
set
\begin{equation} \label{5.25.06.1}
{\rm Sym}_{m}F(f_1, \ldots , f_m):= 
\sum_{\sigma \in \Sigma_m}{\rm sgn}_{\sigma; f_1, ..., f_m}F(f_{\sigma(1)}, \ldots , f_{\sigma(m)}), \qquad 
\end{equation}
where the sign of the transposition of $f_1, f_2$ is $ (-1)^{({\rm deg}(f_1)+1)
({\rm deg}(f_2)+1)}$,
 and the sign of a permutation 
written as a product of transpositions is 
the product of the signs of the corresponding transpositions. 
Observe that this is nothing else but the symmetrization 
of the elements $\overline f_i \in {\cal A}_M^*[-1]$ of 
the shifted de Rham complex of $M$ assigned to the elements $f_i$. 
Set
\begin{equation} \label{12.29.04.1}
\omega(\varphi_0, ..., \varphi_{m}) := 
\frac{1}{(m+1)!}{\rm Sym}_{m+1}
\Bigl( \sum_{k=0}^{m}(-1)^{k} \varphi_0 \wedge \partial \varphi_1 \wedge ... \wedge \partial \varphi_k \wedge 
\overline \partial \varphi_{k+1} \wedge ... \wedge \overline \partial \varphi_{m}\Bigr).
\end{equation} 
There is a degree zero linear map 
$$
\omega: {\rm Sym}^m \Bigl({\cal A}^{\bullet}_M[-1]\Bigr) \lra 
{\cal A}^{\bullet}_M[-1], 
\qquad 
\omega: \varphi_0 \wedge ... \wedge 
\varphi_{m} \lms \omega(\varphi_0, ... , \varphi_{m}).
$$
Its  key property is the following: 
\vskip 3mm
\begin{equation} \label{5-20.1a}
d \omega(\varphi_0, \ldots, \varphi_{m}) \quad = \quad 
(-1)^m \partial \varphi_0 \wedge \ldots \wedge \partial\varphi_{m} + 
\overline \partial \varphi_0 \wedge \ldots \wedge \overline \partial\varphi_{m} +
\end{equation}
\begin{equation} \label{5-20.1}
\frac{1}{m!}{\rm Sym}_{m+1}\Bigl((-1)^{{\rm deg}(\varphi_0)}
\overline\partial \partial \varphi_{0} \wedge 
\omega(\varphi_1, \ldots, \varphi_{m})\Bigr).
\end{equation}

\vskip 3mm
Set
$d^\C:= \partial - \overline \partial$. 
Consider the following differential forms:
\begin{equation} \label{5-x20.1s}
\xi(\varphi_0, \ldots, \varphi_{m}) :=  \frac{1}{(m+1)!}{\rm Sym}_{m+1}
\Bigl( 
\varphi_0 \wedge d^{\C}\varphi_{1}\wedge \ldots \wedge 
d^{\C}\varphi_{m}\Bigr)=
\end{equation}
$$
\frac{1}{(m+1)!}{\rm Sym}_{m+1}\Bigl( \sum_{k=0}^{m} (-1)^k  {m\choose k}
\varphi_0 \wedge \partial\varphi_{1}\wedge \ldots \wedge 
\partial\varphi_{k}\wedge \overline \partial\varphi_{k+1}\wedge \ldots \wedge\overline \partial\varphi_{m}\Bigr).
$$
\begin{equation} \label{5.15.06.7}
\eta(\varphi_0, \ldots, \varphi_{m}):= 
\frac{1}{(m+1)!}{\rm Sym}_{m+1}\Bigl( d^{\C}\varphi_0  
\wedge \ldots \wedge d^{\C}\varphi_m \Bigr) = 
\end{equation}
$$
\frac{1}{(m+1)!}
{\rm Sym}_{m+1}\Bigl( \sum_{k=0}^{m} (-1)^k {m+1\choose k+1}\partial \varphi_0  
\wedge \ldots \wedge\partial \varphi_k \wedge \overline \partial \varphi_{k+1}  
\wedge \ldots \wedge\overline \partial \varphi_m \Bigr).
$$
Every summand in the form $\eta$ appears with the 
coefficient $\pm 1$. We have 
\begin{equation} \label{5-20.1pa}
d^\C \xi(\varphi_0, \ldots, \varphi_{m})  =  \eta(\varphi_0, \ldots, \varphi_{m}).
\end{equation}

\subsection{Green currents} Let 
${M}$ be a compact  Kahler manifold. 
Denote by $({\cal A}_{M}^{*,*}, \partial, \overline \partial)$ 
the Dolbeaut bicomplex on ${M}$. 
Recall the Hodge decomposition $H^*({M},\C) = \oplus_{p,q}H^{p,q}({M})$. 
Choose a representative  
for each $(p.q)$-cohomology class of ${M}$. 
We call such a datum {\it a harmonic splitting of the Dolbeaut bicomplex}, 
and the representatives the {\it harmonic representatives}. 

\vskip 3mm
{\it Green current of a cycle}.  
Let $M$ be a smooth complex projective variety, and $Z \subset M$ 
a complex algebraic cycle of  codimension  $d_Z$ in $M$. 
Then integration over $Z$ is 
a $(d_Z, d_Z)$-current $\delta_Z$, representing 
the cohomology class of the cycle $Z$ in the Dolbeaut bicomplex 
${\cal A}_{M}^{*,*}$. 

On the other hand, given a choice of a splitting $s$ of the 
Dolbeaut bicomplex on $M$, 
we can represent the cohomology class of 
the cycle $Z$ by a smooth differential form ${\rm Har}_Z$ 
on $M$, its {harmonic representative}. 

\begin{definition} \label{5.17.07.3sd} Let $M$ be a smooth 
complex projective variety. A Green current
$G_Z$ corresponding to an algebraic cycle $Z\subset M$ 
and a splitting $s$ 
is a current 
on $M$ satisfying the equation
\begin{equation} \label{7.3.00.1mnsd}
(2\pi i)^{-1}\overline \partial \partial G_Z = 
\delta_{Z} - {\rm Har}_Z.
\end{equation}
\end{definition} 

Here on the right we have a difference of two $(d_Z, d_Z)$-currents 
whose cohomology classes 
provide the class of $Z$.  
Thus their difference is a $\partial $- and $\overline \partial $-closed 
$(d_Z, d_Z)$-current. 
By the $\overline \partial \partial$-lemma there exists 
a solution of the differential 
equation (\ref{7.3.00.1mnsd}). 
It is a $(d_Z-1, d_Z-1)$-current, well defined up to 
elements from  ${\rm Im}\partial + {\rm Im}\overline \partial$.  
There is a solution  
smooth outside of the cycle $Z$.

\vskip 3mm
{\it Green current of a map}. 
Let $f:M\to N$ be a map of compact Kahler manifolds. 
Choose splittings of the Dolbeaut complexes on $M$ and $N$. 
Denote by $\Delta_f$ 
the graph of the map. It is a cycle in $X \times Y$. 
Let ${\rm Har}_{\Delta_f}$ be the harmonic representative 
of the cohomology class of the cycle $\Delta_f$. 
Denote by $\delta_{\Delta_f}$  the $\delta$-current defined by the cycle 
$\Delta_f$. Then a Green current of the map $f$ is a current on $M \times N$ 
satisfying the equation 
\begin{equation} \label{7.3.00.112}
(2\pi i)^{-1}\overline \partial \partial G_Z = 
\delta_{\Delta_f} - {\rm Har}_{\Delta_f}.
\end{equation}

\vskip 3mm
{\it Green current of the diagonal}. 
Let 
$f$ be the identity map of $M$. 
Let $p_1$ and $p_2$ be the  projections of $M \times M$ onto the first and second factors. Let $\Delta_M$ be the diagonal in $M \times M$. 
Then $\delta_{\Delta_M}$ is an $(n, n)$-current on $M \times M$. 
 Choose a splitting of the Dolbeaut bicomplex on $M$.  
Let $\alpha_1, .., \alpha_N$ be the corresponding representatives of a 
Hodge basis in $\oplus_{0< p+q< 2n}H^{p,q}(M, \C)$. 
 Denote by $\alpha^{\vee}_1, .., \alpha^{\vee}_N$  the 
representatives for the 
dual basis: one has
$ 
\int_{M}\alpha_i\wedge \alpha^{\vee}_j = \delta_{ij}. 
$

The 
constant function $1$ is a canonical representative of a class in $H^0$. 
Let $\mu$ be a $2n$-current on $M$ such that $\int_{M}\mu =1$.  
A Green current  $G_{\mu}(x,y)$ is an $(n-1, n-1)$-current on 
$M \times M$ which satisfies 
\begin{equation} \label{7.3.00.1mn}
(2\pi i)^{-1}\overline \partial \partial G_\mu(x,y) = \delta_{\Delta_M} - 
\Bigl( p_1^*\mu + p_2^*\mu  + \sum_{k=1}^N 
 p_1^*\alpha^{\vee}_{k} \wedge p_2^*\alpha_{k}\Bigr).
\end{equation}
and is symmetric with respect to the 
permutation of the factors: 
$
G_{\mu}(x,y) = G_{\mu}(y, x).
$ 
On the right of (\ref{7.3.00.1mn}) stands the
 difference of two $(n,n)$-currents 
whose cohomology classes 
provide the identity map on $H^*(M, \C)$.  
This current is invariant under the 
permutation 
of the factors in $X \times X$.  So symmetrising 
any solution to  (\ref{7.3.00.1mn}) 
we get a symmetric solution. 

\vskip 3mm
There is a specific choice of $\mu$.  
A point $a \in M$ provides  a $2n$-current $\mu:= 
\delta_a$ on $M$: one has 
$\langle \delta_a, \varphi\rangle = \varphi(a)$. 
The 
corresponding Green current is denoted  $G_a(x,y)$. It satisfies the equation 
\begin{equation} \label{7.3.00.1d}
(2\pi i)^{-1}\overline \partial \partial G_a(x,y) = \delta_{\Delta_M} - 
\Bigl( \delta_{\{a\} \times M} + \delta_{M\times \{a\} } + \sum_{k}
 p_1^*\alpha^{\vee}_k \wedge p_2^*\alpha_k \Bigr).
\end{equation}

\vskip 3mm
\subsection{Hodge correlators} 
Let $X$ be an $n$-dimensional
 smooth compact complex Kahler variety. We start by introducing some notation. 
First, we use from now on the notation 
\be \la{8.9.08.1}
{\Bbb H} \quad \mbox{for} \quad  {\Bbb H}_*(X)
\ee
We  
allow the following abuse of notation: we denote sometimes by the same letter 
$\alpha$ a harmonic form and the corresponding element 
of ${\Bbb H}^*(X)$. Strictly speaking, the second element 
is  $\alpha\otimes {\bf E}$ where ${\bf E}$ is 
a degree $-1$ element, 
and thus should have a different notation, $\overline \alpha$. 
Similarly $h$ denotes, depending on the context, 
either a homology class, or a shifted by $1$ homology class. 
The latter  is also denoted 
by $\overline h$ to emphasize the shift 
of the degree. 

We denote by $|\overline \omega|:= {\rm deg}(\omega)-1$ the degree of a form 
$\omega$ in 
the shifted de Rham complex ${\cal A}^*(X)[1]$, and similarly 
$|\overline h|:= {\rm deg}(h)+1$ is the degree in $H_*(X)[-1]$. 
We reserve the notation ${\rm deg}(\omega)$ 
and ${\rm deg}(h)$ for their degrees before the shift. 

\vskip 3mm
A $2n$-current $\mu$ on $X$ is {\it admissible} if it is 
either a smooth volume form 
of total mass $1$ on $X$, 
or the $\delta$-current $\delta_a$ for some point $a \in X$.  
Given an admissible $2n$-current $\mu$,  we are going to define a linear map, 
the (precursor of) {Hodge correlator map}\footnote{We skip the subscript $\mu$ whenever possible.}:
\be \la{1.17.08.10}
{\rm Cor}_{{\cal H}, \mu}: {\cal C}_{{\Bbb H}^*} \otimes 
{\cal H}^{\vee} \lra \C.
\ee
Let us choose a harmonic splitting of the Dolbeaut bicomplex ${\cal A}^{*, *}(X)$. 
Given homogeneous harmonic forms   
$\alpha_i$, consider a  
cyclic element  
\be \la{1.17.08.11}
W = 
{\cal C}(\overline \alpha_0 \otimes \overline \alpha_1 \otimes \ldots \otimes 
\overline \alpha_m).
\ee 
\begin{figure}[ht]
\centerline{\epsfbox{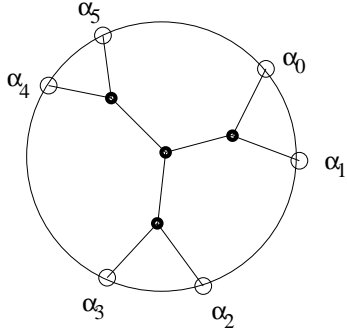}}
\caption{A plane  trivalent tree decorated by harmonic forms $\alpha_i$.}
\label{hf1}
\end{figure} 
  We use the standard terminology about trees and decorated trees, see  
Section 2 of \cite{G1}. 
Let $T$ be a plane trivalent tree decorated by the forms 
$\alpha_0, \ldots , \alpha_m$. 
The external vertices of a plane tree have a cyclic order 
provided by a chosen orientation of the plane, say clockwise. 
We assume that the cyclic order of the forms $\alpha_i$ is
 compatible with the cyclic order of the vertices, as on Fig \ref{hf1}. 
Although 
the decoration depends on the 
presentation of $W$ as a cyclic tensor product of the 
$\alpha_i$'s, 
constructions below depend only on $W$. 
We are going to assign to such a $W$-decorated tree $T$ a  
current $\kappa_W$ on 
\begin{equation} \label{7.1.00.1q} 
X^{\mbox{\{internal vertices of $T$\}}}. 
\end{equation}

Given a finite set ${\cal X} = \{x_1, \ldots, x_{|{\cal X}|}\}$,
 there is a $\Z/2\Z$-torsor, 
the {\it orientation torsor of 
${\cal X}$}. Its elements are expressions $x_1 \wedge 
\ldots \wedge x_{|{\cal X}|}$; 
interchanging two neighbors we change the sign of the expression. 
The {\it orientation torsor ${\rm or}_T$} of a graph $T$ is 
the orientation torsor of the set of edges of $T$. 
An orientation of the plane induces an orientation 
of a  plane trivalent tree.  

Take an internal edge $E$ of the tree $T$. There is a  projection 
\be \la{3.19.8.1}
p_E: X^{\mbox{\{internal vertices of $T$\}}} \quad \lra \quad 
X^{\mbox{\{vertices of E\}}} = X\times X.
\ee
Choose a 
Green current 
$G_\mu(x,y)$ corresponding to an admissible  $2n$-current $\mu$ on $X$.  
. Assign  to the edge $E$ the Green current $G_\mu(x,y)$ 
on the right space in (\ref{3.19.8.1}). 
Since the Green current is symmetric, the output is well defined. 
Denote by 
${G}_E$ its pull back by the map $p_E^* $. 
Since the map $p_E$ is transversal to the wave front 
of the Green current, the pull back  is well defined.

Let  $\{E_0, \ldots, E_{2m}\}$ be the set of  edges of $T$
numbered so that the internal edges are the first $m-2$ of them. Set $k=m-3$. 
Let us introduce degree $-1$ variables ${\bf E}_i$ 
matching the edges $E_i$. 
Let $E_{\alpha_i}$ be the external edge of the tree $T$ 
decorated by the form $\alpha_i$. 
Let $p_{\alpha_i}$ be the projection from (\ref{7.1.00.1q}) 
onto the factor corresponding to the internal vertex of $T$ assigned to the edge 
$E_{\alpha_i}$. 
Abusing notation, we denote below $p^*_{\alpha_{s}}\alpha_{s}$ by 
$\alpha_{s}$.

\vskip 3mm
Given homogeneous forms $\varphi_i$ assigned to the edges $E_i$, let us 
introduce an expression 
\be \la{7.17.08.1}
\omega \Bigl((\varphi_{{0}}\wedge {\bf E}_0) \wedge \ldots \wedge (\varphi_{{k}} 
\wedge {\bf E}_k)\Bigr)  = (-1)^{s} \omega \Bigl(\varphi_{{0}}\wedge \ldots \wedge 
\varphi_{{k}} 
\Bigr) \bigwedge {\bf E}_0 \wedge \ldots \wedge{\bf E}_k, 
\ee  
where $s= ({\rm deg}(\varphi_1) +1) + 2({\rm deg}(\varphi_2) +1) + ... +  
k({\rm deg}(\varphi_k) +1)$. 
Its left hand side  can be understood as follows: 
we apply 
the operators $\partial$ and $\overline \partial$ to the factors, e.g. 
$
\partial(\varphi_{m}  \wedge {\bf E}_{_{m}}):= \partial\varphi_{m}  \wedge {\bf E}_{m}
$ and  
get  
the sign $(-1)^s$ by moving ${\bf E}_i$ to the right. 

Observe that (\ref{7.17.08.1}) does not depend on the order 
of the edges  $E_{0}, \ldots , E_{k}$. 
Set  
\begin{equation} \label{3.8.05.15}
\widetilde \kappa_T(W):= 
\omega \Bigl((G_{E_{0}}\wedge {\bf E}_0) \wedge \ldots \wedge (G_{E_{k}} 
\wedge {\bf E}_k)\Bigr) \bigwedge 
(\alpha_{0}  \wedge {\bf E}_{\alpha_{0}})\wedge \ldots \wedge 
(\alpha_{m}  \wedge {\bf E}_{\alpha_{m}}).
\end{equation}
Here 
$
(\alpha_{0}  \wedge {\bf E}_{\alpha_{0}})\wedge \ldots \wedge 
(\alpha_{m}  \wedge {\bf E}_{\alpha_{m}}): = (-1)^t 
\alpha_{0}  \wedge \ldots \wedge  
\alpha_{m} \bigwedge  {\bf E}_{\alpha_{0}}\wedge \ldots \wedge{\bf E}_{\alpha_{m}},
$ 
 where the sign is obtained by moving ${\bf E}$'s through $\alpha$'s. 
It is invariant with respect to the cyclic shift
of $W$.  
So it is determined by the cyclic word $W$. 
Therefore $\widetilde \kappa_T(W)$ is also 
determined by the cyclic word $W$. 

Moving the odd variables ${\bf E}_{s}$ in (\ref{3.8.05.15}) to the right we get 
\begin{equation} \label{3.8.05.15.a}
\widetilde \kappa_T(W) = \pm \omega(G_{E_{0}}\wedge \ldots \wedge G_{E_{k}} ) 
\bigwedge 
\alpha_{0}  \wedge \ldots \wedge 
\alpha_{m}  \bigwedge  {\bf E}_{0} \wedge \ldots \wedge {\bf E}_{2m}.  
\end{equation}
Now we use an orientation ${\rm Or}_T$ of the tree $T$ 
to get rid of the factor ${\bf E}_{0} \wedge \ldots \wedge {\bf E}_{2m}$, setting 
\begin{equation} \label{3.8.05.15.ab}
\kappa_T(W):= \pm {\rm sgn}(E_0 \wedge \ldots \wedge E_{2m}) 
\wedge \omega(G_{E_{0}}\wedge \ldots \wedge G_{E_{k}} ) 
\bigwedge 
\alpha_{0}  \wedge \ldots \wedge 
\alpha_{m}
\end{equation}
where $\pm$ is the same sign as in (\ref{3.8.05.15.a}), and 
$
{\rm sgn}(E_0 \wedge \ldots \wedge E_{2m}) \in \{\pm 1\}
$  
is the difference between the 
element $E_0 \wedge \ldots \wedge E_{2m}\in {\rm or}_T$ and the 
generator ${\rm Or}_T$ corresponding to the  
clockwise orientation of the plane. 

{\bf Remark}. The form $\omega(G_{E_{0}}\wedge \ldots \wedge G_{E_{k}} ) 
\bigwedge 
\alpha_{0}  \wedge \ldots \wedge 
\alpha_{m}$ is of even degree. So 
one can move the factor 
${\bf E}_{0} \wedge \ldots \wedge {\bf E}_{2m}$ in (\ref{3.8.05.15.a}) to the left 
without getting an extra sign. 

\begin{lemma}
$\kappa_T(W)$ is a current on (\ref{7.1.00.1q}).
\end{lemma}

{\bf Proof}. A Green current $G_a(x,y)$ has the same type of singularity 
at the diagonal and at the 
cycles $x=a$ and $y=a$. These singularities are integrable since 
by the general theory of elliptic equations a solution of 
the equation (\ref{7.3.00.1d}) exists as a current. 

More specifically, for the current $\partial G_a(0,z)$ the singularity at $z=0$ 
is described by the singularity of the Bochner-Martinelly kernel:
$$
\partial G(0, z) \sim \frac{1}{(2\pi i)^n}\frac{\alpha_n(\overline z, 
d\overline z)\wedge dz_1 \wedge \ldots \wedge dz_n}
{(|z_1|^2 + \ldots + |z_1|^2)^n}, \qquad \alpha_n(\overline z, d\overline z) = 
{\rm Alt}_n(\overline z_1 d\overline z_2 \wedge 
\ldots \wedge d\overline z_n).
$$
Here 
$z_1, \ldots , z_n$ are local coordinates near $z=0$. 
Writing it in the spherical coordinates, one sees that this form, 
multiplied by a smooth $1$-form, has an integrable singularity. Similar argument 
is applied to $\overline \partial G(0, z)$, and the form  $G(0, z)$ is even less singular. 
This is another way to see that 
 the form $\kappa_T(W)$ is integrable  at the generic points of the diagonals $x=y$. 

Since diagonals in (\ref{7.1.00.1q}) intersect transversally, 
the form $\kappa_T(W)$ is integrable on intersections of the diagonals. 
In particular, it is integrable at subvarieties  of the diagonals given by 
$x=a=y$.

Although  the singularity at $x=a$ is integrable, it causes 
a problem since $\kappa_T(W)$ contains 
products of several Green forms with the same singularity at $x=a$. 
In fact a single product like $G_a(s_1, x)\wedge 
\partial G_a(s_2, x)\wedge \overline \partial G_a(s_3, x)$ does have a singularity at $x=a$. 
We claim that nevertheless the non-integrable singularity disappears after the skew-symmetrization. 
Indeed, setting $a=0$, one can write a Green current in local coordinates near $x=0$ 
as $G_0(x,y) = B(x) + S(x,y)$ where 
$B(x)$ is the Bochner-Martinelly type current and $S(x,y)$ is smooth. 
Applying the skew-symmetric polydifferential operator $\omega$ we get 
$$
\omega\Bigl(B(x) + S(x,y_1)) \wedge (B(x) + S(x,y_2)) \wedge (B(x) + S(x,y_3)\Bigr) = 
$$
$$
{\rm Alt}_{y_1, y_2, y_3}\omega\Bigl(B(x) \wedge S(x,y_2) \wedge S(x,y_3)\Bigr)
  + \mbox{\rm a smooth form}. 
$$
So the claim follows from the integrabilty of the 
Green current $G_a(x,y)$. 
The lemma is proved. 

\vskip 3mm

\begin{definition} 
Let 
$p_T: X^{\{\mbox{internal vertices of T}\}} \lra \mbox{point}$ 
be the natural projection. 
The  number  ${\rm Cor}_{{\cal H}}(W\otimes {\cal H})$ is 
given by 
$$
{\rm Cor}_{{\cal H}}(W\otimes {\cal H}):= {p_T}_*\Bigl(\sum_{T}\kappa_T(W)\Bigr),
$$
where the sum is over all plane trivalent trees $T$ 
decorated by $W$. 
\end{definition}

\begin{lemma} \label{5.17.07.5}
The Hodge correlator map (\ref{1.17.08.10}) is a linear map of degree zero, 
provided the number of factors in the cyclic 
tensor product is bigger then $3$.  
\end{lemma}

{\bf Proof}. 
A trivalent graph $T$ with $m+1$ external legs has 
$m-1$ internal vertices and $m-2$ internal edges. Let  
$W$ is as in (\ref{1.17.08.11}). 
Then we integrate a form of degree 
$(2n-2)(m-2) +(m-3) + \sum_{i=0}^{m} {\rm deg}(\alpha_i)$ 
over a cycle of dimension $2n(m-1)$. 
The result can be non-zero only if
$
(2n-2)(m-2) +(m-3) + \sum_{i=0}^{m} {\rm deg}(\alpha_i)= 2n(m-1), 
$ 
i.e. $\sum_{i=0}^{m} ({\rm deg}(\alpha_i)-1) = (2n-2)$. 
The left hand side here equals  ${\rm deg}W$. 
The lemma follows.

\begin{proposition} \label{7.3.06.4} Let $\overline \alpha_i \in {\Bbb H}^*_{X}$. 
Then for any $p,q \geq 1$ one has 
\begin{equation} \label{7.3.06.5}
\sum_{\sigma \in \Sigma_{p,q}}\pm 
{\rm Cor}_{\cal H}(\overline \alpha_0 \otimes \overline \alpha_{\sigma(1)} 
\otimes \ldots \otimes \overline \alpha_{\sigma(p+q)}) = 0, 
\end{equation}
where the sum is over all $(p,q)$-shuffles $\sigma \in \Sigma_{p,q}$, and 
the signs are the standard ones. 
\end{proposition}

The proof is identical to the proof of 
the shuffle relation in Section 2 of \cite{G1}.

\vskip 3mm
\paragraph{The Hodge correlator.} 
Let  
${\cal C}$ be the operator of weighted projection 
on the cyclic tensor algebra: 
\be \la{WCP}
{\cal C}(x_1 \otimes \ldots \otimes x_m):= \frac{1}{|{\rm Aut}(W)|}
 (x_{1} \otimes \ldots 
\otimes x_{m})_{\cal C}.
\ee
where $W= (x_{1} \otimes \ldots 
\otimes x_{m})_{\cal C}$, 
and $|{\rm Aut}(W)|$ is the order of its automorphism group. 

\vskip 3mm
Given a collection of harmonic differential 
forms $\omega_i$ and homology classes $h_i$ we define the 
Hodge correlator map of the 
cyclic tensor product of the expressions $\omega_i | h_i$ as a cyclic 
product of their homology factors 
taken with the coefficient 
given by the Hodge correlator map applied to the cyclic product 
of the corresponding  differential form factors, with the
 sign computed via the standard sign rule. Here is an example.
\be \la{HCC}
{\rm Cor}_{\cal H}(\overline \omega_1 | \overline h_1 \otimes \overline \omega_2 | \overline h_2 \otimes 
\overline \omega_3 | \overline h_3 ) :=
\ee
$$
  (-1)^{|\overline \omega_3| (|\overline h_2| + |\overline h_1|) + |\overline h_1||\overline \omega_2|}
{\rm Cor}_{\cal H}(\overline \omega_1 \otimes 
\overline \omega_2\otimes \overline \omega_3)
\cdot ( \overline h_{1}\otimes \overline h_{2} \otimes \overline h_{3})_{\cal C} = 
$$
$$
(-1)^{|\overline \omega_3| (|\overline h_2| + |\overline h_1|) + |\overline h_1||\overline \omega_2|}
(-1)^{{\rm deg}\omega_2}\int_X (\omega_1 \wedge 
\omega_2\wedge \omega_3)
\cdot ( \overline h_{1}\otimes \overline h_{2} \otimes \overline h_{3})_{\cal C}.
$$
Since ${\rm Cor}_{\cal H}(\overline \omega_1 \otimes \ldots \otimes \overline \omega_m)$ 
is graded cyclic invariant, the same is true for  
${\rm Cor}_{\cal H}(\overline \omega_1|\overline h_1 
\otimes \ldots \otimes \overline \omega_m|\overline h_m)$. 
In particular if ${\rm deg}(\omega_i) = {\rm deg}(h_i)$, then 
it is plain cyclic invariant.  

\vskip 3mm
Since $\overline \alpha|h_{\overline \alpha}$ is even, 
properties of the Hodge correlators written 
in this form are more transparent. For example, the shuffle relations look as follows 
(all signs are pluses): 
\begin{equation} \label{7.3.06.5a}
\sum_{\sigma \in \Sigma_{p,q}} 
{\rm Cor}_{\cal H}\Bigl(\overline \alpha_0| \overline h_{\alpha_0}
\otimes \overline \alpha_{\sigma(1)}| \overline h_{\alpha_\sigma(1)} 
\otimes \ldots \otimes \overline \alpha_{\sigma(p+q)}| \overline h_{\alpha_{\sigma(p+q)}}\Bigr)_{\cal C} = 0. 
\end{equation}

\vskip 3mm
 Dualising the map (\ref{1.17.08.10}) we get an element  
${\bf G} \in  {\cal C}_{\Bbb H} \otimes {\cal H}$, 
the precursor of  the Hodge correlator. 
We write it as follows. 
Choose a basis in the 
reduced cohomology ${\Bbb H}^*$ of  $X$. Denote by $\alpha_i$ 
the corresponding harmonic representatives. 
Let $h_{\alpha_i}$ be the dual
 basis of the reduced homology ${\Bbb H}$.

Consider the Casimir element
\be \la{CASx}
{\rm Id}_X:=  1 + \alpha|h_{\alpha} + \alpha|h\otimes \alpha|h
+ \alpha|h\otimes \alpha|h \otimes\alpha|h + \ldots, \qquad
\alpha|h:= \sum_s\overline \alpha_s |\overline h_{\alpha_s},
\ee
where  the sum is over a basis $\{\alpha_s\}$ of harmonic forms. 
Applying operator (\ref{WCP}), we get the cyclic Casimir element:
\be \la{CYCCAS}
{\cal C}
{\rm Id}_X:= {\cal C}\Bigl( 1 + \alpha|h + \alpha|h\otimes \alpha|h 
+ \alpha|h\otimes \alpha|h \otimes\alpha|h + \ldots \Bigr).
\ee
Then one has, where  ${\cal H}$ is the shifted fundamental cohomology 
class (\ref{sfc}), 
\be \la{GGGG}
{\bf G} = {\rm Cor}_{\cal H}\Bigl( {\cal C}{\rm Id}_X \Bigr)\otimes {\cal H}
\in  {\cal C}_{\Bbb H} \otimes {\cal H}.
\ee

We elaborate this as follows. Let   $W = {\cal C}(\overline \alpha_1 \otimes 
\ldots \otimes \overline \alpha_m)$ be a basis in ${\cal C}_{\Bbb H}$. Then 
\be \la{4.25.08.1}
{\bf G} = 
\sum_{W}\frac{1}{|{\rm Aut}(W)|}
{\rm Cor}_{{\cal H}, \mu}(\overline \alpha_1|\overline h_{\alpha_1} 
\otimes \ldots \otimes \overline \alpha_m| 
\overline h_{\alpha_m}) \otimes {\cal H}. 
\ee

 \label{HF2sec}
\section{Motivic rational homotopy type DG Lie algebras} \la{hf3sec}

\subsection{Motivic set-up}

To describe the known 
realizations of the motivic DG Lie algebra ${\rm L}_{\cal M}(X, x)$ 
 let us recall  
the corresponding set-ups (see \cite{G1} for more details).

\vskip 3mm
Let $F$ be a field. We  
work in one of the following categories  
${\cal C}$:

\vskip 3mm 
i)  {\it Motivic}. The hypothetical abelian category of mixed 
motives over a field $F$.

ii) {\it Hodge}. $F=\C$, and  ${\cal C}$ is the category of mixed 
$\Q$- or $\R$-Hodge structures. 

iii) {\it Mixed $l$-adic}. 
$F$ is an arbitrary field such that $\mu_{l^{\infty}} \not \in F$, and 
 ${\cal C}$ 
is 
 the mixed category of 
$l$-adic ${\rm Gal}(\overline F/F)$--modules with a filtration 
$W_{\bullet}$ indexed by integers,  
such that ${\rm gr}^W_n$ is a pure of weight $n$.

iv) {\it Motivic Tate}. $F$ is a number field, 
${\cal C}$ is the abelian category of mixed 
Tate motives over $F$, equipped with the Hodge and $l$-adic 
realization functors, c.f. \cite{DG}. 
 
\vskip 3mm 

The setup i) is conjectural. The other three are well defined.  
A   category ${\cal C}$ from the list above  is an $L$--category, 
where $L = \Q$ in i), iv); $L=\Q$ or $\R$ in ii);  and $L=\Q_l$ in iii). 

Each of the categories has an invertible object, the {\it Tate object}, 
which is denoted, abusing notation, by  
$\Q(1)$ in all cases. 
 We set $\Q(n):= \Q(1)^{\otimes n}$. 
Each object in ${\cal C}$ 
 carries a canonical weight filtration 
$W_{\bullet}$,  morphisms in ${\cal C}$ are strictly 
compatible with this filtration. The weight of $\Q(1)$ is $-2$.

\vskip 3mm

One should have a ${\cal C}$-motivic 
rational homotopy type DG Lie algebra ${\rm L}_{\cal C}(X, x)$ for
 each of the categories ${\cal C}$. It is a DG Lie algebra 
in the category ${\cal C}$, equipped with a 
weight filtration $W_{\bullet}$.

Denote by ${\rm L}^{(l)}(X, x)$ 
the $l$-adic realization 
of the DG Lie algebra ${\rm L}_{\cal M}(X, x)$. The Lie algebra 
of the Galois group 
${\rm Gal}(\overline F/F)$ acts by  
derivations 
of the graded Lie algebra ${\rm L}^{(l)}(X, x)$ 
commuting with the differential and defined modulo homotopy.  
We give a linear algebra description 
of  the Lie algebra of all such derivations of the graded Lie algebra 
${\rm L}^{(l)}(X, x)$. 


\vskip 3mm
{\bf Example: $X$ is a smooth projective curve}. 
 The associated graded 
${\rm gr}^W\pi^{\rm nil}_1(X, x)$ of the 
pronilpotent completion $\pi^{\rm nil}_1(X, x)$ of $\pi_1(X, x)$ 
is canonically 
isomorphic to the quotient of the 
free, graded  by the  weight, Lie algebra  
generated by  $H_1(X)$ of weight $-1$, by the ideal generated by 
the element ${\rm S} = \sum_i[q_i, p_i]$, where $(q_i, p_i)$ 
is a symplectic basis of $H_1(X)$ for the intersection form. 
The  Lie algebra $\pi^{\rm nil}_1(X, x)$  is 
non-canonically isomorphic  to ${\rm gr}^W\pi^{\rm nil}_1(X, x)$. 

On the other hand, in this case 
$$
{\Bbb H} = H_2(X)[-1] \oplus H_1(X)[-1].
$$ 
The DG Lie algebra ${\rm L}_{\cal M}(X, x)$ is a free graded Lie algebra 
generated by ${\Bbb H}$, 
equipped  with the differential 
$$
\delta: H_2(X)[-1] \lra \Lambda^2(H_1(X)[-1]) 
$$
dualising to the product map. 
Since $\delta$ is injective, 
the DG Lie algebra ${\rm L}_{\cal M}(X, x)$ is quasiisomorphic to the  
Lie algebra generated by 
 $H_1(X)$, sitting in degree $0$, with a single relation 
 $$
\delta H_2(X)[-1] = \sum_i[q_i, p_i]=0.
$$

\vskip 3mm
\subsection{Cyclic words and  derivations} 
 In this subsections we present a version of the 
non-commutative symplectic geometry of M. Kontsevich  \cite{K} 
for an algebra with a Poincare duality. 

Let ${\rm H}^*(X)$ be the cohomology motive of an $n$-dimensional 
 regular projective variety $X$. 
It is a pure graded commutative algebra with 
a Poincare duality: 
${\rm H}^0(X)=\Q(0)$, the trace map provides an 
isomorphism $
{\rm H}^{2n}(X) \stackrel{\sim}{\lra} \Q(-n)[-2n]$, 
and there is a perfect pairing
$$
\cup: {\rm H}^*(X) \otimes {\rm H}^*(X) \lra {\rm H}^{2n}(X)
$$
given as product followed by the projection
${\rm H}^*(X) \lra {\rm H}^{2n}(X)$. 
Recall the fundamental class of $X$ shifted by $2$: 
$$
{\cal H}:= {\rm H}^{2n}(X)[2].
$$
It is isomorphic to $\Q(-n)[-(2n-2)]$ via the trace map. 
Thanks to the shift by $2$, the Poincare duality 
is a perfect pairing on ${\rm H}^*(X)[1]$ with values in  ${\cal H}$, and 
there is a $\cap$-product isomorphism
$$ 
\cap: {\rm H}_*(X)[-1] \otimes {\cal H} \stackrel{\sim}{\lra} {\rm H}^*(X)[1].
$$ 
It induces the reduced $\cap$-product  isomorphism 
\begin{equation} \label{polaris1*}
\cap: {\Bbb H}\otimes {\cal H} \stackrel{\sim}{\longrightarrow}  {\Bbb H}^*.
\end{equation}
The canonical pairing 
${\rm H}_*(X)\otimes {\rm H}^*(X) \to \Q(0)$ 
 induces a pairing $\langle *,*\rangle :  
{\Bbb H}\otimes {\Bbb H}^* \to \Q(0)$. 
Set 
\be \la{3.27.08.1}
\langle \ast \cap{\cal H} \cap \ast \rangle: 
{\Bbb H}\otimes {\cal H}\otimes {\Bbb H} 
\lra \Q(0), \qquad 
a \otimes  {\cal H} \otimes  b \to
\langle a, {\cal H} \cap b\rangle, \qquad a, b \in  {\Bbb H}.
\ee 
Since 
${\cal H}$ is  even, we can move it freely. The form (\ref{3.27.08.1}) 
is a symplectic form on ${\rm H}_*(X)[-1]$:
$$
\langle a \cap{\cal H} \cap b \rangle = -(-1)^{|a||b|}
\langle b \cap{\cal H} \cap a \rangle. 
$$ 
Indeed, $X$ is even dimensional, and 
it comes from a symmetric bilinear form on ${\rm H}_*(X)$.

\vskip 3mm
{\it The cyclic tensor product}. 
Recall the cyclic tensor product ${\cal C}_{{\Bbb H}}$ of ${\Bbb H}$. 
It is a sum of pure objects. 
We denote by  $(v_0 \otimes ... \otimes v_m)_{\cal C}$ 
the cyclic tensor product 
of direct summands $v_i$ of ${\Bbb H}$. 

\vskip 3mm
{\it The non-commutative differential}. 
It is a map
$$
{\Bbb D}: {\cal C}_{{\Bbb H}} 
\lra {\rm T}_{{\Bbb H}}\otimes {\Bbb H}, \qquad 
{\Bbb D}(h_0 \otimes ... \otimes h_m)_{\cal C}:= {\rm Cyc}_{m+1}
\left((h_{0}\otimes ... \otimes h_{m-1} )\otimes h_m\right).
$$ 
where ${\rm Cyc}_{m+1}$ is the 
the operator of the graded cyclic shift sum. 
Given  a decomposition ${\Bbb H} = \oplus a_i$ into simple objects,  
we define $\partial_{a_i} F \in {\rm T}_{{\Bbb H}}$ 
via the decomposition 
$$
{\Bbb D}F = \sum_i\partial_{a_i} F \otimes a_i.
$$

We are going to define an action of ${\cal C}_{{\Bbb H}}\otimes {\cal H}$  
by derivations 
of the DGA $({\rm T}_{{\Bbb H}}, \delta)$, and  introduce 
a  DG Lie algebra structure on 
${\cal C}_{{\Bbb H}}\otimes {\cal H} $, making it into an action of 
a DG Lie algebra on a DGA.

\vskip 3mm
{\it An action of 
${\cal C}_{{\Bbb H}}\otimes {\cal H}$ by derivations of  
${\rm T}_{{\Bbb H}}$}.  
We are going to define a map
$$
\theta: {\cal C}_{{\Bbb H}}\otimes {\cal H}\lra {\rm Der}({\rm T}_{{\Bbb H}}). 
$$
Since the algebra ${\rm T}_{{\Bbb H}}$ is free, to define 
its derivation we define it on generators. 
Given  an element $F\otimes {\cal H}
 \in {\cal C}_{{\Bbb H}}\otimes {\cal H}$, 
the derivation $\theta_{F\otimes {\cal H}}$ acts on 
$q \in {\Bbb H}$ by 
$$
\theta_{F\otimes {\cal H}}: q \lms  
\sum_{p}\partial_p F\otimes 
\langle  p \cap {\cal H}\cap q\rangle \in {\rm T}_{{\Bbb H}}.
$$

Here the sum is over a basis $\{p\}$ in ${\Bbb H}$. 
Clearly  map $\theta$ is a degree zero.

\vskip 3mm
{\it A Lie bracket on 
${\cal C}_{{\Bbb H}}\otimes {\cal H}$}. 
Let us define a Lie bracket
$$
\{*,*\}: ({\cal C}_{{\Bbb H}}\otimes {\cal H}) \otimes 
({\cal C}_{{\Bbb H}}\otimes {\cal H}) \lra {\cal C}_{{\Bbb H}}
\otimes {\cal H},
$$
It is handy to have a ``left'' version of the differential, 
which differs from ${\Bbb D}F$ only by signs:
$$
{\Bbb D}^-: {\cal C}_{{\Bbb H}} 
\lra {\Bbb H}\otimes{\rm T}_{{\Bbb H}}, \qquad 
(h_0 \otimes ... \otimes h_m)_{\cal C}:= {\rm Cyc}_{m+1}
\left(h_0\otimes (h_{1}\otimes ... \otimes h_{m} )\right).
$$
It gives rise to ``left'' partial derivatives 
$\partial^-_q G$. Set 
\begin{equation} \label{6.20.00.11cas}
\{F\otimes {\cal H}, G\otimes {\cal H}\}:=
\sum_{p,q}\Bigl(\partial_p F 
 \otimes
\langle p \cap{\cal H} \cap q\rangle 
\otimes \partial_q^- G \Bigr)_{\cal C}
\otimes {\cal H}. 
\end{equation}
The sum is  over bases $\{p\}$, $\{q\}$ of 
${\Bbb H}$. 
Here is an invariant definition of the bracket $\{*,*\}$. Take 
$$
{\Bbb D}F\otimes {\cal H} \otimes 
{\Bbb D}^-G\otimes {\cal H} \subset {\rm T}_{{\Bbb H}} 
\otimes  {\Bbb H}\otimes {\cal H}\otimes {\Bbb H}\otimes 
{\rm T}_{ {\Bbb H}}\otimes {\cal H}. 
$$ 
Then we apply map (\ref{3.27.08.1})
and the product ${\rm T}_{{\Bbb H}} \otimes \Q(0) \otimes {\rm T}_{{\Bbb H}}\to 
{\rm T}_{{\Bbb H}}$, followed by  projection to 
${\cal C}_{{\Bbb H}}$. 
 Similarly the map $\theta_{F \otimes {\cal H}}$ is given as the composition 
${\Bbb D}F\otimes {\cal H} \otimes {\Bbb H} \hra {\rm T}_{{\Bbb H}} 
\otimes  {\Bbb H}\otimes {\cal H}\otimes {\Bbb H} \to {\rm T}_{{\Bbb H}}.$

\bp \la{3.28.08.4}
The bracket $\{*,*\}$ is a degree zero map. It is skew-symmetric:
\be \la{3.28.08.2}
\{F\otimes {\cal H}, G\otimes {\cal H}\} = 
-(-1)^{|F||G|}\{G\otimes {\cal H}, F\otimes {\cal H}\}. 
\ee
It satisfies the Jacobi identity. it is dual to the Lie cobracket 
defined on Fig. \ref{hf7}. 
\ep

{\bf Proof}. The first claim is evident. 
The second is checked using the fact that $\langle p \cap{\cal H} \cap q\rangle$ 
is a symplectic form. The last is easy to prove directly. 
We will deduce it from the following 

\begin{lemma} \label{3/11/07/10} 
The map $\theta: {\cal C}_{{\Bbb H}}\otimes {\cal H}
 \lra {\rm Der}({\rm T}_{{\Bbb H}})$ is injective and respects the brackets.
\end{lemma}
\begin{figure}[ht]
\centerline{\epsfbox{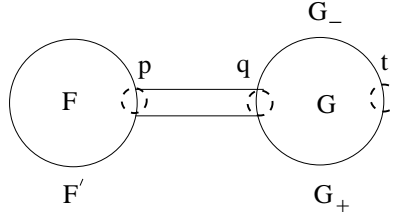}}
\caption{Calculating 
$\theta_{F\otimes{\cal H}} \circ \theta_{G\otimes{\cal H}}$.}
\label{hf4}
\end{figure}
{\bf Proof}. 
Let us calculate how  the super commutator 
$[\theta_{F\otimes{\cal H}}, \theta_{G\otimes{\cal H}}]$ acts on 
$h \subset {\Bbb H}$. Given certain factors $p$ of $F$ and $q, t$ of $G$, 
write  (see Fig \ref{hf4})
$F = (F' \otimes p)_{\cal C}$, and $G = (G_+qG_-t)_{\cal C}$.  
Then 
$$
\theta_{F\otimes{\cal H}} \circ \theta_{G\otimes{\cal H}}: 
h \lms \sum_{t}
(-1)^{|G_+||qG_-t|}\theta_{F\otimes{\cal H}} (G_+qG_-)  
\langle  t\cap {\cal H}\cap h\rangle= 
$$
$$
\sum_{ t, p, q}(-1)^{|G_+||G_-t| +|F||G_+|}
G_+F' \langle p \cap {\cal H}\cap q\rangle 
G_- 
\langle t \cap {\cal H} \cap h\rangle.
$$
Notice that the sign is given by moving $G_+$ through the remaining factor. 
We claim that this equals  $\theta_{\{F\otimes{\cal H}, G\otimes{\cal H}\}}(h)$.  
Indeed, $\theta$ for the element 
$(F'\langle p \cap {\cal H}\cap q\rangle G_-tG_+)_{\cal C}$
maps $h$ to 
$$
(-1)^{|G_+||Ft|}\sum_{t} G_+F'  \langle p \cap {\cal H}\cap q\rangle 
G_-
\langle  t\cap {\cal H}\cap h\rangle.
$$
The Lemma follows from this. Hence the Proposition is proved.

\begin{corollary} \la{3.28.08.5}
The map $\theta$ 
is an injective morphism of graded Lie algebras
\begin{equation} \label{3/11/07/11} 
\theta: {\cal C}_{{\Bbb H}}\otimes {\cal H}
 \lra {\rm Der}({\rm T}_{{\Bbb H}}).
\end{equation} 
\end{corollary}

\vskip 3mm
 {\it The canonical element $\Delta$}. 
Take a plane trivalent tree $T_3$ with a single internal vertex. 
\begin{figure}[ht]
\centerline{\epsfbox{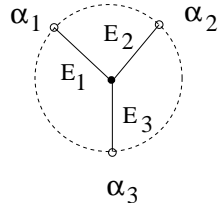}}
\caption{The canonical element $\Delta$.}
\label{hf3}
\end{figure}
Denote by ${\bf T}$ the orientation super line of 
this tree. It is a one dimensional space in the degree $-3$. 
The clockwise orientation of the plane provides a generator  
$
{\bf E}_1\wedge {\bf E}_2\wedge {\bf E}_3 \in {\bf T}. 
$ 
Here ${\bf E}_1, {\bf E}_2, {\bf E}_3$ are the degree $-1$ generators 
assigned to the edges ${E}_1, {E}_2, {E}_3$ of the tree. 
There is a map
 $$
 \otimes^3({H}^*(X)[1]) \lra H^{2n}(X) 
\otimes {\bf T} \stackrel{\sim}{=}{\cal H}[1].
 $$ 
 Namely, we assign the shifted by $1$ cohomology classes 
$\alpha_i\otimes {\bf E}_i$ to the 
external vertices of the edges $E_i$ of the tree $T_3$,  
and set 
$$
(\alpha_1 \otimes {\bf E}_1)\wedge (\alpha_2 \otimes {\bf E}_2)\wedge 
(\alpha_3 \otimes {\bf E}_3) \lms (-1)^{{\rm deg}(\alpha_2)}
(\alpha_1 \wedge \alpha_2 \wedge \alpha_3)_{H^{2n}(X)} \otimes 
{\bf E}_1\wedge {\bf E}_2\wedge {\bf E}_3. 
$$
Here $(\alpha)_{H^{2n}(X)}$ denotes the projection of a cohomology class 
$\alpha$ onto $H^{2n}(X)$. 
This map is graded cyclic invariant. 
 It gives rise to a map 
 $\Q(0)\lra \otimes^3({H}_*(X)[-1])\otimes {\cal H}[1]$.  
 Killing $H_0(X) \oplus H_{2n}(X)$ and projecting to the 
 cyclic envelope we obtain a canonical injective  map 
 $$
 \Delta: \Q(0)[-1]\lra {\cal C}({\Bbb H} 
 \otimes {\Bbb H} \otimes {\Bbb H})\otimes {\cal H}
  \subset {\cal C}_{{\Bbb H}}\otimes {\cal H}. 
 $$ 
 We often identify it with its image, also denoted by $\Delta$. 
 
\begin{proposition} \label{3/11/07/4}
 (i) One has $\{\Delta, \Delta\}=0.$

(ii) The map $\delta: \ast \lms \{\Delta, \ast\}$ given by 
the   commutator with $\Delta$   
is a derivation: 
$\delta^2=0$. 
 
(iii) The bracket $\{*,*\}$ and the differential 
$\delta$ provide ${\cal C}_{{\Bbb H}}\otimes {\cal H}$ with  
a DG Lie algebra structure. 
\end{proposition}

 {\bf Proof}. 
 Since ${\rm deg}~\Delta =1$, the map $\delta$ is of degree $1$. 
Now (i) is equivalent to (ii) thanks to the Jacobi identity. 
(iii) obviously follows from (ii). 
 The part (i) is easy to check directly. 
The part (ii) is also follows from the explicit formula 
for the map $\delta$ given in Lemma \ref{3.28.08.10} below (whose proof is independent of the proposition). 
The proposition is proved.

\vskip 3mm
There is a differential 
$\delta:= \theta_{\Delta}$ on ${\rm T}_{{\Bbb H}}$. 
It obviously provides ${\rm T}_{{\Bbb H}}$ with a DGA structure. The  differential 
on ${\cal C}_{{\Bbb H}}\otimes {\cal H}$ coincides with the one 
inherited from the differential  
on ${\rm T}_{{\Bbb H}}$. 
\vskip 3mm
\begin{lemma} \label{3.28.08.10}
The map $\theta$ is an action of the DG Lie algebra 
${\cal C}_{{\Bbb H}}\otimes {\cal H}$ on the DGA 
${\rm T}_{{\Bbb H}}$. 
\end{lemma}

{\bf Proof}. We have to show that 
${\cal C}_{{\Bbb H}}\otimes {\cal H}\otimes {\rm T}_{{\Bbb H}} 
\lra {\rm T}_{{\Bbb H}}$ commutes with $\delta$. 
This is very similar to, and in fact deduced from 
the fact that $\{\Delta, \{F\otimes {\cal H}, G\otimes {\cal H}\}\}$ 
satisfies Jacobi identity.

\begin{corollary} \label{3.28.08.10}
If $\delta(F \otimes {\cal H})=0$, then 
 $\theta_{F \otimes {\cal H}}$ is a derivation 
of the graded  algebra 
${\rm T}_{{\Bbb H}}$ which commutes 
with the differential $\delta$. 
\end{corollary}

\vskip 3mm
The DGA algebra ${\rm T}_{{\Bbb H}}$ 
has a DG Hopf algebra structure with the standard (coshuffle) 
coproduct $\nu$. The kernel of  $\nu - 
(1 \otimes {\rm Id} + {\rm Id} \otimes 1)$ is a 
DG Lie algebra ${\cal L}ie_{{\Bbb H}}$. It is  the free graded Lie algebra generated 
by ${\Bbb H}$. 
The differential 
$\delta$ on ${\rm T}_{{\Bbb H}}$ induces a differential on 
${\cal L}ie_{{\Bbb H}}$. 

Recall the space ${\cal C}{\cal L}ie_{\Bbb H} \otimes {\cal H}$ defined in 
Section 2.2.
\bl \la{popo}
${\cal C}{\cal L}ie_{\Bbb H} \otimes {\cal H}$ is  
the subspace of the DG Lie algebra ${\cal C}_{\Bbb H}\otimes {\cal H}$ 
preserving ${\cal L}ie_{{\Bbb H}}$. It is a DG Lie subalgebra. 
\el

{\bf Proof}. An easy exercise, left to the reader.

\bl \la{6.15.08.1} There is an isomorphism of DG Lie algebras
\be \la{6.15.08.2}
{\cal C}{{\cal L}ie}_{\Bbb H}\otimes {\cal H} \stackrel{\sim}{\lra} 
{\rm Der}^S{{\cal L}ie}_{\Bbb H}. 
\ee
\el

{\bf Proof}. The space ${\Bbb H}^*$ has a symplectic structure with values in 
${\cal H}$. So thanks to \cite{K}, there is an isomorphism of Lie algebras
(\ref{6.15.08.2}). One easily checks that it commutes 
with the differentials. The lemma is proved.

\vskip 2mm
{\bf Proof of Proposition \ref{1.17.08.10as1sd}}. 
(i) It is equivalent to Lemma \ref{6.15.08.1}. 
  
(ii) This was done above. 
The Proposition is proved.

\label{HF3sec}

\section{Hodge correlator classes and their interpretations} 
\label{hf6sec}

In Section \ref{hf6.a} we introduce the 
Hodge correlator class. In Sections \ref{hf6.1ref}-\ref{hf5.2ref} 
we recall some standard 
background material on cohomology of DG Lie algebras 
and DG spaces arising from DG Lie algebras. 
Using this, in Section 4.4  we define versions of the Hodge correlator 
classes:  the {\it Hodge vector fields}. In Section \ref{hf5.4ref} we 
 define Hodge vector fields explicitly as (the sum of) correlators 
of the Feynman integral from Section \ref{hf1.6}. We interpret them as 
 vector fields on the formal neighborhood of the trivial local system, providing 
a homotopy action of the Hodge Galois group there.

\vskip 3mm

\subsection{The Hodge correlator class} \la{hf6.a}

The Hodge correlator ${\bf G}$, see (\ref{GGGG}), 
 depends on the choice of a Green current. 
To define the latter, we  
choose  a splitting of the Dolbeaut 
complex into the space of harmonic forms and its orthogonal complement. 
In particular the volume form $\mu$  provides a choice of  the
 harmonic $2n$-form. These choices determine 
differential equation (\ref{7.3.00.1mn}) 
for the Green current. 

\vskip 3mm
\begin{theorem}\label{3/11/07/101q}
(i) The Hodge correlator ${\bf G}$ is a degree $0$ element of 
${\cal C}{\cal L}ie_{\Bbb H}\otimes {\cal H}$.

(ii) One has 
${\delta} {\bf G}=0$. 

(iii) 
Altering either a splitting of the Dolbeaut complex or 
a Green current we get an element $\widetilde {\bf G}$, 
which differs from ${\bf G}$ by a coboundary: 
there is an element 
$
{\bf B} \in {\cal C}{\cal L}ie_{\Bbb H}\otimes {\cal H} 
$ such that 
$$
\widetilde {\bf G} - {\bf G} = {\delta}{\bf B}.
$$

Thus there is  a well defined 
${\delta}$-cohomology class of ${\bf G}$, called the 
{\em Hodge correlator class}:
 $$
{\bf H}_X \in H_0^\delta({\cal C}{\cal L}ie_{\Bbb H}\otimes {\cal H}).
$$
\end{theorem}

We prove this Theorem in Section \ref{hf4.1ref}.

\vskip 3mm

Theorem \ref{3/11/07/101q} 
and Corollary \ref{3.28.08.10} plus Lemma \ref{popo} 
immediately imply Proposition \ref{IDAD}, that is: 

\vskip 3mm
{\it The element ${\bf G}$ provides a derivation 
of the graded Lie algebra ${\cal L}ie_{\Bbb H}$ 
commuting with the differential $\delta$ in ${\cal L}ie_{\Bbb H}$. 
Its homotopy class 
depends only on the Hodge correlator class ${\bf H}_X$}.

\vskip 3mm
\subsection{Derivations and cohomology of DG Lie algebras}  \label{hf6.1ref}
Let $({\cal G}, \delta)$ be a DG Lie algebra, 
i.e. a Lie algebra in the tensor category of complexes. 
There are two slightly different versions of the (DG) Lie algebra cohomology.

\vskip 3mm 
1. The standard cochain complex of ${\cal G}$ with coefficients in 
a DG ${\cal G}$-module $M$ is  
given by 
$$
C^*({\cal G}, M):= S^*({{\cal G}[1]}^{\vee})\otimes M.
$$ 
Here the grading is provided by the gradings 
on ${\cal G}$ and $M$, and the differential is 
the sum of the Chevalley differential 
$\delta_{\rm Ch}$ and the differentials induced by the differentials on  
${\cal G}$ and $M$.  
The  cohomology of this complex is denoted by ${\Bbb H}^{*}({\cal G}, M)$. 

Let ${\cal G}_i$ be DG Lie algebras, and 
$M_i$ is DG Lie modules over ${\cal G}_i$, $i=1,2$.  
Given a morphism of DG Lie algebras 
$\varphi: {\cal G}_1 \to {\cal G}_2$ and a $\varphi$-equivariant morphism 
$\psi: M_1 \to M_2$ we get  maps 
\be \la{4.12.08.1}
{\Bbb H}^*({\cal G}_2, M_2) \stackrel{\varphi}{\lra} 
{\Bbb H}^*({\cal G}_1, M_2) \stackrel{\psi}{\longleftarrow} 
{\Bbb H}^*({\cal G}_1, M_1). 
\ee
If $\varphi$ and $\psi$ are quasiisomorphisms,  
then maps  (\ref{4.12.08.1}) are also isomorphisms. 
This is proved by a spectral sequence 
argument applied to maps of the bicomplexes induced by 
$\varphi$ and $\psi$. 

\vskip 3mm
{\it Motivation}. Let ${\cal G}$ be a Lie algebra. Denote by 
${\rm Der}^{\rm Inn}({\cal G})$ the Lie algebra of inner derivations of ${\cal G}$. 
Then the Lie algebra of outer derivations 
of ${\cal G}$ is given by 
$$
{\rm Der}^{\rm Out}({\cal G}):= \frac{{\rm Der}({\cal G})}{{\rm Der}^{\rm Inn}({\cal G})} 
\stackrel{\sim}{=}H^1({\cal G}, {\cal G}). 
$$
Here ${\cal G}$ acts on itself by the adjoint action. 

\vskip 3mm
Now let $({\cal G}, \delta)$ be a DG Lie algebra
 which, considered as a graded Lie algebra, 
is a free non-abelian graded Lie algebra. 
Then  ${\Bbb H}^1({\cal G}, {\cal G})$ can be presented as follows: 

The group 
 $H^1({\cal G}, {\cal G})$ is the only non-zero  Lie algebra cohomology group. 
 It is graded space, and the differential  $\delta$ 
provides a differential $\Delta$ on $H^1({\cal G}, {\cal G})$. 
The spectral sequence argument tells that
$
{\Bbb H}^i({\cal G}, {\cal G}) = 
H^{i-1}_{\Delta}\Bigl(H^1({\cal G}, {\cal G})\Bigr).
$ 
It follows that 
$$
{\Bbb H}^1({\cal G}, {\cal G}) = 
H^{0}_{\Delta}\Bigl(H^1({\cal G}, {\cal G})\Bigr) = 
$$
$$
\frac{\mbox{outer derivations of the 
graded Lie algebra ${\cal G}$ commuting 
with the differential $\delta$}}{\mbox{derivations homotopic to zero}}.
$$
 
\vskip 3mm
2. We  need a cohomology group responsible for   
derivations rather then outer derivations.  

Let us define the reduced cochain complex 
of a DG Lie algebra ${\cal G}$ with coefficients in 
a DG ${\cal G}$-module $M$ as  
$$
\widetilde C^*({\cal G}, M):= S^{>0}({\cal G}[1]^{\vee})\otimes M.
$$ 
It is quasiisomorphic to the complex
$
{\rm Cone}\Bigl( S^*({\cal G}[1]^{\vee})\otimes M \lra M\Bigr).
$ 
Here the map is induced by the augmentation map $S^*({\cal G}[1]^{\vee})\to \Q$. 
The reduced cohomology 
$\widetilde {\Bbb H}^*({\cal G}, M)$ are the cohomology of this complex. 

Given ${\cal G}_i$-modules $M_i$, $i=1,2$, a morphism of DG Lie algebras 
$\varphi: {\cal G}_1 \to {\cal G}_2$, and a $\varphi$-equivariant morphism 
$\psi: M_1 \to M_2$ we get  morphisms
\be \la{4.12.08.1a}
\widetilde {\Bbb H}^*({\cal G}_2, M_2) \stackrel{\varphi}{\lra} 
\widetilde {\Bbb H}^*({\cal G}_1, M_2) \stackrel{\psi}{\longleftarrow} 
\widetilde {\Bbb H}^*({\cal G}_1, M_1). 
\ee
If $\varphi$ and $\psi$ are quasiisomorphisms, the same is true 
for the morphisms (\ref{4.12.08.1a}). 

 Let us assume that ${\cal G}$, considered as a graded Lie algebra, 
is a free non-abelian graded Lie algebra. 
Then $\widetilde H^1({\cal G}, {\cal G})$ is the only non-zero cohomology group, identified with the derivations of ${\cal G}$. 
It is a graded space, the differential  $\delta$ 
provides a differential $\Delta$ on $\widetilde H^1({\cal G}, {\cal G})$, and 
$$
\widetilde {\Bbb H}^i({\cal G}, {\cal G}) = 
H^{i-1}_{\Delta}\Bigl(\widetilde H^1({\cal G}, {\cal G})\Bigr).
$$
It follows that 
$$
\widetilde {\Bbb H}^1({\cal G}, {\cal G}) = 
H^{0}_{\Delta}\Bigl(\widetilde H^1({\cal G}, {\cal G})\Bigr) = 
$$
$$
\frac{\mbox{derivations of the 
graded Lie algebra ${\cal G}$ commuting 
with the differential $\delta$}}{\mbox{derivations homotopic to zero}}.
$$

\vskip 3mm
\subsection{The DG space ${\cal G}[1]$} \la{hf5.2ref}

The category of affine DG spaces is the opposite  
to the category of graded 
commutative algebras equipped  
with a degree one derivation $v$ with $v^2=0$. 
Such a $v$ is called  
a {\it homological vector field} on a DG space. 

Let ${\cal G}$ be a DG Lie algebra over a field $k$. 
Then ${\cal G}[1]$ has a structure of a DG 
space with the 
 algebra of functions 
\begin{equation} \label{CSq1}
{\cal O}_{{\cal G}[1]}:= S^*({\cal G}[1]^{\vee}) = C^*({\cal G}, k).
\end{equation}
It is identified with 
the standard Chevalley complex of the DG Lie algebra ${\cal G}$. 
Its differential provides a homological vector field on ${\cal G}[1]$. 

Alternatively, the 
homological vector field on ${\cal G}[1]$ is also known as the 
{\it Chern-Simons vector field $Q_{CS}$}. 
It is a  quadratic vector field   given by 
\begin{equation} \label{CSq}
\stackrel{\cdot }{\alpha } = Q_{CS}(\alpha ) 
:= d\alpha + \frac{1}{2}[\alpha, \alpha].
 \end{equation}
It is an odd vector field: the commutator $[\alpha, \alpha]$ vanishes if 
$\alpha$ is even. One has 
\begin{equation} \label{CS}
Q_{CS}^2 = \frac{1}{2}[Q_{CS}, Q_{CS}]=0.
\end{equation}
Indeed, $Q_{CS}$ is a sum of linear and quadratic terms, $Q_1+Q_2$. 
Then $Q_2 \circ Q_2 =0$ by Jacobi identity, $Q_1 \circ Q_1 - d^2 =0$, 
and $Q_2 \circ Q_1+ Q_1 \circ Q_2  =0$ 
is equivalent to the Leibniz rule for $d$.

The Lie algebra $\widetilde {\rm Vect}_{{\cal G}[1]}$ 
of formal vector fields on ${\cal G}[1]$ vanishing at the origin is 
a Lie subalgebra of derivations of the graded commutative algebra (\ref{CSq1}). 
Since the latter algebra is free, a derivation is determined by 
its action on the generators, and we have 
$$
\widetilde {\rm Vect}_{{\cal G}[1]}= 
\widetilde C^*({\cal G}, {\cal G})[1].
$$
The graded  Lie algebra structure on 
$\widetilde C^*({\cal G}, {\cal G})[1]$ is given by 
the commutator of vector fields. 
Since $Q_{CS}$ is a homological
 vector field, we get 
 a differential
$$
D:= [Q_{CS}, \ast], \qquad D^2=0.
$$
One has 
$$
\widetilde {\Bbb H}^{*+1}({\cal G}, {\cal G}) = H_D^*(\widetilde 
{\rm Vect}_{{\cal G}[1]}) = 
H_D^{*+1}(\widetilde C^*({\cal G}, {\cal G})).
$$
In particular, 
$$
\widetilde {\Bbb H}^1({\cal G}, {\cal G}) = 
H^{0}_{D}\Bigl(\widetilde 
{\rm Vect}_{{\cal G}[1]}\Bigr) = 
$$
$$
\frac{\mbox{Vector fields on ${\cal G}[1]$, vanishing at $0$, commuting with the 
Chern-Simons vector field 
$Q_{CS}$}}{\mbox{Commutators with $Q_{CS}$}}.
$$
We  
identify these classes with the isomorphism classes of deformations of the 
DG space ${\cal G}[1]$ over the one dimensional odd line ${\rm Spec}(\Z[\varepsilon])$: 
a degree zero  vector field $Q$ commuting with $Q_{CS}$ provides a homological 
vector field $Q_{CS} + \varepsilon Q$, and vice versa. 

\vskip 3mm
Assume now that the DG Lie algebra ${\cal G}$ 
has an even invariant  non-degenerate scalar product $Q(\ast,\ast)$. 
Then the  space ${\cal G}[1]$ has an even 
Poisson structure provided by the bivector corresponding to the 
induced scalar product in ${\cal G}^{\vee}$.  
The vector field $Q_{CS}$ is a Hamiltonian vector field with the 
Hamiltonian given by the Chern-Simons functional, a cubic polynomial function on  
${\cal G}[1]$: 
$$
CS(\alpha):= \frac{1}{2}(\alpha, d\alpha) + 
\frac{1}{6}(\alpha, [\alpha, \alpha]).
$$
One has $\{CS, CS\}=0$ in agreement with  (\ref{CS}).

\vskip 3mm
{\bf Example}. Let $X$ be a compact K\"ahler manifold, and 
${\cal G}$ a Lie algebra with an invariant non-degenerate scalar product 
$Q(\ast, \ast)$. 
Then ${\cal A}^*(X) \otimes {\cal G}$ is a 
DG Lie algebra with a scalar  product of degree $2-2n$, 
where $n = {\rm dim}_\C X$:
$$
(\alpha, \beta):= (-1)^{{\rm deg}(\alpha)}\int_X Q(\alpha \wedge \beta), 
$$ 
and ${\cal A}^*(X)[1] \otimes {\cal G}$ 
is a DG space  with an even symplectic  structure.

\vskip 3mm
\subsection{Hodge vector fields} \la{hf5.3sec}
Let $X$ be a compact K\"ahler manifold.   
The Formality Theorem \cite{DGMS} asserts that assigning to a cohomology class 
of $X$ its harmonic representative, we get a quasiisomorphism 
of the DG commutative bigraded algebras $H^{*,*}_X \lra {\cal A}^{*,*}(X)$.

Denote by $\overline {\cal A}^{*,*}_X$ the quotient of ${\cal A}^{*,*}(X)$ 
by the space spanned by the constants and the Kahler volume form, which are 
the harmonic representatives $H^0_X \oplus H^{2n}_X $. So
the Formality Theorem  implies 

\bc \la{3.20.09.1}
The DG commutative bigraded algebras $\overline {\cal A}^{*,*}_X$ and $\overline H^{*,*}_X$ 
are quasiisomorphic. 
\ec

Let ${\cal G}$ be a Lie algebra. 
The bigraded DG Com's in Corollary  \ref{3.20.09.1} 
provide bigraded DG Lie algebras 
$$
{\cal G}_{{\rm A}} : = \overline {\cal A}^{*,*}_X  
\otimes {\cal G}, \qquad 
 {\cal G}_{{\rm H}}:= \overline H^{*,*}_X \otimes {\cal G}.
$$
We assume that ${\cal G}$ is a complexification of a real Lie algebra. 
Since both DG Com's from Corollary  \ref{3.20.09.1} have real structures, 
there is a complex involution $\sigma$  acting on these DG Lie algebras, 
and hence on their cohomology.  

Below we assume that ${\cal G}$ is a Lie algebra
with non-degenerate invariant scalar product $Q(\ast, \ast)$.

Given a point $a \in X$, and using the isomorphism of DG spaces (\ref{4.20.08.5}), 
we transform the Hodge correlator class into a cohomology class
$$
{{\rm Q}_{\rm Hod}} = {\rm Q}_{X, {\cal G},a}\in \widetilde 
{\Bbb H}^1({\cal G}_{{H}}, {\cal G}_{{\rm H}}),
\qquad 
 \sigma({{\rm Q}_{\rm Hod}}) = - {{\rm Q}_{\rm Hod}}.
$$
Theorem \ref{PPPq}  implies that it is functorial:  
A map $f:X \to Y$ of Kahler manifolds gives rise to a DG Lie algebra map 
$
f^*: {\cal G}_{{\rm H}^*(Y)} \lra {\cal G}_{{\rm H}^*(X)}, 
$ 
which intertwines, up to a homotopy,  
derivations ${\rm Q}_{X, {\cal G}, a}$ and ${\rm Q}_{Y, {\cal G}, a}$:
$$
 f^* \circ {\rm Q}_{Y, {\cal G}, f(a)}  \quad 
\mbox{is homotopic to} \quad {\rm Q}_{X, {\cal G}, a} \circ f^*. 
$$

\vskip 3mm
Corollary  \ref{3.20.09.1} implies that the DG Lie algebras ${\cal G}_{{\rm A}}$ and 
${\cal G}_{\rm H}$ are quasiisomorphic. So there is an isomorphism
$$
\widetilde {\Bbb H}^1({\cal G}_{\rm H}, {\cal G}_{\rm H})  = 
\widetilde {\Bbb H}^1({\cal G}_{{\rm A}}, {\cal G}_{{\rm A}}).  
$$
This leads to 
\bt \la{4.13.08.1a}
Given a point $a \in X$, there is a functorial cohomology class
\be \la{4.13.08.1asd}
{\Bbb V}_{\rm Hod}\in \widetilde {\Bbb H}^1({\cal G}_{{\rm A}}, 
{\cal G}_{{\rm A}}) \stackrel{\sim}{=} H^{0}_{\Delta}\Bigl(\widetilde H^1({\cal G}_{{\rm A}}, {\cal G}_{{\rm A}})\Bigr).
\ee
\et

{\bf Remark}. The 
action of the group $\C^*_{\C/\R}$ on the Dolbeaut bicomplex of $X$ 
providing the Hodge bigrading gives rise to an action 
on $\widetilde H^1({\cal G}_{{\rm A}}, {\cal G}_{{\rm A}})$. It 
does not commute with the differential $\delta$. Thus it is not the action 
of the pure Hodge Galois group.

\vskip 3mm
\subsection{Constructing the Hodge vector fields} \la{hf5.4ref}

Below we  give 
an explicit construction of the Hodge vector fields ${\rm Q}_{\rm Hod}$. 

Let us show that, given a point $a \in X$, 
 a version of the Hodge correlator construction provides a 
degree zero linear map:
\be \la{2.25.08.2}
{\rm S}^*({\cal G}_{\rm H}[1])\otimes H_{2n}(X)[-2]
 \lra \C. 
\ee
 Observe that $H_{2n}(X)$ is in  degree $-2n$, and is canonically 
isomorphic $\C$ since $X$ is oriented.

Choose a Green current 
$G_a(x,y)$ on $X \times X$. 
Take a plane trivalent tree $T$. 
Decorate its external vertices by harmonic forms $\alpha_0, \ldots , \alpha_m$. 
Let $Q' \in S^2{\cal G}\subset {\cal G}\otimes{\cal G}$ be the element
provided by the form on ${\cal G}$. 
We assign to every internal edge $E$ of $T$ the element
$G_a(x,y)\otimes Q' $ 
in the space of ${\cal G}$-valued distributions on 
$X^{\{\mbox{vertices of $E$}\}}$. Using  
 the polydifferential operator $\omega$ and, at every internal vertex, 
 the canonical trilinear form 
$$
\langle l_1, l_2, l_3 \rangle = Q(l_1, [l_2, l_3]): \Lambda^3{\cal G} \to \C, 
$$
we cook up a differential form on 
$
X^{\{\mbox{internal vertices of $T$}\}}.
$  
It provides a current there. Integrating it, and taking then the sum over all 
plane trivalent trees, we get the map (\ref{2.25.08.2}). 

One checks that the map (\ref{2.25.08.2}) is a degree zero map. It
 provides a function on the space 
${\cal G}_{{\rm H} }[1]$: 
its value on an ${\cal G}$-valued harmonic form $\alpha$ 
is obtained by evaluating the map (\ref{2.25.08.2}) on 
$\alpha \otimes \ldots \otimes \alpha$. 
It is a formal power series 
on the space ${\cal G}_{{\rm H} }[1]$, 
i.e. a formal sum of its homogeneous components.  Since 
${\cal G}_{{\rm H} }[1]$ 
 is a  Poisson space, it gives rise to  
a formal Hamiltonian vector field,  a {\it Hodge vector field}: 
 \be \la{2.25.08.2a}
{\rm Q}_{\rm Hod}  \in S^*({\cal G}_{{\rm H}}[1]^{\vee}) \otimes  
{\cal G}_{{\rm H}}[1].
\ee

\vskip 3mm

Theorem \ref{18.4.08.10}ii) plus Theorem \ref{3/11/07/101q} imply that 
the Hodge vector field ${\rm Q}_{\rm Hod}$ commutes 
with the vector field $Q_{CS}$.
Furthermore, a
 different choice of the splitting of the Dolbeaut complex and the Green current  
results in a Hodge vector field $\widetilde {\rm Q}_{\rm Hod}$ 
such that 
$
\widetilde {\rm Q}_{\rm Hod} - {\rm Q}_{\rm Hod} \in {\rm Im} [Q_{CS}, \ast].
$ 
 Functoriality follows similarly from Theorem \ref{PPPq}.

\vskip 3mm
The cohomology 
$\widetilde {\Bbb H}^1({\cal G}_{{H}}, {\cal G}_{{H}})$ are bigraded 
by the Hodge bigrading.

\bl \la{2.25.08.6}
The bidegree $(p,q)$ component of 
${\rm Q}_{\rm Hod}$ is a polynomial vector field on ${\cal G}_{{\rm H}}[1]$. 
It can be 
 non-zero only if 
$p,q>0$. 
\el

{\bf Proof}. Let $(p_i, q_i)$ be the Hodge bidegree of a form decorating $i$-th 
an external vertex of a trivalent tree $T$. Then it contributes to the 
$(p,q)$-component of  
${\rm Q}_{\rm Hod}$, where 
$p = \sum_{i=0}^m p_i -n$ and $q = \sum_{i=0}^m q_i -n$. 
 Since $p_i, q_i \geq 0$ and $p_i+q_i=2n$, 
for a given $(p,q)$   
the number $m+1$ of vertices of the tree $T$ is bounded from above. 
So the $(p,q)$-component is given by a finite sum. 
\vskip 3mm

\subsection{The Hodge Galois group action on close to trivial local systems}
\la{hf6.5ref}

\vskip 3mm
Theorem \ref{4.13.08.1a} provides 
a formal vector field ${\Bbb V}_{\rm Hod}$ on the DG space 
${\cal G}_{{\cal A}^*(X, a)}$. It is  well defined up to a homotopy, 
commutes with the Chern Simons vector field $Q_{CS}$, 
and covariant, up to homotopy,  
for the maps of Kahler manifolds. 

Below we show that the 
vector fields    
${\Bbb V}_{\rm Hod}$ 
give rise to a 
functorial vector field ${\Bbb V}_{X, {\cal G}, a}$ on the 
formal neighborhood of the trivial $G$-local system in the space of 
complex local systems on $X$. 
It provides an action of the Lie 
algebra ${\rm L}_{\rm Hod}$: 
the generator $G$ acts by the vector field ${\Bbb V}_{X, {\cal G}, a}$.

\vskip 3mm
Given a homological vector field $Q$ on a DG space $V$, 
let $V^Q \subset V$ be 
the set of its 
 fixed points. The tangent space $T_f V$ 
at a point $f$ is a graded vector space. 
The linearization of $Q$ at a point $f \in  V^Q$ 
provides a linear operator $Q= Q_{(f)}$ in $T_f V$, 
described as follows. The vector fields 
vanishing at a point form a Lie subalgebra. So  
the value $[Q, Q']_f$ of the commutator at $f$ 
depends only on the value $Q'_f$ of $Q'$ at $f$. 
Set $Q_{(f)}(v):= [Q, Q']_f$ where $Q'_f=v$. 
Since $Q^2=0$,  $Q_{(f)}$ is a differential 
in  $T_f V$. We denote it simply by $Q$. 

Let $H$ be a vector field commuting with $Q$. 
Then at a $Q$-fixed  point $f$  
one has $Q(H_f) =0$. 
Since  $[Q, R]_f = Q(R_f)$, 
adding to $H$ a commutator $[Q, R]$ with $Q$ and 
restricting the result to $f \in V^Q$ 
we do not change the cohomology class 
of $H_f$ in $H^*_Q(T_f V)$. 

\vskip 3mm
In our case 
$Q$ is the Chern-Simons vector field $Q_{CS}$ 
on the DG space ${\cal G}_{{\rm A}}[1]$, and 
$H = {\Bbb V}_{\rm Hod}$.  
The set ${\rm Con}_{\cal F}$ of fixed points of the vector field $Q_{CS}$ 
on ${\cal G}_{{\rm A}}[1]$ 
is the set of flat DG connections. 
The subset 
${\rm Con}^1_{\cal F}$ of ${\cal G}$-valued 
$1$-forms $\alpha \in {\rm Con}_{\cal F}$ is 
identified with the flat connections $d + \alpha$. 
The flat connection $d + \alpha$ provides 
a local system ${\cal V}_\alpha$. 
The complex given by the operator $Q_{CS}$ acting on the tangent space 
to $\alpha \in {\rm Con}^1_{\cal F}$ is identified with the 
de Rham complex of the local system ${\cal V}_\alpha$. Thus 
one has 
\be \la{4.20.08.1}
H_{Q_{CS}}^*(T_\alpha {\cal G}_{{\rm A}}[1]) = H^*(X, {\cal V}_\alpha).
\ee

\bp \la{4.19.08.3}
The vector field ${\Bbb V}_{\rm Hod}$  
provides a flow on the formal neighborhood of the trivial 
local system  
on $X$. 
\ep

{\bf Proof}. 
Thanks to (\ref{4.20.08.1}), the formal vector field  ${\Bbb V}_{\rm Hod}$ restricted to 
a flat connection $d+ \alpha$  give 
rise to a cohomology class of the corresponding local system:
$
{\Bbb V}_{\rm Hod}^\alpha \in H^1(X, {\cal V}_\alpha)$. One easily  
sees that the class ${\Bbb V}_{\rm Hod}^\alpha$ does not change 
under the action of the gauge transformations.


\vskip 3mm

The tangent space to a local system ${\cal V}_\alpha$ 
in the space of all ${\cal G}$-local systems on $X$ is identified with 
$H^1(X, {\cal V}_\alpha)$. So we get a  vector field 
${\Bbb V}_{X, {\cal G}, a}$ 
on the formal neighborhood of a 
 trivial local system on $X$. It is functorial 
for the maps $f: Y \to X$ of Kahler manifolds: 
$
f^*{\Bbb V}_{X, {\cal G}, f(a)} = {\Bbb V}_{Y, {\cal G}, a}. 
$ 



\label{HF6sec}


\subsection{Generalizations: Hodge correlators for local systems}
Let ${\cal E}$ be a real 
polarized variation of Hodge structures over a compact Kahler manifold $X$. 
Consider the smooth de Rham complex of the local system of the endomorphisms of ${\cal E}$
$$
{\cal G}_{{\cal E}}:=  {\cal A}^{*,*}(X) \otimes {\rm End}{\cal E}.
$$
It has two standard differentials $D$ and $D^\C$. 
The $dd^\C$-lemma holds for these differentials. 

The complex ${\cal G}_{{\cal E}}$ has  a natural bigraded DG Lie algebra structure. 
The Hodge vector field ${\rm Q}_{\rm Hod}$  has a straitforward analog in this context: 

\bt \la{4.13.08.1fr}
1. Given a point $a \in X$, there is a cohomology class
$$
{\rm Q}_{{\cal E},a}\in \widetilde {\Bbb H}^1({\cal G}_{{\cal E}}, 
{\cal G}_{{\cal E}}),
\qquad 
 \sigma({\rm Q}_{{\cal E},a}) = - {\rm Q}_{{\cal E},a}.
$$
2. It is functorial:  
A map $f:X \to Y$ of Kahler manifolds gives rise to a DG Lie algebra map 
$
f^*: {\cal G}_{{\cal E}} \lra {\cal G}_{
f^*{\cal E}}, 
$ 
which intertwines, up to a homotopy,  derivations ${\rm Q}_{{\cal E}, f(a)}$ 
and ${\rm Q}_{f^*{\cal E}, a}$:
$$
 f^* \circ {\rm Q}_{{\cal E}, f(a)}  \quad 
\mbox{is homotopic to} \quad {\rm Q}_{f^*{\cal E}, a} \circ f^*. 
$$
\et
There is an explicit construction 
of  the class ${\rm Q}_{{\cal E},a}$ generalizing the one from Section \ref{hf5.4ref}.

 \label{HF7sec} 
\section{Functoriality of Hodge correlator classes} \la{hf4sec}

In this Section we prove main properties 
of the Hodge correlators for compact Kahler manifolds. They tell that 
 the Hodge correlator cohomology class 
is well defined 
(i.e. does not depend on the choices made in its definition) 
and functorial. In Section \ref{hf4.1} we 
collect some technical results. The reader can skip them, and return when needed.

\subsection{Preliminary results} \la{hf4.1}

\paragraph{\it A formula for $\delta$.} 
Let $\omega\in H^k$. 
 Define  $s_\omega\in H_{2n-k}$  by the Poincare duality operator: 
$$
s_{\omega} := H_{2n} \cap \omega, \qquad 
\langle s_{\omega}, \alpha\rangle := 
\langle H_{2n}, \omega \cup \alpha\rangle = \int_X\omega \cup \alpha.
$$
Recall 
$
(H^{2n} \cap h)\cup \beta:= H^{2n} \cdot (h, \beta)$ and 
$\langle h_1 \cap H^{2n} \cap h_2\rangle := 
(h_1, H^{2n} \cap h_2). 
$ 
One has  $\langle {H}^{2n}\cap  s_\omega \rangle = \omega$, and hence 
$\langle {\cal H}\cap  \overline s_\omega \rangle = 
\overline\omega$.  We use dual bases $\{\alpha_s\}$ and $\{h_{\alpha_s}\}$, so that 
$(h_{\alpha_s}, \alpha_{s'}) = \delta_{s s'}$  and  
notation $\alpha|h:= \sum_s\overline \alpha_s |\overline h_{\alpha_s}$, 
as in  Section \ref{hf2sec}.
Set 
\be \la{thatf}
\delta (\overline \alpha | \overline h): =  
(-1)^{|\overline \alpha|} \overline \alpha  | \delta \overline h.
\ee

\bl \la{4.21.08.11}
One has 
\begin{equation} \label{9.29.07.5m}
\delta(\alpha | h) = -
\alpha| h  
\otimes  \alpha| h.
\end{equation}

Equivalently,
 \be \la{4.21.08.1a}
\delta:\overline s_{\omega}  \lms \int_X \alpha|h \otimes \alpha|h\otimes \overline \omega. 
\ee
\el

{\bf Proof}. One has 
\be \la{4.21.08.1}
\delta: x \lms \sum_{\alpha_3}
\int_X (\alpha|h \otimes 
\alpha | h\otimes \overline \alpha_3) 
\cdot \langle \overline h_{\alpha_3}\cap {\cal H}\cap 
x\rangle.
\ee
Furthermore, 
$ 
\langle \overline a\cap {\cal H}\cap  \overline b \rangle:= 
(-1)^{{\rm deg}(a)}\langle a\cap {H}^{2n}(X)\cap  b \rangle.
$ 
Since 
$$
\sum_{\alpha_3}
\overline \alpha_3 
\cdot \langle \overline h_{\alpha_3}\cap {\cal H}\cap 
\overline s_{\omega}\rangle = \overline \alpha_3 
\cdot \langle \overline h_{\alpha_3}, \overline {\omega}\rangle = \overline {\omega}.
$$
we get (\ref{4.21.08.1a}). Since $h_{\alpha^{\vee}_k}= s_{\alpha_k}$, 
we see that (\ref{4.21.08.1a})  is equivalent to (\ref{9.29.07.5m}). 
Indeed, 
$$
\sum_s (-1)^{|\overline \alpha_s^{\vee}|}\overline \alpha_s^{\vee}| 
\delta h_{\overline \alpha_s^{\vee}} \stackrel{(\ref{4.21.08.1a})}{=}
\sum_s  (-1)^{|\overline \alpha_s^{\vee}|}\overline \alpha_s^{\vee}\int_X 
\overline \alpha_s \otimes \alpha|h \otimes \alpha|h = 
 -\alpha|h \otimes \alpha|h. 
$$
The lemma is proved.


\vskip 3mm
Given a tree $T$, denote by $T^0$ the tree 
obtained by removing all external edges of  $T$. 

We extend the definition of the Hodge correlators to the case 
when the forms $\alpha_i$ may be arbitrary forms by applying the same 
recipe as in Section \ref{hf2sec}.

\vskip 3mm
\paragraph{\it Two useful integrals.} Let $T$ be a plane trivalent tree 
decorated by $\alpha_1 \otimes \ldots \otimes \alpha_m$, 
no assumptions made about the forms $\alpha_i$. 
Let us assign to internal edges $E$ of the tree $T$ 
some forms/currents $g_E$ on $X ^{\mbox{{\rm vert}(E)}}$. 
Consider the following integrals assigned to this data:
\footnote{Here and later on the symmetrization ${\rm Sym}$ is taken with respect to the 
shifted by $1$ degree $|g|:= {\rm deg}(g) +1$}
\begin{equation} \label{inteta}
\eta_T(\alpha_1 \otimes \ldots \otimes \alpha_m):= 
\int_{X ^{\mbox{{\rm vert}($T^0$)}}}d^\C 
g_{1} \wedge \ldots \wedge 
d^\C g_{{2m-1}} \wedge \alpha_1 \wedge \ldots \wedge \alpha_m, 
\end{equation}
\begin{equation} \label{inteta}
\xi_T(\alpha_1 \otimes \ldots \otimes \alpha_m):= 
{\rm Sym}_{g_1, ..., g_{2m-1}}\int_{X ^{\mbox{{\rm vert}($T^0$)}}}
g_{1} \wedge d^\C g_{2} \ldots \wedge 
d^\C g_{{2m-1}} \wedge \alpha_1 \wedge \ldots \wedge \alpha_m.
\end{equation}
We assume that products of currents $g_i$ make sense.  
Write $a \sim_{\Q^*} b$ if $a = \lambda b$ for $\lambda \in \Q^*$. 
Then 
\begin{equation} \label{5.26.07.1}
\kappa_T (d^\C\beta \otimes \alpha_1 \otimes 
\ldots \otimes  \alpha_k) \sim_{\Q^*} \eta_T(\beta \otimes \alpha_1 \otimes 
\ldots \otimes  \alpha_k)\qquad \mbox{assuming $d^\C \alpha_i=0$}, 
\end{equation}
where we assume that the forms $\alpha_i$ are 
homogeneous with respect to the Hodge bidegree $(p,q)$. 

\begin{lemma} \label{corT1}  
If $d^\C \alpha_i=0$ for all $i$ and 
$m>3$, then 
$\eta_T(\alpha_1 \otimes \ldots \otimes \alpha_m)= 0$. 
\end{lemma}

{\bf Proof}. The integrand is $d^\C$ of the form 
$g_{1} \wedge d^\C g_{2} \wedge \ldots \wedge 
d^\C g_{{2m-1}} \wedge \alpha_1 \wedge \ldots \wedge \alpha_m$.  
\vskip 3mm

\begin{figure}[ht]
\centerline{\epsfbox{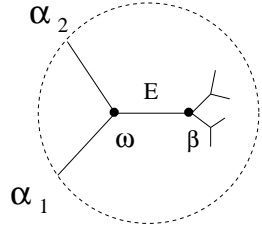}}
\caption{Only such diagrams can contribute.}
\label{fey1}
\end{figure}

Denote by $p_1, p_2$ the projections of $X ^{2}$ to the factors. 
\begin{lemma} \label{7.24.07.1}
a) Let $T$ be a plane trivalent tree 
decorated by $W= {\cal C}(\alpha_1 \otimes \ldots \otimes \alpha_m)$, 
$d^\C \alpha_i=0$. We assign to an internal edge $E$ of $T$ a decomposable current $p_1^*\omega \wedge p_2^*\beta$, where 
$d^\C\omega=0$, and to the other 
internal edges $F$ of $T$ some 
currents $g_F$. We 
assume that the corresponding 
form $\kappa_T(W)$, see (\ref{3.8.05.15.ab}), is  a current. 
Then its  integral 
is zero 
unless the following condition holds:
\begin{equation} \label{9.24.07.5}
\mbox{$E$ 
is an external edge of $T^0$, and $\omega$ labels 
an external vertex of $E$ in $T^0$, as on Fig. \ref{fey1}}.
\end{equation}

b) Let us assume (\ref{9.24.07.5}). 
Denote by $T_\beta$ and $T_\omega$ the two trees  obtained from $T$ 
by cutting the edge $E$, so that  
if $T_\omega$ is decorated by 
$\omega \wedge \alpha_{1}\wedge \alpha_2$ (see Fig. \ref{fey2}). 
The integral is proportional to 
\be \la{CONT1}
\int_{X }\alpha_{1}\wedge \alpha_2\wedge \omega 
 ~\cdot~ \int_{X }\kappa_{T'_E}(d^\C\beta\otimes 
\alpha_3 \otimes \ldots \otimes \alpha_{m}), 
\ee
unless $E$ is the only internal edge of $T$, in which case the contribution is 
\be \la{CONT2}
\int_{X }\alpha_{1}\wedge \alpha_2\wedge \omega 
 ~\cdot~ \int_{X }(\beta\wedge \alpha_3 \wedge \alpha_{4}).
\ee
\end{lemma}

\begin{figure}[ht]
\centerline{\epsfbox{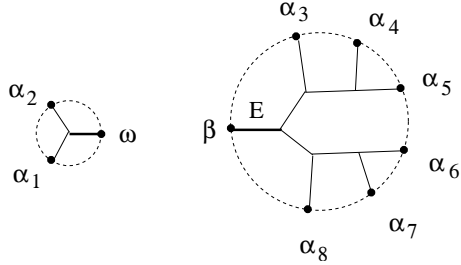}}
\caption{A possible non-zero contribution of 
a decomposable current $p_1^*\omega \wedge p_2^*\beta$.}
\label{fey2}
\end{figure}

{\bf Proof}. The external edges  of the trees $T_\beta$ and $T_\omega$ 
inherit decorations  $\omega$ and  $\beta$ or $d^\C\beta$. 
If $T^0=E$, we get integral (\ref{CONT2}).

Assume that the tree $T^0$ has more then one edge. 
Then we may assume that 
$E$ is decorated by $\omega$ and  $d^\C\beta$. Indeed, 
there is a unique internal edge of $T$ such that the 
polydifferential operator $\omega$ acts on the current assigned to this edge 
as the identity, and this edge can be chosen arbitrarily. 
So we may assume that $E$ is not this edge.   Moreover, employing 
$\xi_T(W)$  instead of  $\kappa_T(W)$,  
we change the integral by a binomial coefficient, depending on the type of $W$, 
see (\ref{5-x20.1s}). 
Thus we may assume that $\omega$ acts 
to the current assigned to this edge by $d^{\C}$. It remains to notice that 
$d^{\C}\omega = 0$.

We claim that if $T_\omega$ has internal edges, then 
the corresponding Hodge correlator type integral is zero. 
Indeed, the contribution to the integral 
of the tree $T$ with the edge $E$ decorated as above 
 is 
a sum of rational multiples of 
$$
\xi_{T_\beta}(d^\C \beta \otimes ...)  
\cdot \eta_{T_\omega}(\omega \otimes ...), \quad \mbox{or} \quad 
\eta_{T_\beta}(d^\C \beta \otimes ...)  
\cdot \xi_{T_\omega}(\omega \otimes ...).
$$
Since  $d^\C \alpha_i=0$ for every $i$,  $\eta_{T_\beta}(d^\C \beta \otimes ...)  =0$. 
Indeed, if the tree $T_\beta$ has internal 
edges, this follows from Lemma \ref{corT1}. Otherwise the integral is 
$
\int_X d^\C \beta \wedge \alpha_1 \wedge \alpha_2= 
\int_X d^\C (\beta \wedge \alpha_1 \wedge \alpha_2) = 0.
$  
Since the tree $T_\omega$ has internal edges, 
$\eta_{T_\omega}(\omega \otimes ...)=0$  by Lemma \ref{corT1}.  
The claim is proved. 

So $T_\omega$ has no internal edges, and we get integral (\ref{CONT1}). 
The Lemma is proved.

\vskip 3mm

We use a notation $\stackrel{(-)}{\partial}$ 
meaning ``either ${\partial}$ or 
${\overline \partial}$''.  
Consider the following integrals
\begin{equation} \label{intetaq}
{\rm Sym}_{2m-2}\int_{X ^{\mbox{{\rm vert}($T^0$)}}}
\stackrel{(-)}{\partial} g_{1} \wedge \ldots \wedge 
\stackrel{(-)}{\partial} 
g_{{2m-1}} \wedge \alpha_1 \wedge \ldots \wedge \alpha_m.
\end{equation}
\begin{equation} \label{intetaqq}
{\rm Sym}_{2m-2}\int_{X ^{\mbox{{\rm vert}($T^0$)}}}
g_{1} \wedge \stackrel{(-)}{\partial} g_{2} \wedge \ldots \wedge 
\stackrel{(-)}{\partial} 
g_{{2m-2}} \wedge {\partial}{\overline \partial} 
g_{{2m-1}}\wedge \alpha_1 \wedge \ldots \wedge \alpha_m.
\end{equation}

\begin{proposition} \label{17:24i} If the tree $T^0$ is not empty
 and $d\alpha_i=0$, then 
integrals (\ref{intetaq}) and (\ref{intetaqq}) are zero. 
\end{proposition} 

\begin{lemma} \label{17:24}
For any forms $\varphi_i$  and a 
form $\omega$ such that $\partial\omega = \overline \partial\omega = 0$,  one has   
\begin{equation} \label{nab}
({\rm i}) ~ {\rm Sym}_2\int_{X }\omega \wedge \stackrel{(-)}{\partial} \varphi_1 \wedge \stackrel{(-)}{\partial} \varphi_2 = 0, \qquad 
({\rm ii}) ~ {\rm Sym}_3\int_{X }\stackrel{(-)}{\partial} \varphi_1 
\wedge \stackrel{(-)}{\partial} \varphi_2 \wedge 
\stackrel{(-)}{\partial}  \varphi_3  = 0. 
\end{equation}
\end{lemma}

{\bf Proof}. 
(i) There are three cases, depending 
on how many ${\overline \partial}$'s we have in the integrand. 
The two ${\overline \partial}$'s case is 
reduced by complex conjugation to the case without ${\overline \partial}$'s.  
In the latter case the integral equals 
$\int_{X }\partial (\omega \wedge \varphi_1 \wedge \partial \varphi_2) = 0$. 
Finally, in the one  ${\overline \partial}$ case 
$$
(\partial \varphi_1 \wedge \overline \partial \varphi_2 + 
 \overline \partial \varphi_1 \wedge \partial \varphi_2)\wedge \omega  = 
\partial \left(\varphi_1 \wedge \overline \partial \varphi_2  
\wedge \omega\right) + 
\overline \partial \left(
 \varphi_1 \wedge \partial \varphi_2 \wedge \omega\right). 
$$
These are complete derivatives. So integrating
 over $X $ we get zero. 

(ii) By using the  complex conjugation we may assume that  we have 
no or one ${\overline \partial}$. In the 
first case
 the form ${\partial} \varphi_1 \wedge {\partial} \varphi_2 \wedge 
{\partial}  \varphi_3$ is a complete derivative. 
In the second   the integrand is 
$$
(\partial \varphi_1 \wedge \overline \partial \varphi_2 + 
 \overline \partial \varphi_1 \wedge \partial \varphi_2)\wedge \partial \varphi_3 
+ \partial \varphi_1 \wedge \partial \varphi_2 \wedge \overline \partial \varphi_3
=
$$
$$ 
\partial \left(\varphi_1 \wedge \overline \partial \varphi_2  
\wedge \partial \varphi_3\right) + 
\overline \partial \left(
 \varphi_1 \wedge \partial \varphi_2 \wedge \partial \varphi_3 \right) + 
\partial (\varphi_1 \wedge \partial \varphi_2 \wedge \overline \partial \varphi_3).
$$
So it is also a complete derivative. So 
integrating over $X$ we get zero. The lemma is proved.

\vskip 3mm
{\bf Proof of Proposition \ref{17:24i}}. Applying  
$\stackrel{(-)}{\partial}$ 
to a current $g(x,y)$, we 
think about it as a sum of two operators, 
$\stackrel{(-)}{\partial_x} + \stackrel{(-)}{\partial_y}$, 
acting on $x$ and $y$ parts of $g(x,y)$. 

We start from integral (\ref{intetaq}). 
Take an external edge $E$ of the tree $T^0$. Let $v$ be  
its external vertex (or one of them if $T^0$ has a single edge). 
Then by Lemma \ref{17:24}(i), or using 
$\int_{X }\alpha_1\wedge \alpha_2 \wedge \stackrel{(-)}{\partial} \gamma =0$, 
the integral over the copy of $X $ assigned to $v$ is zero 
for the component  $\stackrel{(-)}{\partial_x}g(x,y)$ 
where $x$ runs through the copy of $X$ labeled by $v$. 
So we have only the component $\stackrel{(-)}{\partial_y}g(x,y)$. 
Using Lemma \ref{17:24} we run a similar argument, until we  
find an internal vertex of $T^0$ 
such that the corresponding integral 
assigned to this vertex is of type (\ref{nab}), 
and thus vanishes by Lemma \ref{17:24}(ii).  

Now let us consider (\ref{intetaqq}). Let $E$ be the edge 
contributing the Laplacian $\partial \overline \partial g_E$. 
Cutting this edge, we get two trees $T_E^\pm$. 
Denote by $v_\pm$ the vertex of the tree edge $E$ entering the tree $T_E^\pm$. 
Assume that 
the edge $F$ contributing just $g_F$ is in the tree $T_E^-$. 
Then if the current $g_E$ was hit by 
$\stackrel{(-)}{\partial_{v_+}}$, the integral related to $T_E^+$ is zero. 
This is proved just like we proved that 
(\ref{intetaqq}) is zero: we start from the integral related 
to the vertex $v_+$. Otherwise the current $g_E$ was hit by 
the Laplacian ${\overline \partial_{v_-}}{\partial_{v_-}}$. In this case 
the integral for the vertex $v_+$ is equivalent, by the Stokes formula, 
to the integral (ii)  in (\ref{nab}). 
Proposition \ref{17:24i} is proved.

\vskip 3mm
\bl \la{6.16.08.1}
The Hodge  correlator integral for a cyclic word 
${\cal C}(\theta_0 \otimes \ldots \otimes \theta_m)$ 
where $d\theta_i = d^\C\theta_i =0$, and 
$\theta_i = dd^\C \eta_i$ for at least two 
$\theta_i$'s, is zero. 
\el

{\bf Proof}. Suppose that $\theta_i = dd^\C  \eta_i$ 
for $i=a,b$. Using the Stokes formula, we move 
$dd^\C $ from  $\theta_a$ to the  Green current 
assigned to an edge $E$ of a tree $T$.  Formula (\ref{7.3.00.1mn}) for 
$dd^\C G_E$ gives us  the usual three terms: 
the Casimir, the volume and the $\delta$-terms. 
The last vanishes after taking the sum over all trees. The volume term always vanishes. So we are left with the Casimir term. 
Running the same argument for $dd^\C\eta_b$  
we get zero since there are no edges $F$ contributing $G_F$  rather then $d^\C G_F$  
to the operator $\omega$ left. 
The lemma is proved. 
\vskip 3mm

\subsection{Proof of Theorem \ref{3/11/07/101q}} \la{hf4.1ref}

 {\it Proof of (i)}  The degree is zero by Lemma \ref{5.17.07.5}. 
The second claim is equivalent to the shuffle 
relations for the Hodge correlators from Proposition \ref{7.3.06.4}. 
The part (i) is proved.

\vskip 3mm
\noindent
{\it Proof of (ii)}. Recall the cyclic Casimir element ${\cal C}{\rm Id}_X$, see 
(\ref{CYCCAS}).  
\begin{proposition} \label{9.24.07.71} One has 
\be \la{7.18.08.10}
\delta  {\rm Cor}_{\cal H}\Bigl( {\cal C}{\rm Id}_X \Bigr) = 
\sum_{T} \int_{X ^{\mbox{{\rm vert}$(T^0)$}}} d \kappa_T \Bigl( ~{\cal C}{\rm Id}_X \Bigr) = 0.
\ee
\end{proposition}

{\bf Remark}. Let $\gamma_0, \ldots, \gamma_m$ be harmonic representatives 
of 
cohomology classes in $H^{>0}(X)$ such 
that $\sum_i ({\rm deg}(\gamma_i)-1) = 2{\rm dim}X-3$. 
Then the  current 
$
\kappa_T(\gamma_1 \otimes  \ldots\otimes  \gamma_m)
$  has degree one less then 
the maximal possible. So (\ref{7.18.08.10}) just means that 
\begin{equation} \label{9.24.07.7}
\delta {\bf G} = \sum_{W}\frac{1}{{\rm Aut}(W)|}
\sum_T\int_{X ^{\mbox{{\rm vert}$(T^0)$}}}
d\kappa_T(\gamma_0|h_{\gamma_0} \otimes  \ldots\otimes  \gamma_m|h_{\gamma_m}) = 0. 
\end{equation}
The first sum is over a basis $W= {\cal C}(\gamma_0\otimes \ldots\otimes \gamma_m)$ 
in the degree 
$-1$ subspace of 
   ${\cal C}_{\Bbb H^*} \otimes {\cal H}$.

\vskip 3mm
{\bf Proof of Proposition \ref{9.24.07.71}}. Let us calculate the double sum  before the integration. 
We use the same terminology  and method 
as in Theorem 6.8 in \cite{G1}. We employ formula 
(\ref{5-20.1}) 
for $dw$, and use (\ref{7.3.00.1mn}) to calculate $(2\pi i)^{-1}
\overline \partial \partial G(x,y)$. 
There are three terms in the latter formula, called the $\delta$-term, the 
volume term and the Casimir term. 
The $\delta$-term vanishes just as in proof of Theorem 6.8. 
The volume term vanishes 
due to the absence of $H^0$-factors
 in the decoration and Lemma \ref{5.17.07.5}.  
By 
Lemma \ref{7.24.07.1}, the Casimir term for an edge $E$ can be non-zero only if 
 $E$ is an external edge of $T^0$. 

Let us assume that $E = E_0$ is an external edge of $T^0$, and  
$F_0$, $F_1$ are the two external edges incident to a vertex $v$ of $E_0$. 
We assume that ${\rm sgn}({\bf E}_0 \wedge {\bf F_0} \wedge {\bf F_1})=1$. 
Let 
$T'$ be the tree obtained by removing the edges  
${F}_0, {F}_1$, and   ${\rm Or}_{T'}$ 
is an 
orientation of $T'$.  
 Then  
\be \la{3.31.08.1}
{\rm sgn}({\rm Or}_{T'}) = 
{\rm sgn}({\rm Or}_{T'}\wedge {\bf F}_0 \wedge {\bf F}_1).
\ee 

Let us calculate the contribution of the 
Casimir term assigned to the edge $E_0$. 
We assume  that the external edges $F_i$ are decorated by $\gamma_i$, 
and the internal edges are  $E_0, ..., E_k$. 

\bl \la{CASCONT} The contribution of 
Casimir term for the edge $E_0$, after summation 
over $\gamma_0, \gamma_1$, is 
\be \la{4.25.08.11}
\sum_{\gamma_i}
{\rm Cor}_{\cal H}\Bigl(\delta(\alpha| h) 
\otimes \gamma_2|h_{\gamma_2} \otimes ... \otimes \gamma_m|h_{\gamma_m}\Bigr).
\ee
\el

{\bf Proof}. 
The element  ${\cal C}(\overline h_{\gamma_0} 
\otimes \ldots \otimes \overline h_{\gamma_m})$ 
appears in (\ref{9.24.07.7}) with the coefficient obtained 
by taking  
 the sum over trees of the integrals of 
$$
d\omega \Bigl((G_{E_0} \wedge {\bf E}_0) \wedge 
\ldots \wedge (G_{E_k} \wedge {\bf E}_k)\Bigr)
 \wedge 
(\gamma_0 \wedge {\bf F}_0)|\overline h_{\gamma_0}  \otimes 
(\gamma_1 \wedge{\bf F}_1)|\overline h_{\gamma_1} \otimes 
 \ldots.
$$ 
Observe that formula (\ref{5-20.1}) just means that 
$$
d\omega \Bigl((\varphi_0 \wedge {\bf E}_0) \wedge \ldots \wedge 
(\varphi_k \wedge {\bf E}_k)\Bigr)
 = 
$$
\be \la{DEQW}
- (-1)^{{\rm deg}\varphi_0}(\overline \partial \partial \varphi_0  \wedge {\bf E}_0) 
\wedge \omega \Bigl((\varphi_1 \wedge {\bf E}_1) \wedge \ldots \wedge 
(\varphi_k \wedge {\bf E}_k)\Bigr) + \ldots 
\ee 
Here $...$ means sum of the similar terms for $\varphi_i$, $i>0$, plus the 
two summands (\ref{5-20.1a}). 

The Green current has an even degree. So taking into account 
 differential equation (\ref{7.3.00.1d}) for the Green current, 
we conclude that the  
contribution of the Casimir term for the edge $E_0$ equals 
\be \la{DEQW1}
\sum_s(\alpha^{\vee}_s\wedge \alpha_s) \wedge {\bf E}_0 \bigwedge 
\omega \Bigl((G_{E_1} \wedge {\bf E}_1) \wedge \ldots \wedge (G_{E_k} \wedge {\bf E}_k)\Bigr) 
 \bigwedge 
(\gamma_0 \wedge {\bf F}_0)|\overline h_{\gamma_0} \otimes
(\gamma_1 \wedge {\bf F}_1)|\overline h_{\gamma_1} \otimes
 \ldots .
\ee
Notice that the $-$ signs in (\ref{DEQW}) and in front of the Casimir term in 
(\ref{7.3.00.1d}) 
canceled each other. Further, the form $\alpha^{\vee}_s\wedge \alpha_s$ is even, 
and $\omega(...)$ in (\ref{DEQW1}) is odd. 
So moving $\alpha^{\vee}_s\wedge \alpha_s \wedge {\bf E}_0$ 
we get 
$$
-\omega\Bigl((G_{E_1} \wedge {\bf E}_1) \wedge \ldots \wedge (G_{E_k} \wedge {\bf E}_k)\Bigr)
 \bigwedge 
$$
$$
(\gamma_0 \wedge {\bf F}_0)|\overline h_{\gamma_0}  
\otimes (\gamma_1 \wedge {\bf F}_1)|\overline h_{\gamma_1} 
 \bigwedge \sum_s 
(\alpha^{\vee}_s \wedge 
\alpha_s) \wedge {\bf E}_0 \bigwedge (\gamma_2 \wedge {\bf F}_2)|\overline h_{\gamma_2} \otimes
 \ldots .
$$

\bl \la{7.11.08.10}
$$
\int_X
\sum_{\gamma_0, \gamma_1}
(\gamma_0 \wedge {\bf F}_0)|\overline 
h_{\gamma_0} \otimes (\gamma_1 \wedge {\bf F}_1)|\overline h_{\gamma_0} \bigwedge \sum_s 
(\alpha^{\vee}_s \wedge 
\alpha_s) \wedge {\bf E}_0  = 
- \delta( \alpha|h) 
\cdot {\bf E}_0 \wedge {\bf F}_0 \wedge {\bf F}_1. 
$$
\el

{\bf Proof}. Let ${\rm pr}_{\cal H}$ be the orthogonal 
projection onto the space of harmonic forms. Then  
\begin{equation} \label{9.24.07.7s}
{\rm pr}_{\cal H}(\beta) = \sum_s (\beta, \alpha_s^{\vee}) \alpha_s. 
\end{equation}
Thus  $
\sum_s\int_X (\gamma_0 \wedge  \gamma_{1}\wedge \alpha^{\vee}_s)\cdot  
\alpha_s = {\rm pr}_{\cal H}(\gamma_0 \wedge  \gamma_{1})$. The Lemma follows from (\ref{9.29.07.5m}). 

\vskip 3mm
This 
proves Lemma \ref{CASCONT}. Notice that the two $-$ signs in the last two formulas canceled 
each other. 
Thus we proved Proposition \ref{9.24.07.71} 
and hence the part (ii).

\vskip 3mm
\noindent
{\it Proof of (iii)}. Recall that there are the 
following choices  in the definition of the 
Green current: 

\begin{enumerate}

\item
 A solution of equation (\ref{7.3.00.1mn}) for the Green current  $G$ 
is unique up to  

(a) currents of type $\partial \alpha +  \overline \partial \beta$. 

(b) harmonic forms of type $(n-1, n-1)$. 

\item
Differential equation (\ref{7.3.00.1mn}) 
depends on the splitting  
of the Dolbeaut complex. 

\end{enumerate}

Let us investigate how these choices  affect the resulting Hodge 
correlator ${\bf G}$. 

\vskip 3mm

1a) 
The difference of the Hodge correlators 
of ${\cal C}(\alpha_1 \otimes \ldots \otimes \alpha_m)$ defined using 
$G(x,y)$ and $G(x,y) + \stackrel{(-)}{\partial}\gamma$   
is  a linear combination of  integrals like (\ref{intetaq}) 
where $g_{i}$ stands for $G_E$ or $\gamma$. 

If the tree $T^0$ is empty, the difference is zero by the very definition. 
If the tree $T^0$ is not empty, it is zero by Proposition \ref{17:24i}.

\vskip 3mm

1b) Let us add a closed form $\gamma$ to a Green current. 
Using the Kunneth formula for $H^*(X\times X)$ and 1a) we may assume
that it is decomposable: $\gamma = \omega_1 \times \omega_2$, $d\omega_i=0$. 
Then its contribution is zero by Lemma \ref{corT1}, 
unless the tree $T^0$ has a single vertex. 
To handle the latter case we use Lemma \ref{7.24.07.1}. 
Thus if the tree $T^0$ has a single vertex,
 the contribution is shown on Fig \ref{hf15}. It equals
$
{\cal C}\left(\delta h_{\omega_1} \otimes \delta h_{\omega_2} \right). 
$ 
This is obviously a coboundary. The part 1. of (iii) is proved.  

\vskip 3mm
\begin{figure}[ht]
\centerline{\epsfbox{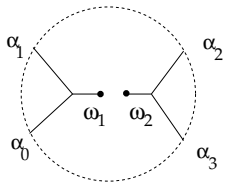}}
\caption{The contribution of $\gamma = \omega_1 \times \omega_2$.}
\label{hf15}
\end{figure}

\vskip 3mm

2. We want to compare two splittings $s$ and 
 $\widetilde s$ of the Dolbeaut complex. We assume that they differ 
only by choice of another representative $\widetilde \omega^{\vee}$ 
of a $(p, q)$-harmonic form $\omega^{\vee}$, which is dual to 
a harmonic form $\omega$.  
By the $\overline \partial\partial$-lemma 
there exists a $(p-1, q-1)$-form $\beta$ such that 
$\overline \partial\partial\beta = 
\widetilde \omega^{\vee} - \omega^{\vee}$. 
Let $\widetilde G = \widetilde G(x,y)$ be a Green current 
for the splitting $\widetilde s$. 
One can choose it as 
\begin{equation} \label{3.17.07.2}
\widetilde G(x,y) = G(x,y) + p_1^*\omega
\wedge p_2^*\beta, 
\end{equation}
so the second term is decomposable. 
Let ${\rm Cor}^{\widetilde G}_{\cal H}$ 
be the  
Hodge correlator map for the 
Green current $\widetilde G$. 
Let us calculate 
the difference $\widetilde {\bf G} - {\bf G}$ 
between the Hodge correlators for 
the Green currents $\widetilde G(x,y)$ and $G(x,y)$.

\vskip 3mm
The integrand of the Hodge correlator is obtained 
by applying the polydifferential operator  $\omega$  to the 
(pull backs of) Green currents assigned to the edges of $T^0$, and multiplying the result by the closed forms assigned to the external edges of $T$.  
The difference $\widetilde {\bf G} - {\bf G}$ is a sum of such 
integrals, where one or more of the Green 
currents $G(x,y)$ are replaced by $p_1^*\omega \wedge p_2^*\beta$. 
Let us pick an edge $E$ of $T$ 
where  $p_1^*\omega \wedge p_2^*\beta$ was employed. 
By Lemma \ref{7.24.07.1}
the contribution of such a pair  $(E, T)$ 
to $\widetilde G = \widetilde G(x,y)$ is zero 
unless the condition (\ref{9.24.07.5}) holds.

\vskip 3mm
Consider expressions of the following types ($s_{\omega}$ was introduced in 
the end of Section 5.1):
\begin{equation} \label{10.7.09.2}
(i) ~ \delta (d^\C\beta | s_\omega), \qquad 
(ii)~  dd^\C \beta | s_\omega, \qquad 
(iii) ~  \alpha_i | h_{\alpha_i}. 
\end{equation} 
We make   cyclic words 
out of them. Since $\omega, s_\omega, \beta$ have the same parity 
of the grading, these expressions 
are plain cyclic invariant. 
Let us introduce a notation
$$
\Delta (\gamma | h): =  d \gamma| h +
(-1)^{|\gamma|} \gamma | \delta h. \qquad \mbox{So} \quad  \Delta^2 =0.
$$
It agrees with  (\ref{thatf}) since there we have $d \gamma=0$.  

\vskip 3mm
Denote by ${\rm Cor}_{\cal H}^{\ast}$ the Hodge correlator map 
defined by using the polydifferential operator  $\xi$ instead of 
$\omega$, see (\ref{5-x20.1s}). The relationship between ${\rm Cor}_{\cal H}^{\ast}$ and 
${\rm Cor}_{\cal H}$ 
is the following: 
$$
{\rm Cor}^{\ast}_{\cal H} = \mu \circ {\rm Cor}_{\cal H} \circ \mu^{-1}.  
$$
Here $\mu$ is the operator multiplying a 
homogeneous element of bidegree $(p,q)$ by $p+q-2\choose p-1$.
Set 
$$
{\bf A}^\ast := {\bf A'}^\ast + {\bf A''}^\ast := 
\sum_Z{\rm Cor}^{\ast, G}_{\cal H}({\cal C}(Z)) ~+~ 
{\rm Cor}^{\ast, G}_{\cal H}(\alpha|h \otimes 
\alpha|h \otimes \omega) 
{\rm Cor}^{\ast, G}_{\cal H}(\beta \otimes \alpha|h \otimes 
\alpha|h). 
$$ 
Here $Z$ runs through a basis in the space 
of all cyclic products 
of expressions $\Delta (d^\C \beta | s_w)$ and $\alpha_s|h_{\alpha_s}$. 
The second term in ${\bf A}^\ast$
corresponds to the diagram on Fig. \ref{hf15} where $\omega_1=\omega$ and $\omega_2 = \beta$.

\bl \la{10.17.07.1}
One has the following, where $n$ is the number of factors ${\rm Id}_X$:

\be \la{AAA}
{\bf A'}^* = \sum_{n\geq 1} 
\frac{1}{n} {\rm Cor}^{\ast, G}_{\cal H}\Bigl(
\Delta (d^\C \beta | s_w)
 \otimes {\rm Id}_X \otimes \ldots \otimes 
\Delta (d^\C \beta | s_w)
 \otimes {\rm Id}_X\Bigr). 
\ee
\el

{\bf Remark}. The coefficient $n$ is the order of the 
automorphism group of the summand. 

\vskip 3mm
{\bf Example}. Choose a basis $\{A_i\}$ 
in the tensor algebra of the vector space with the 
basis $\{\alpha_k|h_{\alpha_k}\}$. 
 The beginning of ${\bf A'}^*$  looks as follows: 
$$
\sum_{i}{\rm Cor}^{\ast, G}_{\cal H}\Bigl(\Delta (d^\C \beta | s_w)
\otimes A_i\Bigr) + 
\frac{1}{2} \sum_{i \not = j}{\rm Cor}^{\ast, G}_{\cal H}
\Bigl(\Delta (d^\C \beta | s_w) \otimes A_i \otimes 
\Delta (d^\C \beta | s_w) \otimes A_j\Bigr) + 
$$
$$
\frac{1}{2} \sum_{i}{\rm Cor}^{\ast, G}_{\cal H}\Bigl(
\Delta (d^\C \beta | s_w)
 \otimes A_i \otimes 
\Delta (d^\C \beta | s_w) \otimes A_i\Bigr) + \ldots 
$$
Here the coefficient $1/2$ in front of the second term appears since 
the cyclic word for the pair $(A_i, A_j)$ is the same as for the one  
$(A_j, A_i)$. The coefficient $1/2$ in front of the next term appears 
for a different reason: 
the corresponding cyclic word has an automorphism group of order two.

\vskip 3mm
{\bf Proof of Lemma \ref{10.17.07.1}}. Given a cyclic word 
$$
Z = {\cal C}\Bigl(\Delta (d^\C \beta | s_w)
 \otimes A_{1} \otimes \ldots \otimes 
\Delta (d^\C \beta | s_w) \otimes A_{n}\Bigr),
$$ write $n=ab$ where 
$\Z/a\Z$ is the automorphism group of $Z$ -- so $A_{i+a} = A_i$.  
Then $Z$ enters with the coefficient $1/a \cdot 1/b = 1/n$.
Indeed, the coefficient $1/a$ appears since $a = |{\rm Aut}(Z)|$, and 
the coefficient $1/b$ is needed  since the cyclic 
shift which moves $A_k$ to $A_{k+i}$, where $0 \leq i \leq b-1$,   
produces a cyclic word  
which equals to $Z$. The lemma is proved. 
\vskip 3mm

\begin{proposition} \label{10.7.09.1}
One has $\widetilde {\bf G}^* - {\bf G}^* = {\bf A}^*$. 
\end{proposition}

{\bf Proof}. 
1) Denote by ${\cal C}\widetilde {\rm Id}_X$ the cyclic Casimir element written by using the 
forms $\widetilde \alpha_s$ instead of $\alpha_s$. 
Since $\widetilde \omega^{\vee} - \omega^{\vee} = dd^\C\beta$, we have 
\be \la{ww1w}
{\rm Cor}^{\ast, \widetilde G}_{\cal H} 
\Bigl({\cal C}\widetilde {\rm Id}_X - 
{\cal C}{\rm Id}_X\Bigr) =   
\sum_{X}{\rm Cor}^{\ast, \widetilde G}_{\cal H}
({\cal C}(X)). 
\ee 
The sum on the right of (\ref{ww1w}) is over a basis $X$ in the cyclic tensor algebra 
of the vector space with a basis $dd^\C\beta|s_\omega$, $\{\alpha_s|h_{\alpha_s}\}$. 
The equality is proved by an argument similar to the one in the proof of Lemma \ref{10.17.07.1}. 


2) Let us calculate 
\be \la{po}
\sum_{X}
\Bigl({\rm Cor}^{\ast, \widetilde G}_{\cal H}({\cal C}(X)) - 
{\rm Cor}^{\ast, G}_{\cal H}({\cal C}(X))\Bigr).  
\ee 
By Lemma \ref{7.24.07.1} 
we can get non-zero contributions 
only from an external edges $E$ of  $T^0$ such that

{\it there is a vertex of $E$ shared  by 
two external edges, $F_-$ and $F_+$, 
and these edges are decorated by $\alpha$'s, 
not by $dd^\C\beta$ -- we denote them by $\alpha_{-}$ and $\alpha_{+}$, see Fig \ref{hf16}}.
\begin{figure}[ht]
\centerline{\epsfbox{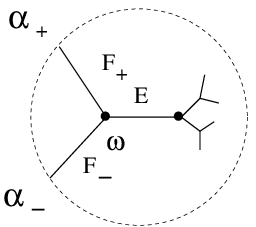}}
\caption{}
\label{hf16}
\end{figure}

If the tree $T^0$ has a single vertex,  
we get the second term in ${\bf A}^\ast$. 

\vskip 3mm
\bl \la{7.15.08.1}
Assume that the tree $T^0$ has more then one vertex.
 Then the contribution to (\ref{po}) 
of an external edge $E$ of the tree $T^0$ as above 
is described 
by the following procedure:

\begin{itemize} 
\item 
Cut the edge $E$ and throw away the small tree $T''_E$, 

\item 
Decorate 
the new external edge $E'$ by $\delta(d^\C \beta|  s_\omega)$. 
\end{itemize}
\el

The same recipe is applied to every edge $E$ where 
$G_E$ is replaced by  $\widetilde G_E$. 
Using the same argument as 
in the proof of Lemma \ref{10.17.07.1}, we 
get 
\be \la{rerer}
(\ref{po}) = \sum_{Y}{\rm Cor}^{\ast, \widetilde G}_{\cal H}
({\cal C}(Y))
\ee
where $Y$ is a basis in the cyclic tensor algebra  
of the vector space spanned by $\delta(d^\C\beta| s_{\omega})$ 
and $\{\alpha_s|h_{\alpha_s}\}$. 
Proposition \ref{10.7.09.1} would follow from this by using the 
Lemma \ref{10.17.07.1} 
argument again. 

\vskip 3mm
{\bf Proof of Lemma \ref{7.15.08.1}}. It is 
similar to the proof of Proposition \ref{9.24.07.71}. 
The key point is this:\footnote{This is the reason to use 
the version ${\rm Cor}^*$ of the Hodge correlator, see Remark at the end of this Subsection.}   
\be \la{7.16.08.10}
\int\xi(\widetilde G_{E_{0}}\wedge \ldots \wedge \widetilde G_{E_{k}} )\wedge \ldots  = 
\int\xi(\widetilde G_{E_{0}}\wedge \ldots \wedge \widetilde G_{E_{k-1}} )
\wedge d^\C\widetilde G_{E_{k}}\wedge \ldots . 
\ee
We arrange the edges so that $E= E_k$, one has 
 ${\rm sgn}({\bf F}_- \wedge {\bf F}_+ \wedge {\bf E}) = 1$, 
and $F_{\pm}$ is decorated by $\alpha_{\pm}$. Lemma \ref{7.15.08.1} 
boils down to the computation of 
$$
(d^\C(\widetilde G_{E} - G_{E}) \wedge {\bf E}) \bigwedge (\alpha_- \wedge {\bf F}_-)|
\overline h_{\alpha_{-}} 
\otimes (\alpha_+ \wedge {\bf F}_+)|\overline h_{\alpha_{+}} 
$$
Using  (\ref{4.21.08.1a}), 
we prove Lemma \ref{7.15.08.1}.  
Proposition \ref{10.7.09.1} is proved.


\vskip 3mm 
Set  
\be \la{BB} 
{\bf B} :=  \sum_{n\geq 1}
\frac{1}{n}{\rm Cor}^{G}_{\cal H}\Bigl(d^\C \beta | s_w \otimes 
{\rm Id}_{X} \otimes \Delta (d^\C \beta | s_w)
 \otimes 
{\rm Id}_{X} \otimes  \ldots \otimes  
\Delta (d^\C \beta | s_w)
\otimes {\rm Id}_{X} \Bigr). 
\ee 
Here there are $(n-1)$ factors 
$\Delta (d^\C \beta | s_w) \otimes {\rm Id}_{X}$. 
We denote by ${\bf A}'$ the expression similar to ${\bf A}^{' *}$ 
defined using the Hodge correlator ${\rm Cor}$ instead of ${\rm Cor}^*$. 

\begin{proposition} \label{10.7.09.1df}
One has ${\bf A}'= \delta{\bf B}$. 
\end{proposition}

{\bf Proof}.
Applying the differential $\delta$ to 
$d^\C\beta|s_w$ in ${\bf B}$ we get  a part of the sum 
${\bf A}^{'}$:
$$
{\bf A}^{'}_1:= \sum_{n>0}{\rm Cor}_{\cal H}
\Bigl(d^\C\beta | \delta s_{\omega} \otimes {\rm Id}_{{X}} \otimes 
\Delta (d^\C \beta | s_w) \otimes 
{\rm Id}_{{X}} \otimes \ldots\Bigr).
$$ 

\bl \la{7.16.08.1} Applying  $\delta$ 
to every factor $h_{\alpha_i}$ in 
${\bf B}$ 
except $s_\omega$,  and  taking the sum 
we get 
$$
{\bf A}^{'}_2:= \sum_{n>0}{\rm Cor}_{\cal H}
\Bigl(dd^\C\beta | s_{\omega} \otimes {\rm Id}_{{X}} \otimes 
\Delta (d^\C \beta | s_w) \otimes 
{\rm Id}_{{X}} \otimes \ldots\Bigr).
$$
\el

\vskip 3mm
{\bf Proof}. The claim we calculate the following sum, 
similar to the one defining ${\bf B}$: 
$$
0= \sum_{n>0}\int \sum_Td\kappa_T
\Bigl(d^\C\beta | s_{\omega} \otimes {\rm Id}_{{X}} \otimes 
\Delta (d^\C \beta | s_w) \otimes 
{\rm Id}_{{X}} \otimes \ldots\Bigr).
$$
Let us calculate  
$d\kappa_T(d^\C\beta|s_\omega \otimes {\rm Id}_{{X}} \otimes \ldots)$. 
Applying $d$ to $d^\C\beta|s_\omega$ we 
get $dd^\C\beta|s_\omega$. 
The other  
decorating forms are $d^\C$-closed. 
Thus differentiating the  Green currents 
in $\kappa_T(d^\C\beta|s_\omega 
 \otimes {\rm Id}_{{X}} \otimes \ldots)$, 
we get, similarly to Lemma \ref{CASCONT}, an 
expression containing the orthogonal 
projections onto the harmonic forms of products of three kinds: 
$d^\C \beta \wedge \alpha_\ast$, $
dd^\C \beta \wedge \alpha_\ast$, or 
$\alpha_{i} \wedge \alpha_{i+1}$. 
The first two of them are zero since we integrate complete differentials.
For example, using $d^\C \alpha_j=0$, 
$$
{\rm pr}_{\cal H}
(d^\C \beta \wedge \alpha_\ast) = 
\int_{X }d^\C \beta\wedge  \alpha_\ast \wedge  \alpha_k = 
\int_{X }d^\C (\beta\wedge  \alpha_\ast \wedge  \alpha_k) = 0.  
$$
The sum of remaining integrals 
coincides with the negative (because $d^\C \beta | s_{\omega}$ is odd) of the 
expression obtained by applying the differential 
$\delta$ to every factor $h_{\alpha_i}$ in $\alpha|h_{\alpha_i}$ 
except $s_\omega$, and taking the sum.  
Lemma \ref{7.16.08.1}, hence Proposition \ref{10.7.09.1df}, and 
thus the part 2) of (iii) are proved. 

\vskip 3mm
Clearly ${\bf A}^{'} = {\bf A}^{'}_1+ {\bf A}^{'}_2$.  
Furthermore, ${\bf A}'= \delta{\bf B}$ implies ${\bf A}^{'*}= \delta{\bf B}^*$ 
for the element ${\bf B}^*$ defined using ${\rm Cor}^*$. 
Thus Theorem \ref{3/11/07/101q} is proved.

\vskip 3mm 
{\bf Remark}. To write down a formula for the coboundary, 
we employ two forms of the Hodge correlator: the one ${\rm Cor}_{\cal H}$ using 
the polydifferential operator $\omega$, and the one ${\rm Cor}^*_{\cal H}$ using $\xi$. 
Dealing with $dd^\C G_E$ it is convenient to use $\omega$, 
while handling the effect of replacing $\widetilde G_E$ by $G_E$ 
it is handy to employ $\xi$, because of formula (\ref{7.16.08.10}).

\vskip 3mm

\subsection{Functoriality: Proof of Theorem 
\ref{PPPq}} \la{hf4.2ref}

\vskip 3mm
We work in the category of 
compact complex manifolds. Let $f: X \to Y$ be a map of manifolds. 
There are two standard functors: 
the pull back 
$f^*: H^*(Y) \to H^*(X)$,  and the dual map $f_*: H_*(X) \to H_*(Y)$. 
We define a map $f^!: H_*(Y) \to H_*(X)$ by using the Poincare duality:
\be \la{6.17.08.1}
f^!(c) \cap {\cal H}_X:= f^*(c \cap {\cal H}_Y).
\ee
Therefore one has 
\be \la{6.17.08.4}
\langle c_0 \cap {\cal H}_Y \cap f_*h\rangle_Y = 
\langle f^*(c_0 \cap {\cal H}_X) \cap h\rangle_Y \stackrel{(\ref{6.17.08.1})}{=} 
\langle f^!c_0 \cap {\cal H}_X \cap h\rangle_X. 
\ee
Let $f_!: H^*(X) \to H^*(Y)$ be the dual map. There is the projection formula
\be \la{PROJJ}
f_!(\alpha \cup f^*\gamma) = f_!\alpha \cup \gamma.  
\ee

We reserve the notation $f^*$ for the inverse 
image of smooth forms, and $f_!=f_*$ for 
the direct image of distributions -- the adjoint map. 
They are well defined in the category of manifolds. 


We denote by $f^!$ the inverse image of distributions. 
It is not always defined: a sufficient condition is that 
the conormal bundle to the graph of $f$ is transversal 
to the wave front of the distribution. 
\vskip 3mm

A map $f: X \to Y$   
induces a map 
$
f_*: {\Bbb H}_*(X) \to 
{\Bbb H}_*(Y),
$ 
 and hence a push forward 
$f_*{\bf H}_{X, x}$.  
\begin{lemma} \label{6.13.08.1} The induced map $f_*:
 {\cal C}{\cal L}ie_{{\Bbb H}_*(X)} \lra  {\cal C}{\cal L}ie_{{\Bbb H}_*(Y)}$ 
is a map of complexes. 
\end{lemma}

{\bf Proof}. Given a finite dimensional vector space $V$, denote by 
${\rm Id}_V \in V^* \otimes V$ the Casimir element corresponding to the identity map. 
Given a linear map $f: V \to W$, let $f^*:W^* \to V^*$ be the dual map. Set
$$
 F_* = (f, {\rm Id}): V^* \otimes V \lra V^* \otimes W, 
\quad F^* = ({\rm Id}, f^*): W^* \otimes W \lra V^* \otimes W. 
$$
Then one has 
\be \la{FCASIMIR}
F_*{\rm Id}_V = F^*{\rm Id}_W.
\ee

\vskip 3mm

Let $\delta_X$ (respectively $\delta_Y$) be the differential  
in ${\cal C}{\cal L}ie_{{\Bbb H}_*(X)}$ (respectively ${\cal C}{\cal L}ie_{{\Bbb H}_*(Y)}$). 
We denote by $\{\alpha = \alpha_{(s)}\}$ a basis in $s_X{\Bbb H}^*(X)$, by 
$\{a = a_{(s)}\}$ the dual basis in 
${\Bbb H}_*(X)$, and by $\{\gamma\}$ and $\{c\}$ a basis and the dual basis 
in $s_Y{\Bbb H}^*(Y)$ and ${\Bbb H}_*(Y)$. 

\vskip 3mm
{\it Conventions.} 1. We assume that 
if both $\alpha$ and $a$ enter to a formula, 
we have the corresponding Casimir element 
$\alpha| a:= \sum_s \alpha_{(s)}\otimes a_{(s)}$ there.  Similarly with $\gamma$ and $c$. 

2. We often ignore signs in the Hodge correlator integrals 
since they were carefully explained 
in the Section \ref{hf4.1ref}. 

\vskip 3mm
Given a generator $h \in {\Bbb H}_*(X)$ of the Lie algebra 
${\cal L}ie_{{\Bbb H}_*(X)}$, one has 
$$
f_*\delta_X(h) = 
\sum_{{\cal C}(\alpha_0\otimes \alpha_1\otimes \alpha_2)}
\int_X(\alpha_0\wedge \alpha_1\wedge \alpha_2) 
\langle a_2 \cap {\cal H}_X \cap h\rangle_X \cdot 
f_*a_0 \otimes f_*a_1.
$$
$$
\delta_Y(f_*h) = 
\sum_{{\cal C}(\gamma_0\otimes \gamma_1\otimes \gamma_2)}
\int_Y(\gamma_0\wedge \gamma_1\wedge  \gamma_2) 
\langle c_2 \cap {\cal H}_Y \cap f_*h\rangle_Y \cdot 
c_0 \otimes c_1.
$$
Here and below the sum is over bases 
${\cal C}(\alpha_0\otimes \alpha_1\otimes \alpha_2)$, 
${\cal C}(\gamma_0\otimes \gamma_1\otimes \gamma_2)$. 
One has 
$$
\delta_Y(f_*h) \stackrel{(\ref{6.17.08.4})}{=}
\sum_{}
\int_Y(\gamma_0\wedge \gamma_1\wedge  \gamma_2) 
\langle f^!c_2 \cap {\cal H}_X \cap h\rangle_X \cdot 
c_0 \otimes c_1. 
$$
Applying the Casimir elements identity $\gamma_2| f^!c_2 = f_!\alpha_2| a_2$, see 
(\ref{FCASIMIR}), 
we write this as 
$$
\sum_{}
\int_Y(\gamma_0\wedge \gamma_1\wedge  f_!\alpha_2) 
\langle a_2 \cap {\cal H}_X \cap h\rangle_X \cdot 
c_0 \otimes c_1 \stackrel{(\ref{PROJJ})}{=} 
$$
$$
\sum_{}
\int_X(f^*(\gamma_0\wedge  \gamma_1)\wedge \alpha_2) 
\langle a_2 \cap {\cal H}_X \cap h\rangle_X \cdot 
c_0 \otimes c_1 \stackrel{(\ref{FCASIMIR})}{=} 
$$
$$
\sum_{}
\int_X(\alpha_0\wedge \alpha_1\wedge  \alpha_2) 
\langle a_2 \cap {\cal H}_X \cap h\rangle_X \cdot 
f_*a_0 \otimes f_*a_1 = f_*\delta_X(h). 
$$
The lemma is proved.
\vskip 3mm

Let us choose direct sum decompositions provided by the map $f^*: H^*(Y) \lra H^*(X)$:
$$
H^*(Y) = {\rm Ker}f^* \oplus \widetilde {\rm CoIm}f^*, \qquad 
 H^*(X) = {\rm Im}f^* \oplus \widetilde {\rm CoKer}f^*. 
\qquad f^*: \widetilde {\rm CoIm}f^* 
\stackrel{\sim}{\lra}{\rm Im}f^*.
$$
Using Theorem \ref{3/11/07/101q} we reduce the functoriality proof to the case 
when the splittings $s_X$ and $s_Y$ of the Dolbeaut complexes of 
$X$ and $Y$ are compatible in the following sense: 
$$
f^*\Bigl(s_Y\widetilde {\rm CoIm}f^*\Bigr)  = s_X{\rm Im}f^*. 
$$
We will assume this below. 

Denote by 
$\{\gamma''_{(j)}\}$ (respectively $\{\gamma'_{(i)}\}$) a basis in $s_Y{\rm Ker}f^*$
 (respectively $s_Y{\rm CoIm}f^*$),  and by 
$\{c''_{(j)}\}$ (respectively $\{c'_{(i)}\}$) 
the dual bases 
in the homology. In particular, $\{c'_{(i)}\}$ is 
a basis in $f_*{\Bbb H}_*(X)$. 
The basis $\{\gamma'_{(i)}\}$ 
gives rise to a basis $\{\alpha'_{(i)}:= f^*\gamma'_{(i)}\}$ in  $s_X{\rm Im}f^*$.

\vskip 3mm
The induced by $f$ maps 
$$
f_1: X \times X \lra Y \times X, \qquad    f_2: Y \times X \lra Y \times Y, 
$$
provide the 
following two currents on $Y  \times X$: 

(i) the push forward $f_{1*}G_X(y_1, x_2)$ of the 
Green current $G_X(x_1, x_2)$;

(ii) the restriction $f_2^!G_Y(y_1, x_2)$ 
of the Green current $G_Y(y_1, y_2)$. 

\noindent
To check that $f^!_2G_Y$ is defined, decompose 
$f$ into a composition of a closed embedding 
and a projection $X \times {\Bbb P}^n \to X$.

Denote by $p_Y$ and $p_X$ the projections  from $Y \times X$ onto $Y$ and $X$. 
The following Lemma is another incarnation of identity (\ref{FCASIMIR}). 

\begin{lemma} \label{9.23.07.1}  
$f_{1 \ast}G_X - f^!_{2}G_Y$  is 
a sum of decomposable forms on $Y  \times X $:
\be \la{6.16.08.4}
f_{1 \ast}G_X - f^!_{2}G_Y  ~=~  
\sum_{j}p_Y^* \gamma''^{\vee}_{(j)} \wedge p_X^*\varphi_{(j)} 
- \sum_{k}p_Y^* \psi_{(k)}\wedge  p_X^*\alpha''_{(k)} 
+  \mbox{\rm a closed form}.  
\ee
\end{lemma}

{\bf Remark}. There are three terms in this formula. The first two are of the form  
\be \la{6.16.08.4a}
\mbox{(a harmonic form) $\times$ (an arbitrary form)}  +
\mbox{(an arbitrary form) $\times$ (a harmonic form)}.  
\ee

{\bf Proof}. 
Let $\Delta_f$ be the graph of the map 
$f:X \to Y$, and ${\rm Har}_D$ is a harmonic representative 
of the cohomology class of a cycle $D \subset X \times Y$.  
The currents (i), (ii)
 satisfy differential equations
$$
(4\pi i)^{-1}dd^\C G' = 
\delta_{\Delta_f} - f_{1 \ast}{\rm Har}_{\Delta_X}, \qquad
(4\pi i)^{-1}dd^\C  G'' = 
\delta_{\Delta_f} - f_2^*{\rm Har}_{\Delta_Y}. 
$$
One has 
${\rm Ker}f_! = ({\rm Im}f^*)^{\perp}$, where the orthogonal is for
 the form on the cohomology. 
Recall the dual basis $\{\alpha_{(k)}^{\vee}\}$. 
So $\{\alpha_{(k)}''^{\vee}\}$ is a basis in ${\rm Ker}f_!$. 
By the $dd^\C$-lemma 
there exists a form $\psi_{(k)}$ such that 
\be \la{CFAsa}
f_*\alpha_{(k)}''^{\vee} = (4\pi i)^{-1}dd^\C \psi_{(k)}. 
\ee
Since the cohomology class $f^*[\gamma''_{(j)}]=0$, 
by the $dd^\C$-lemma 
there exists a form $\varphi_{(j)}$ such 
that 
\be \la{CFA}
f^*\gamma''_{(j)} = (4\pi i)^{-1}
dd^\C\varphi_{(j)}. 
\ee 
Thus   
$$
f^!_{2}{\rm Har}_{\Delta_Y} -  f_{1 \ast}{\rm Har}_{\Delta_X}
= (4\pi i)^{-1}
\Bigl(\sum_j p_Y^*\gamma''^{\vee}_{(j)}\wedge p_X^*dd^\C \varphi_{(j)} 
- \sum_k p_Y^*dd^\C \psi_{(k)}\wedge p_X^*\alpha''_{(k)}\Bigr). 
$$
The lemma follows by the $dd^\C$-lemma. 
\vskip 3mm
\begin{figure}[ht]
\centerline{\epsfbox{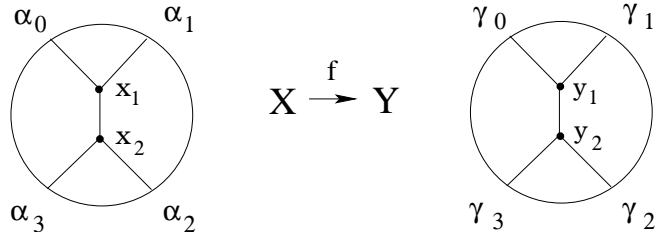}}
\caption{Proof of functoriality: the running example.}
\label{hf11}
\end{figure}

Let us denote by 
${\rm T}$ 
the tensor algebra of the vector space 
spanned by 
$$
\gamma_{(s)}|c_{(s)}  
\quad \mbox{and} \quad \delta(f_*d^\C\varphi_{(j)}|c_{(j)}'').
$$ 
We denote by ${\rm T}_k$ (respectively ${\rm T}_{\geq k}$) 
its component of the tensor degree $k$ (respectively $\geq k$). 

We denote by ${\rm T}^0$ the subalgebra of ${\rm T}$ spanned by 
$\gamma_{(s)}|c_{(s)}$. 
Recall ({\ref{CASx}) the element ${\rm Id}_{Y}$.

\bt \la{7.12.08.1} Let $\{A\}$ be a basis in ${\rm T}_{2}$, and 
 $\{B\}$ a basis in  ${\rm T}_{\geq 3}$. 
Then 
$$
f_*\delta_X(h) - \delta_X(f_*h) = 
$$
\be \la{7.12.08.3a}
\delta\sum_{m}
{\rm Cor}_{{\cal H}, Y}\Bigl({\rm Id}_{Y}
\otimes d^\C \varphi_{i}|c_i'' 
\otimes {\rm Id}_{Y} \otimes \delta (d^\C \varphi_{i}|c_i'') \otimes \ldots \otimes   
  {\rm Id}_{Y} \otimes \delta (d^\C \varphi_{i}|c_i'') \otimes 
\gamma_m\Bigr) 
\langle f^!c_m \cap {\cal H}_X \cap h\rangle
\ee
\be \la{7.12.08.3}
+ \delta\sum_{k; ~A, B}
{\rm Cor}_{{\cal H}, Y}\Bigl(A
\otimes \psi_{(k)} + B
\otimes d^\C\psi_{(k)}\Bigr) 
\langle a_{(k)}''^{\vee} \cap {\cal H}_X \cap h\rangle.
\ee
\et

{\bf Proof}. 
We call the external edge of a tree $T$ 
decorated by the form $\alpha_m$ 
{\it the special leg}, and show it by a thick edge on pictures. 
Its vertices are called {\it special vertices}. 
Every edge of $T$ is oriented
 ``towards the special edge''. 
Given a vertex $v$ of $T$, there is a unique edge $\stackrel{\to}{E_v}$ 
incident to $v$ and closest to the special edge. 

One has 
\be \la{CVA}
f_*\delta_X(h) = \sum_m{\rm Cycle}_{m+1}\sum_W
\frac{1}{|{\rm Aut}W|}{\rm Cor}_{{\cal H}, X}(\alpha_0|a_0 
\otimes \ldots \otimes \alpha_{m-1}|a_{m-1} 
\otimes \alpha_{m}) \langle a_m \cap {\cal H}_X \cap h\rangle . 
\ee 
We are going to transform the sum of integrals (\ref{CVA}) over products of copies of $X$ 
into a similar sum of integrals over products of copies of $Y$,
getting $\delta_X(f_*h) + (\ref{7.12.08.3}) + (\ref{7.12.08.3a})$.

We use the tree on Fig. \ref{hf11} decorated by 
a cyclic word ${\cal C}(\alpha_0 \otimes \ldots \otimes \alpha_3)$
as a running example.

\vskip 3mm
{\it Step 1 -- replacing $\alpha_i|a_i$ by $f^*\gamma'_i|c_i'$ for 
non-special external vertices of $T$}. 
Taking the sum over a basis 
${\cal C}(\alpha_0 \otimes \ldots \otimes \alpha_3)$, we get
$$
\sum 
{\int}_{X\times X}\Bigl(\alpha_0\wedge\alpha_1\wedge G_X(x_1, x_2) 
\wedge\alpha_2\wedge \alpha_3\Bigr) \langle a_3\cap {\cal H}_X \cap h\rangle \cdot
f_*a_0 \otimes f_*a_1 \otimes f_*a_2  = 
$$
$$
\sum_{}
{\int}_{X\times X}\Bigl(f^*\gamma'_0\wedge f^*\gamma'_1\wedge G_X(x_1, x_2) 
\wedge f^*\gamma'_2\wedge \alpha_3\Bigr) \langle a_3\cap {\cal H}_X \cap h\rangle \cdot
c'_0 \otimes c'_1 \otimes c'_2. 
$$
The first sum is  over $a_i'$ only, 
since $a''_i$ is in ${\rm Ker}f_*$. Since 
$\alpha'_i = f^*\gamma'_i$, we get the second sum. 

The general case is completely similar. 

\vskip 3mm
{\it Step 2 --  
replacing the Casimir factor 
$f^*\gamma'_i | c_i'$ by the one $f^*\gamma_i | c_i$}.  
By (\ref{CFA})
$$
f^*\gamma_i | c_i = f^*\gamma'_i | c_i' + 
(4\pi i)^{-1}dd^\C \varphi_{i} | c''_{i}, \qquad 
dd^\C\varphi_{i} | c''_{i}:= 
\sum_s dd^\C \varphi_{i, (s)} \otimes c''_{i, (s)}. 
$$
So we have to subtract Casimir factors 
$(4\pi i)^{-1}dd^\C \varphi_i | c''_i$. 
If we have at least two such factors  in the integral, 
we get zero after the integration thanks to 
Lemma \ref{6.16.08.1}.

\vskip 3mm

{\it Step 3 -- Using the projection formula}. 
Since we can have at most one factor like 
$dd^\C \varphi_{i} | c''_{i}$, we 
can take any non-special external vertex of the tree $T^0$, and use the projection formula 
for the integral related to this vertex. In our running example 
such a vertex is $x_1$, and we get
\be \la{6.16.08.3}
\stackrel{1}{\sim}\sum_{}
{\int}_{Y\times X}\Bigl(\gamma_0\wedge \gamma_1\wedge f_{1 \ast}G_X(y_1, x_2) 
\wedge f^*\gamma_2\wedge \alpha_3\Bigr) \langle a_3\cap {\cal H}_X \cap h\rangle 
\cdot
c_0 \otimes c_1 \otimes c_2.  
\ee
Here we may have 
$-(4\pi i)^{-1}dd^\C f_*\varphi_i|c''_i$ instead of 
$\gamma_i|c_i$, or $-(4\pi i)^{-1}dd^\C \varphi_i|c''_i$ instead of 
$f^*\gamma_i|c_i$. 
 Hence put $\stackrel{1}{\sim}$ instead of $=$  
in front of the formula.  Below we transform 
formula (\ref{6.16.08.3}), keeping  in mind the terms with $\varphi_i$'s. 

We perform this procedure for  all non-special external vertices of $T^0$.

\vskip 3mm
{\it Step 4 -- replacing $f_{1 \ast}G_X$ by 
$f^!_{2}G_Y$}. In our running example, 
using 
Lemma \ref{9.23.07.1}, we replace $f_{1 \ast}G_X(y_1, x_2)$ by 
$f^!_{2}G_Y(y_1, x_2)$ plus the terms in (\ref{6.16.08.4}). 
Observe the following: 

(a) The closed form does not change the Hodge correlator class -- see Step 7.

(b)The term 
$\sum_{k}p_Y^* \psi_{(k)}\wedge  p_X^*\alpha''_{(k)}$ gives
 zero after the integration, since
$$
\int_X \alpha''_{(k)} \wedge f^*\gamma'_2\wedge \alpha_3 = 
\int_Y f_*(\alpha''_{(k)}) \wedge \gamma'_2\wedge f_*(\alpha_3)=0. 
$$

(c) We may get a non-zero contribution of the term 
$\sum_{j}p_Y^* \gamma''^{\vee}_{(j)} \wedge p_X^*\varphi_{(j)}$. 

So we get 
\be \la{6.16.08.3a}
\stackrel{2}{\sim}\sum_{}
{\int}_{Y\times X}\Bigl(\gamma_0\wedge \gamma_1\wedge f^!_{2}G_Y(y_1, x_2) 
\wedge f^*\gamma_2\wedge \alpha_3\Bigr) \langle a_3
\cap {\cal H}_X \cap h\rangle \cdot
c_0 \otimes c_1 \otimes c_2.
\ee
Here $\stackrel{2}{\sim}\sum_{}$ means that we 
keep in mind the new terms with $\varphi_{(j)}$   
coming from (\ref{6.16.08.4}) -- see (c). 

In general we perform this procedure at  the edge $\stackrel{\to}{E_v}$ for 
every non-special external vertex $v$ of $T^0$. 
The term in (\ref{6.16.08.4a}) with a  harmonic form 
assigned to the other vertex then $v$  gives zero contribution. 
In the running example this is the case (b) above. Otherwise this follows from 
Lemma \ref{7.24.07.1}. The closed form (case (a)) also contributes
 zero  -- see Step 7.

\vskip 3mm
{\it Step 5 -- using the projection formula}. The current $f^!_{2}G_Y(y_1, x_2)$ 
is given by the form 
$f^*_{2}G_Y(y_1, x_2)$ at the generic point. 
We can use projection formula (\ref{PROJJ}). In the running example, using 
 it for the vertex $x_2$ we get 
\be \la{6.16.08.3b}
(\ref{6.16.08.3a}) =\sum_{}
{\int}_{Y\times Y}\Bigl(\gamma_0\wedge \gamma_1\wedge G_Y(y_1, y_2) 
\wedge \gamma_2\wedge f_{2\ast}\alpha_3\Bigr) \langle a_3\cap {\cal H}_X \cap h\rangle \cdot
c_0 \otimes c_1 \otimes c_2. 
\ee
Here one of the forms $\gamma_i$ 
can be replaced by the correction term 
$-(4\pi i)^{-1}f_*dd^\C\varphi_i$ from Step 2.

In general we do this for every external edge $E$ of $T^0$.

\vskip 3mm
{\it Step 6 --  Going down the tree}. 
We move down the tree towards the special edge: 

(a)  By using the  projection formula argument at a vertex $v$ (like in Step 5), and 

(b) By 
replacing $f_{1 \ast}G_X$ by 
$f^!_{2}G_Y$ at the corresponding edge $\stackrel{\to}{E_v}$ (like in Step 4). See Fig. \ref{hf20}, where 
we use the projection formula for the copy of $X$ assigned to
 the vertex $v_3$ as follows:  
$$
\int_{X_3} \omega\Bigl(f_2^!G_Y(y_1, x_3) 
\wedge f_2^!G_Y(y_2, x_3) \wedge G(x_3, x_4)\Bigr) =
\int_{Y_3} \omega\Bigl(G_Y(y_1, y_3) 
\wedge G_Y(y_2, y_3) \wedge f_{1\ast}G_X(y_3, x_4)\Bigr). 
$$
\begin{figure}[ht]
\centerline{\epsfbox{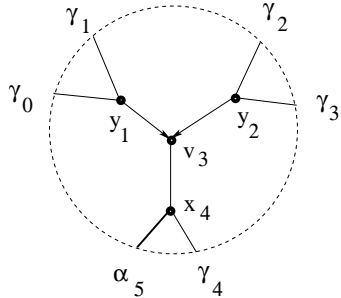}}
\caption{Using the projection formula at the vertex $v_3$.}
\label{hf20}
\end{figure}

Notice that if $\stackrel{\to}{E_v}$ is not an external edge of $T^0$, the latter procedure 
does not introduce  extra terms. 
This follows easily from Lemmas \ref{corT1} and \ref{7.24.07.1}. 

\vskip 3mm
{\it Step 7 --  independence of the closed form in (\ref{6.16.08.4})}. 
Adding to $f^!_2G_Y(y_1, x_2)$ a 
 closed form does not change the Hodge correlator class.
Indeed, one easily sees by using the projection formula that 
it results to adding a closed form 
 to a Green current $G_Y$. 
By Theorem \ref{3/11/07/101q} this does not 
change the Hodge correlator class.

\vskip 3mm

{\it Step 8.} The previous steps end up with an expression 
\be \la{CVA1}
\sim \sum{\rm Cor}_{{\cal H}, Y}(\gamma_0|c_0 
\otimes \ldots \otimes \gamma_{m-1}|c_{m-1} 
\otimes f_*\alpha_{m}) \langle a_m \cap {\cal H}_X \cap h\rangle . 
\ee 
Similarly to the proof of Lemma \ref{6.13.08.1}, we write it as 
\be \la{CVA2}
\stackrel{3}{\sim} \sum{\rm Cor}_{{\cal H}, Y}(\gamma_0|c_0 
\otimes \ldots \otimes \gamma_{m-1}|c_{m-1} 
\otimes \gamma_{m}) \langle f^!c_m \cap {\cal H}_X \cap h\rangle . 
\ee 
Here $\stackrel{3}{\sim}$ means equality modulo the extra terms provided by 
$f_*\alpha_{(k)}''^{\vee} \otimes a_{(k)}''^{\vee}$ and  (\ref{CFAsa}). 

\noindent
By (\ref{6.17.08.4}) we have
  $\langle f^!c_m \cap {\cal H}_X \cap h\rangle_X = 
\langle c_m \cap {\cal H}_X \cap f_*h\rangle_Y$. Thus 
the sum in (\ref{CVA2}) equals  $\delta_Yf_*(h)$. 

\vskip 3mm
It remains to take care of the extra terms, and show that their 
sum is a coboundary.  
We can get them from the following sources, marked $\stackrel{i}{\sim}$ in the proof: 

\vskip 3mm
(1) At most one term $-(4\pi i)^{-1}dd^\C f_*\varphi_i|c''_i$, 
see Step 2, for an external non-special edge of $T$. 

(2) Extra terms $\sum_{j}\gamma''^{\vee}_{(j)} \wedge d^\C f_*\varphi_{(j)}$ 
from Step 4(c) for non-special external vertices of $T^0$. 

(3) The extra term for the special external vertex of $T^0$ at Step (8), provided by 
$f_*\alpha_{(k)}''^{\vee} \otimes a_{(k)}''^{\vee}$.

\vskip 3mm
{\it Step 9 -- Calculating the extra terms}. 
By Lemma \ref{6.16.08.1}, 
terms with 
more then one factor $f_*dd^\C\eta_{i}$ are zero. 
Thus in the case (3) there are no extra terms (1).

\begin{figure}[ht]
\centerline{\epsfbox{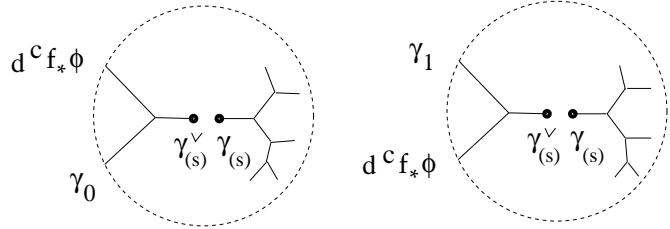}}
\caption{The sum of these contributions is zero.}
\label{hf14}
\end{figure}

(1). Moving $-(4\pi i)^{-1}d$ (or $-(4\pi i)^{-1}dd^\C$ 
if the tree $T$ has a single internal edge)
 from $-(4\pi i)^{-1}dd^\C f_*\varphi_i|c''_i$ 
to the operator $\omega$ applied to the Green currents of  the
internal edges of the tree $T$, and using the formula for $d\omega$ 
(respectively $dd^\C G$), 
we conclude that, by Lemma \ref{7.24.07.1},  only the Casimir terms for
 the external edges of $T^0$ can 
give non-zero contribution. Denote by $T_E''$ the small tree obtained 
by cutting the external edge $E$ of $T^0$. 
Then there are two cases:

(i) None of the edges of the tree $T_E''$ is  decorated by $f_*d^\C\varphi$ 
(respectively $f_*\varphi$). 
Then by Lemma \ref{CASCONT} the corresponding contribution is $\delta (\gamma_i|c_i'')$.

(ii) One of the edges of the tree $T_E''$ is  decorated by $f_*d^\C\varphi$ 
(respectively $f_*\varphi$), see Fig. \ref{hf14}. 
Then the contribution is zero. Indeed, there are two possible decorations,  
see 
Fig. \ref{hf14}. Their contribution is the simplest shuffle relation:
$$
{\rm Cor}_{\cal H}(\gamma_1|c_1\otimes f_*d^\C|c'' \otimes \gamma_{(s)}^{\vee}) 
+ {\rm Cor}_{\cal H}(f_*d^\C|c'' \otimes \gamma_1|c_1\otimes \gamma_{(s)}^{\vee}) = 0.
 $$

\vskip 3mm
(2) The contribution of a term (2) is 
$\delta (d^\C f_*\varphi_i|c''_i)$. 
 This is checked using Lemma \ref{CASCONT}. 

Taking into account the case (1)(i), we arrive at the term 
(\ref{7.12.08.3a}). 

\vskip 3mm
(3) Let us handle first the running example. The extra term 
is 
\be \la{6.16.08.3ba}
\sum_{k}\sum_{}
{\int}_{Y\times Y}\Bigl(\gamma_0\wedge \gamma_1\wedge G_Y(y_1, y_2) 
\wedge \gamma_2\wedge (4\pi i)^{-1}dd^\C \psi_{(k)}\Bigr) 
\langle a_{(k)}''^{\vee}\cap {\cal H}_X \cap h\rangle \cdot
c_0 \otimes c_1 \otimes c_2. 
\ee
Moving $(4\pi i)^{-1}dd^\C$ to $G_Y(y_1, y_2)$, and using the defining equation 
for $(4\pi i)^{-1}dd^\C G_Y$ we observe the following: 
The $\delta$-term vanishes after the summation over the (two) trivalent trees. 
The volume term vanishes. The Casimir term 
$\sum_s \gamma_{(s)}^{\vee}\otimes \gamma_{(s)} $ leads to an equality 
$$
(\ref{6.16.08.3ba}) = \sum_{s, k}\sum_{}
{\int}_{Y}\Bigl(\gamma_0\wedge \gamma_1\wedge \gamma_{(s)}^{\vee}\Bigr) 
\cdot {\int}_{Y}\Bigl(\gamma_{(s)} 
\wedge \gamma_2\wedge \psi_{(k)}\Bigr) 
\langle a_{(k)}''^{\vee}\cap {\cal H}_X \cap h\rangle \cdot
c_0 \otimes c_1 \otimes c_2 = 
$$
$$
\sum_{s, k}\sum_{\gamma_2|c_2}
{\rm Cor}_{\cal H}\Bigl((-1)^{|\gamma_{(s)}|}\gamma_{(s)}| 
\delta h_{\gamma_{(s)}} \otimes 
\gamma_2|c_2\otimes \psi_{(k)}\Bigr) 
\langle a_{(k)}''^{\vee}\cap {\cal H}_X \cap h\rangle.
$$
\begin{figure}[ht]
\centerline{\epsfbox{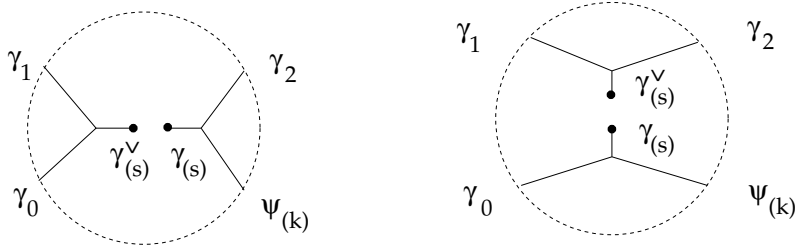}}
\caption{Correlator integrals for the extra term (3) in the running example.}
\label{hf13}
\end{figure}
It is illustrated by the diagram on the left on Fig \ref{hf13}. 
Adding a similar contribution of the right diagram on  Fig \ref{hf13}
we conclude that the total contribution is 
$$
\delta \Bigl(\sum_{s, k}\sum_{}
{\rm Cor}_{\cal H}\Bigl(
\gamma_0|c_0\otimes \gamma_1|c_1
\otimes \psi_{(k)}\Bigr) 
\langle a_{(k)}''^{\vee}\cap {\cal H}_X \cap h\rangle\Bigr).
$$

In general, if the tree $T^0$ has more then one edge, 
 the extra term (3) for integral (\ref{CVA1}) is 
$$
\delta \Bigl(\sum_{s, k}\sum_{}
{\rm Cor}_{\cal H}\Bigl(
\gamma_0|c_0\otimes \ldots \otimes \gamma_{m-1}|c_{m-1} 
\otimes d^\C\psi_{(k)}\Bigr) 
\langle a_{(k)}''^{\vee}\cap {\cal H}_X \cap h\rangle\Bigr).
$$
The difference with the case when the tree $T^0$ has one edge 
is that we have now $d^\C\psi_{(k)}$ instead of 
$\psi_{(k)}$. The proof is similar: 
we move the differential $(4\pi i)^{-1}d$ from $d^\C\psi_{(k)}$ to 
the operator $\omega$ applied to  the Green currents assigned to the 
internal edges of the tree $T$, and use the formula for $d\omega$. 
The Casimir term at certain internal edge $E$ of the tree $T$ 
is the only one which may contribute.  
If the edge $E$ is an 
internal edge of the tree $T^0$ its contribution is zero by 
Lemma \ref{7.24.07.1}. Otherwise the contribution is 
just as in Lemma \ref{CASCONT}.

Finally, for certain external edges $E$ of $T$ 
the factor $\gamma_i|c_i$ attached to $E$ may be replaced by the one of type (2).  
Then the Casimir term for the external edge of $T^0$ incident to $E$ 
provides zero contribution in the above argument. 
Thus we arrive  at the term (\ref{7.12.08.3}).

Theorem \ref{7.12.08.1} is proved.

 \label{HF4sec}


\end{document}